\documentclass[11pt]{amsart}
\usepackage[hmargin=3cm,vmargin=3cm]{geometry}

\usepackage{combelow} 

\usepackage[latin1,utf8]{inputenc} 
\usepackage[english]{babel}
\usepackage[backend=bibtex,doi=false,url=false,isbn=false,style=alphabetic,maxnames=6]{biblatex}
\usepackage{etex}
\AtBeginBibliography{\small}

\usepackage{xspace,amssymb,amsfonts,euscript}
\usepackage{amsthm,amsmath,amscd,stmaryrd,latexsym,youngtab,verbatim,mathrsfs}
\usepackage[nohug]{diagrams}

\usepackage{graphicx}
\usepackage[all]{xy}
\SelectTips{cm}{}

\usepackage{tikz, tikz-cd}
\usetikzlibrary{matrix, arrows, calc}
\tikzset{frontline/.style={preaction={draw=white,-,line width=6pt}},}  

\usepackage[dvips]{epsfig}
\usepackage{pinlabel}
\usepackage{psfrag}

\newcommand{\ig}[2]{\vcenter{\xy (0,0)*{\includegraphics[scale=#1]{fig/#2}} \endxy}}

\IfFileExists{comments.sty}{\usepackage{comments}}

\newcommand{\MH}[1]{\textcolor[rgb]{.6,.3,1}{MH: #1}}

\usepackage{enumitem}

\RequirePackage{color}
\definecolor{myred}{rgb}{0.75,0,0}
\definecolor{mygreen}{rgb}{0,0.5,0}
\definecolor{myblue}{rgb}{0,0.25,0.65}
\definecolor{references}{rgb}{0,0,1}

\newtheorem{thm}{Theorem}[section]
\newtheorem{lemma}[thm]{Lemma}
\newtheorem{theorem}[thm]{Theorem}

\newtheorem{proposition}[thm]{Proposition}

\newtheorem{corollary}[thm]{Corollary}

\newtheorem{conjecture}[thm]{Conjecture}

\newtheorem*{prop*}{Proposition}
\newtheorem*{lemma*}{Lemma}

\theoremstyle{definition}

\newtheorem{definition}[thm]{Definition}

\newtheorem{notation}[thm]{Notation}

\newtheorem{example}[thm]{Example}

\theoremstyle{remark}
\newtheorem{remark}[thm]{Remark}

\numberwithin{equation}{section}

    \def\SM{{\mathbb{S}}}

    \def\AC{{\mathcal{A}}}
    
\def\CB{{\mathbf C}}  
\def\CC{{\mathcal{C}}}
\def\DB{{\mathbf D}}    \def\DC{{\mathcal{D}}}
    
    \def\FC{{\mathcal{F}}}

    \def\JC{{\mathcal{J}}}
\def\KB{{\mathbf K}}    \def\KC{{\mathcal{K}}} 
\def\LB{{\mathbf L}}    \def\LC{{\mathcal{L}}}
    \def\MC{{\mathcal{M}}}
    \def\NC{{\mathcal{N}}}
    \def\OC{{\mathcal{O}}}
\def\PB{{\mathbf P}}    \def\PC{{\mathcal{P}}}
\def\QB{{\mathbf Q}}

    \def\YC{{\mathcal{Y}}}

\def\a{\alpha}
\def\b{\beta}

\def\d{\delta}

\def\e{\varepsilon}
\renewcommand{\k}{\mathbbm{k}}

\def\L{\Lambda}

\let\phi=\varphi

\usepackage{bbm}
\def\C{{\mathbbm C}}

\def\R{{\mathbbm R}}
\def\Z{{\mathbbm Z}}
\def\Q{{\mathbbm Q}}
\def\1{\mathbbm{1}}
\newcommand{\one}{\1}

\newcommand{\ptrz}{\operatorname{pTr}_0}
\newcommand{\Sym}{\operatorname{Sym}}
\newcommand{\Hnorm}{\HKR^{\text{norm}}}
\newcommand{\Pnorm}{\PC^{\text{norm}}}

\newcommand{\Lab}{\LB}
\newcommand{\cyc}{\operatorname{Cyc}}

\newcommand{\seq}[1]{\operatorname{Seq}(#1)}

\newcommand{\FY}{\operatorname{FY}}
\newcommand{\DGR}{\operatorname{DGR}}
\newcommand{\diff}{\Delta_{+}}

\newcommand{\yring}[1]{\Q[\yy]_{#1}}
\newcommand{\xyring}[1]{\Q[\xx,\yy]_{#1}}
\newcommand{\IY}{\operatorname{Ideal}}

\newcommand{\HY}{\operatorname{HY}}
\newcommand{\CY}{\operatorname{CY}}
\newcommand{\components}{\pi_0}
\renewcommand{\split}{\operatorname{split}}
\newcommand{\bbeta}{\boldsymbol{\beta}}
\newcommand{\ttheta}{\boldsymbol{\theta}}
\newcommand{\Fac}{\mathbf{Fac}}

\newcommand{\derived}{\DC}
\newcommand{\ii}{\mathbf{i}}
\newcommand{\jj}{\mathbf{j}}
\newcommand{\kk}{\mathbf{k}}
\newcommand{\lseq}{\mathbf{l}}
\newcommand{\Homc}{\underline{\Hom}}
\newcommand{\Endc}{\underline{\End}}
\newcommand{\xx}{\mathbf{x}}
\newcommand{\yy}{\mathbf{y}}

\newcommand{\w}{\omega}

\renewcommand{\mod}{\textrm{-mod}}

\newcommand{\ft}{\operatorname{FT}}

\newcommand{\HHH}{\operatorname{HHH}}
\newcommand{\HKR}{H_{KR}}
\newcommand{\CKR}{\operatorname{CKR}}
\newcommand{\HH}{\operatorname{HH}}

\renewcommand{\setminus}{\smallsetminus}

\newcommand{\Endg}{\End^\Z}

\newcommand{\smMatrix}[1]{\left[\begin{smallmatrix}#1\end{smallmatrix}\right]}
\newcommand{\Proj}{\operatorname{Proj}}
\newcommand{\sgn}{\operatorname{sgn}}
\newcommand{\ev}{\operatorname{ev}}
\newcommand{\thth}{\boldsymbol{\theta}}

\newcommand{\sqmatrix}[1]{\left[\begin{matrix} #1\end{matrix}\right]}

\renewcommand{\to}{\rightarrow}

\renewcommand{\sl}{\mathfrak{sl}}

\def\gmod{\text{-gmod}}
\def\gbimod{{\text{-gbimod}}}

\newcommand{\refequal}[1]{\xy {\ar@{=}^{#1}
(-1,0)*{};(1,0)*{}};
\endxy}

\newcommand{\Hom}{\operatorname{Hom}}

\newcommand{\End}{\operatorname{End}}

\newcommand{\Ext}{\operatorname{Ext}}

\newcommand{\Id}{\operatorname{Id}}

\newcommand{\inv}{^{-1}}

\newcommand{\Cone}{\operatorname{Cone}}
\newcommand{\rev}{\operatorname{rev}}

\newcommand{\Bim}{{\rm Bim }}

\newcommand{\SBim}{\SM\Bim}

\newcommand{\wt}{\operatorname{wt}}

\newcommand{\FT}{\operatorname{FT}}

\newcommand{\Br}{\operatorname{Br}}

\DeclareMathOperator{\dinv}{dinv}
\DeclareMathOperator{\Hilb}{Hilb}
\DeclareMathOperator{\Alt}{Alt}

\newcommand{\AJ}{\mathcal{J}}
\newcommand{\Ht}{\widetilde{H}}

\begin{document}

\begin{abstract} We define a deformation of the triply graded Khovanov-Rozansky homology of a link $L$ depending on a choice of parameters $y_c$ for each component of $L$, which satisfies link-splitting properties similar to the Batson-Seed invariant.  Keeping the $y_c$ as formal variables yields a link homology valued in triply graded modules over $\Q[x_c,y_c]_{c\in \components(L)}$.  We conjecture that this invariant restores the missing $Q\leftrightarrow TQ\inv$ symmetry of the triply graded Khovanov-Rozansky homology, and in addition satisfies a number of predictions coming from a conjectural connection with Hilbert schemes of points in the plane. We compute this invariant for all positive powers of the full twist and match it to the family of ideals appearing in Haiman's description of the isospectral Hilbert scheme.
\end{abstract}

\title{Hilbert schemes and $y$-ification of Khovanov-Rozansky homology}

\author{Eugene Gorsky}
\address{University of California, Davis}
\address{International Laboratory of Representation Theory and Mathematical Physics, Higher School of Economics, Moscow, Russia}
\email{egorskiy@math.ucdavis.edu}

\author{Matthew Hogancamp} \address{University of Southern California, Los Angeles}
\email{hogancam@usc.edu}

\maketitle

\setcounter{tocdepth}{1}
\tableofcontents

\section{Introduction}
Recent conjectures of the first author, Negu\cb{t} and Rasmussen \cite{GNR}, and of Oblomkov-Rozansky \cite{OR1,OR2,OR3} relate the Khovanov-Rozansky homology of links to the algebraic geometry of the Hilbert scheme of points on the plane.  To each $n$-strand braid $\b$ these conjectures associate a complex of coherent sheaves $\FC_\b$ on the Hilbert scheme $\Hilb^n(\C^2)$ whose cohomology recovers the Khovanov-Rozansky homology of the link $L=\hat\b$ obtained by closing $\b$. 

In this paper, we make a first step towards the proof of this conjectures. Following the ideas in \cite{GNR}, we construct a complex of {\em quasicoherent} sheaves on the Hilbert scheme for each braid, and compute it explicitly for torus links $T(n,kn)$ for all $n,k>0$. We expect them to be coherent, but do not prove it here. We postpone the detailed discussion of this construction until \S \ref{subsec:introHilb}, but first point out a major obstacle that one needs to overcome. 

First, note that if $L\subset \R^3$ is an oriented link then Khovanov-Rozansky homology of $L$ is a triply graded module over $\Q[x_c]_{c\in \components(L)}$, where $\pi_0(L)$ denotes the set of components of $L$.  If $L$ is presented as the closure of an $n$-strand braid $\b$ then this action is inherited from an action of $\Q[x_1,\ldots,x_n]$.  To match this structure on the Hilbert scheme side, we remark that the sheaves predicted by the conjectures in \cite{GNR,OR1,OR2,OR3} are supported on the scheme $\Hilb^n(\C^2,\C)$ of points on the $x$-axis $\C\cong \{y=0\}$.  Then the Hilbert-Chow map
$\Hilb^n(\C^2,\C)\to \Sym^n(\C\times \{0\})$ endows the cohomology of every sheaf on $\Hilb^n(\C^2,\C)$ with the structure of a module over symmetric functions in $n$ variables $x_1,\ldots,x_n$.  It is expected that the cohomologies agree not only as vector spaces, but as modules over polynomial rings (see \S \ref{subsec:introHilb} for more).

From the viewpoint of algebraic geometry, many sheaves on $\Hilb^n(\C^2,\C)$ naturally arise by restricting sheaves on the full Hilbert scheme $\Hilb^n(\C^2)$.  The following question naturally arises: is there a deformation of Khovanov-Rozansky link homology, taking values in triply graded modules over $\Q[x_c,y_c]_{c\in \components(L)}$, which corresponds to sheaves on the full Hilbert scheme $\Hilb^n(\C^2)$?  The purpose of this paper is to construct precisely such a link homology theory and show that it indeed has a strong relationship with Hilbert schemes.


Let $L\subset \R^3$ be an oriented link, presented as the closure of an $n$-strand braid $\b$.  Let $\CKR(\b)$ denote the complex which computes the Khovanov-Rozansky homology of $L=\hat\b$, as constructed for instance in \cite{Kh07}.  In \S \ref{sec:homology} (see also \S \ref{subsec:introdescription} of the introduction) we construct an explicit complex of the form
\[
\CY(\b)=\CKR(\b)\otimes_\Q \Q[y_c]_{c\in \pi_0(L)}, \qquad d_{\CY(\b)} = d_{\CKR}\otimes 1 + \sum_c h_c\otimes y_c,
\]
and prove the following.
\begin{theorem}
\label{th:intro invariance}
The homology $\HY(\b)$ of $\CY(\b)$ depends only on the link $L=\hat\b$ up to isomorphism and overall shift of $\Q[x_c,y_c]_{c\in \pi_0(L)}$-modules.
\end{theorem}
We refer to $\HY(\b)$ as $y$-ified Khovanov-Rozansky homology; it is triply graded and if $\b$, $\b'$ represent the same link then $\HY(\b)\cong \HY(\b')$ up to overall shift which can be fixed by an appropriate normalization (\S \ref{subsec:conventions}). We will henceforth write $\HY(L)$ for $\HY(\b)$.  Setting the $y_c=0$ before taking homology recovers Khovanov-Rozansky homology, and specializing $y_c$ to some other scalars gives a family of homology theories $\HY(L,\Q_\nu)$, parametrized by points $\nu\in \Q^{\pi_0(L)}$, with link splitting properties similar to the Batson-Seed invariant (\S \ref{subsec:introSplitting} and \S \ref{subsec:introcomparison}).

We are able to compute the $y$-ified homology for several families of links and compare with predictions coming from Hilbert schemes.  Our next main theorem is the first result of its kind, giving a direct link between Soergel bimodules and Hilbert schemes.

\begin{definition}
Let $\xx=(x_1,\ldots,x_n)$, $\yy=(y_1,\ldots,y_n)$, $\ttheta=(\theta_1,\ldots,\theta_n)$ denote sets of formal variables, with tridegrees $\deg(x_i)=(2,0,0)$, $\deg(y_i)=(-2,0,2)$, and $\deg(\theta_i)=(2,-1,0)$.  The $\theta_i$ are regarded as odd variables, and we let $\Q[\xx,\yy,\ttheta]$ denote the super polynomial ring $\Q[\xx,\yy]\otimes_\Q\Lambda[\ttheta]$.

Let $\AJ_n\subset \Q[\xx,\yy,\ttheta]$ denote the ideal generated by the anti-symmetric polynomials with respect to the $S_n$-action which simultaneously permutes all three sets of variables.
\end{definition}

\begin{remark}
The three gradings are referred to as the Soergel bimodule grading $\deg_Q$, the Hochschild degree $\deg_A$, and the homological degree $\deg_T$.
\end{remark}
\begin{remark}
The super-polynomial ring $\Q[\xx,\yy,\ttheta]$ is the $y$-ified homology of the $n$-component unlink, and the $y$-ified homology of every $n$-component link is naturally a triply graded $\Q[\xx,\yy,\ttheta]$-module.
\end{remark}

\begin{theorem}
\label{th: intro ft}
The $y$-ified homology of the torus link $T(n,nk)$ is isomorphic to the ideal $\AJ_n^k\subset \Q[\xx,\yy,\ttheta]$ as triply graded $\Q[S_n]\ltimes \Q[\xx,\yy,\ttheta]$-modules, for all $n,k\geq 0$.  These isomorphisms are compatible with the multiplication maps
\[
\HY(T(n,nk))\otimes \HY(T(n,nl))\rightarrow \HY(T(n,nk+nl)).
\]
\end{theorem}
Let us explain the statement a bit more.  The $T(n,nk)$ torus link is the closure of $\ft_n^k$, where $\ft_n=(\sigma_1\cdots\sigma_{n-1})^n$ denotes the full twist braid, expressed here in terms of the usual braid generators $\sigma_i$.  The full-twist is central in group, hence the homology $\HY(\ft_n^k)$ admits an action of the braid group, given by conjugation (\S \ref{subsec: sn action}).  The above asserts that this action factors through the symmetric group.  The algebra structure on $\bigoplus_{k\geq 0}\HY(\ft_n^k)$ is defined via a general construction which we recall in \S \ref{subsec:localpicture}.

We have the following important corollary for the usual Khovanov-Rozansky homology.

\begin{corollary}\label{cor:KRFT}
The Khovanov-Rozansky homology of $T(n,nk)$ is isomorphic to the quotient $\AJ_n^k/(\yy)\AJ_n^k$, as triply graded $\Q[S_n]\ltimes \Q[\xx,\ttheta]$-modules, for all $n,k\geq 0$.  These isomorphisms are compatible with the multiplication maps
\[
\HKR(T(n,nk))\otimes \HKR(T(n,nl))\rightarrow \HKR(T(n,nk+nl)).
\]
\end{corollary}

\begin{remark}
Earlier work of the second author and Ben Elias \cite{EH,Hog17b} computes the homology $\HKR(T(n,nk))$ as a triply graded vector space, for all $n,k\geq 0$.  In particular these computations imply that the homology is supported only in even homological degrees.  This fact is crucial in our proof of Theorem \ref{th: intro ft}.
\end{remark}

\begin{remark}\label{rmk:introJn}
Let $J_n$ denote the $\deg_A=0$ component of $\AJ_n$.  In other words $J_n\subset \Q[\xx,\yy]$ is the ideal spanned by the antisymmetric polynomials.  Our proof of Theorem \ref{th: intro ft} relies on some deep results of Haiman, particularly that $J_n^k$ is free as a $\Q[\yy]$-module and equals the intersection
\[
J_n^k=\bigcap_{i\neq j}(x_i-x_j,y_i-y_j)^k\subset \Q[\xx,\yy].
\]
These ideals are intimately related to Hilbert schemes by work of Haiman \cite{haiman2001hilbert}.  See \S \ref{subsec:introHilb}.
\end{remark}

Let us also mention another major motivation for this work.  The Poincar\'e series $\PC_L(Q,A,T)$ of Khovanov-Rozansky homology $\HKR(L)$ is conjecturally unchanged by the swapping $Q\leftrightarrow TQ\inv$ and fixing $AQ^{-2}$.  This symmetry would categorify the well known $Q\leftrightarrow -Q\inv$ symmetry of the HOMFLY-PT polynomial.  This symmetry cannot be realized on the level of homology, since the $\Q[x_c]$-module structure on $\HKR(L)$ breaks the symmetry between $Q^2$ and $T^2Q^{-2}$.  We originally sought to construct the $y$-ified homology in order to restore this missing symmetry.  We do not prove any results in this direction, but all available computations suggest the following.

\begin{conjecture}\label{conj:introsymmetry}
The triply graded homology $\HY(L)$ is unchanged up to isomorphism by the regrading according to
\[
i(1,0,0)+j(-2,1,0)+k(-1,0,1)\ \ \leftrightarrow \ \  k(1,0,0)+j(-2,1,0)+i(-1,0,1).
\]
Further, this symmetry exchanges the $x_i$ and $y_i$ actions, and is equivariant with respect to the $\theta_i$ action for all $i$.
\end{conjecture}

\subsection{Link splitting properties}
\label{subsec:introSplitting}

Let $L=\hat\b$ be an oriented link, presented as the closure of $\b\in \Br_n$, with labelled components $L = L_1\cup \cdots \cup L_r$.  For each $\nu\in \Q^{r}$, let $\HY(\b,\Q_\nu)$ be the homology of $\Q\otimes_{\Q[y_1,\ldots,y_r]} \CY(\b)$, in which we have specialized the $y_i$ to scalars $y_i=\nu_i\in \Q$.  We also write $\HY(L;\Q_\nu)=\HY(\b,\Q_\nu)$.  Then
\[
\HKR(L)\cong \HY(L;\Q_0).
\]
Thus, we may think of $\HY(L)$ as a family of link homology theories parametrized by $y_1,\ldots, y_r$.  For $\nu\neq 0$, specialization of the $y_i$ variables collapses the trigrading on $\HY(L,\Q_\nu)$ to a bigrading, via $(\deg_A,\deg_Q+\deg_T)$, and the $\deg_Q$ grading on $\CY(L)$ induces a filtration on $\HY(L,\Q_\nu)$ (see \S \ref{subsec:homology with coefficients} for details).

The existence of $\CY(L)$ implies the existence of certain interesting homological operations in $\HKR(L)$. In particular, we prove the following.

\begin{proposition}
\label{prop:intro monodromy}
Let $d=d_0+\sum_i h_iy_i$ denote the differential on $\CY(L)$.  Then $h_i$ yields a well-defined operator on $\HKR(L)$ of degree $(2,0,-1)$ which corresponds to ``monodromy'' of the operator $x_c$ around the component $c$. 
\end{proposition}

Similar such monodromy maps appear in other contexts in the work of Batson-Seed \cite{Batson}, Sarkar \cite{Sarkar}, and others \cite{Zemke,BLS17}.  The monodromy maps arise as the differentials in the $E_2$ page of a spectral sequence abutting to $\HY(L;\Q_\nu)$ for $\nu\neq 0$.

\begin{theorem}\label{thm:introSS}
For each $r$-component link $L$ and each $\nu\in \Q^r$ there exists a spectral sequence with $E_2$ page $\HKR(L)$ and $E_\infty$ page the associated graded of $\HY(L;\Q_\nu)$ with respect to the $\deg_Q$ filtration.
\end{theorem}

The Batson-Seed link-splitting phenomenon has a counterpart in our setting as well.

\begin{theorem}
\label{th:intro splitting}
Let $L$ be written in terms of its components $L=L_1\cup\cdots \cup L_r$.  If $\nu_i\neq \nu_1$ for all $i\neq 1$, then
\begin{equation}\label{eq:splitComponent}
\HY(L;\Q_{\nu})\cong \HY(L_1)\otimes_\Q \HY(L_2\cup\cdots\cup L_n)
\end{equation}
as bigraded vector spaces.
\end{theorem}
We remark that the right-hand side of \eqref{eq:splitComponent} is the $y$-ified homology of link in which the component $L_1$ has been ``unlinked'' from the remaining components.  In particular if all $\nu_i$ are distinct, then $\HY(L,\Q_\nu)$ is isomorphic to $\HY(\split(L),\Q_\nu)$, where $\split(L)$ is the split union of the components of $L$.  As a consequence, we obtain the following corollary which is a direct generalization of the Batson-Seed spectral sequence in Khovanov homology \cite{Batson}.

\begin{corollary} 
\label{cor:intro ss}
There is a spectral sequence with $E_2$ page isomorphic to $\HKR(L)$ and $E_{\infty}$ page isomorphic to (the $\deg_Q$-associated graded of) $\HKR(\split(L))$.
This spectral sequence preserves the $(\deg_A,\deg_Q+\deg_T)$ bigrading.
\end{corollary}

\label{subsec:introsplitmaps}
As mentioned above, we have $\HY(L;\Q_\nu)\cong \HY(\split(L);\Q_\nu)$ for almost all $\nu$.  In fact more is true: there exists a map $\CY(L)\rightarrow \CY(\split(L))$ at the chain level which becomes an equivalence after specialization to $\nu$.

\begin{theorem}
\label{th: intro splitting injective}
Suppose that a link $L$ can be transformed to a link $L'$ by a sequence of  crossing changes  between different components. Then there is a homogeneous ``link splitting map'' $\Psi: \HY(L)\to \HY(L')$ which preserves the $\Q[\xx,\yy,\boldsymbol{\theta}]$--module structure.  If, in addition, $\HY(L)$ is free as a $\Q[\yy]$-module then $\Psi$ is injective.  If the crossing changes only involve components $i$ and $j$, then the link splitting map becomes a homotopy equivalence after inverting $y_i-y_j$, where $i$ and $j$ are the components involved.
\end{theorem}
\begin{remark}
The map $\Psi$ depends on the sequence of crossing changes; for instance switching a crossing then switching back corresponds to multiplication by $y_i-y_j$.  In general, switching a crossing from positive to negative has degree $(0,0,0)$, while switching from negative to positive has degree $(-2,0,2)$.
\end{remark}

The condition that $\HY(L)$ is free over $\Q[\yy]$ is surprisingly common.

\begin{definition}
A link will be called \emph{parity} if $\HKR(L)$ is supported only if even homological degrees (or odd homological degrees).
\end{definition}
For example, recent computations of Ben Elias, the second author, and Anton Mellit \cite{EH,Hog17b,Mellit} show that all positive torus links are parity.

\begin{theorem}
\label{th: intro parity is flat}
If an $r$-component link $L$ is parity then $\HY(L)\cong \HKR(L)\otimes_\Q \Q[\yy]$ is a free $\Q[\yy]$-module.  Consequently any link splitting map identifies $\HY(L)$ with a $\Q[\xx,\yy,\ttheta]$-submodule of $\HY(\split(L))$.
\end{theorem}
\begin{remark}
The above implies that if an $r$-component link $L$ is parity then the Poincar\'e series of $\HY(L)$ and $\HKR(L)$ agree up to a factor $(1-t)^r$, where $t=T^2Q^{-2}$.
\end{remark}

As a special case, we see immediately that $\HY(T(n,nk))$ is isomorphic to an ideal in $\Q[\xx,\yy,\ttheta]$, since it is parity and the $\split(L)$ is the $n$-component unlink in this case.  More generally, we have the following.

\begin{corollary}
\label{th: intro splitting ideal}
Suppose that $\beta$ is a pure braid on $n$ strands such that 
the braid closure $L:=\hat\b$ is parity, and let $L'$ denote the $n$-component unlink. Then the map $\Psi:\HY(L)\to \HY(L')=\Q[\xx,\yy,\boldsymbol{\theta}]$ is injective, so $\HY(L)$ is isomorphic to a certain ideal in $\Q[\xx,\yy,\boldsymbol{\theta}]$. 
\end{corollary}

This corollary allows one to give a more concrete description of $\HY(L)$ (and, as a consequence, of $\HKR(L)$) for many interesting braids.  In addition to Theorem \ref{th: intro ft} on the homologies of $\ft_n^k$, we also prove the following.

\begin{proposition}
\label{prop: intro JM}
Let $\LC_n=\sigma_{n-1}\cdots \sigma_2\sigma_1^2\sigma_2\cdots \sigma_{n-1}$ denote the Jucys-Murphy braid. Then the link splitting morphism $\Psi$ identifies $\HY(\LC_n)$ with an intersection of ideals:
\[
\HY(\LC_n)\cong \bigcap_{i=1}^{n-1}(x_i-x_n,y_i-y_n,\theta_i-\theta_n)\subset \Q[\xx,\yy,\ttheta].
\]
\end{proposition}

\subsection{Description of the homology}
\label{subsec:introdescription}
In this section we give an elementary description of our homology theory.  Let $L\subset \R^3$ be an oriented link, and choose a braid representative $\b\in \Br_n$ for $L$.  Let $R=\Q[x_1,\ldots,x_n]$, and recall that a Soergel bimodule $B\in \SBim_n$ is, in particular, a graded $(R,R)$-bimodule.  Let $F(\b)\in \KC^b(\SBim_n)$ be the Rouquier complex (\S \ref{sub:soergel rouquier}) associated to $\b$, and let $\CKR(\b):=\HH(F(\b))$ be the complex obtained by applying the Hochschild cohomology functor to the terms of $F(\b)$.  Then the homology of $\HH(F(\b))$ is isomorphic to the Khovanov-Rozansky homology of $L$ by \cite{Kh07}.  Next we discuss in these terms how to construct the $y$-ified homology. 


Let $D$ be the link diagram determined by $\b$, and let $V\subset D$ be the set of crossings, i.e.~4-valent vertices.  For each point $p\in D$ away from the crossings, we have an endomorphism $x_p: \HH(F(\b))\rightarrow \HH(F(\b))(2)$ corresponding to ``multiplication by $x$ at $p$.''  To be more precise, let us factor $\b$ as $\b\cong \b'\cdot 1_n\cdot \b''$, where $1_n$ denotes the trivial braid, and $p$ lies on the $i$-th strand of $1_n$.  Then $x_p$ corresponds to $\Id_{F(\b')}\otimes x_i\otimes \Id_{F(\b'')}$ with respect to the isomorphism
\[
F(\b)\cong F(\b')\otimes_R R\otimes_R F(\b''),
\] 
where $R=F(1_n)$ is the trivial Soergel bimodule.  If $p,q$ are points on $D$ which are on the same component of $D\setminus V$, then $x_p=x_q$.  If $p$ and $q$ are on the same link component, then $x_p\simeq x_q$.  For each pair $p,q$ of points which are separated by a single crossing, we let $h_{p,q}\in \End^{2,-1}(F(\b))$ denote a chosen homotopy, satisfying
\[
d\circ h_{p,q} + h_{p,q}\circ d = x_p-x_q.
\]
Let $\components(L)$ be the set of components of $L$, and choose $c\in \components(L)$.  Choose a sequence $p_1,\ldots,p_r$ of points on $c$, away from the crossings, such that $p_i$ and $p_{i-1}$ are separated by exactly one crossing for all $i$ (indices read modulo $r$), and set
\begin{equation}
\label{def hc}
h_c:=\sum_{i=1}^r h_{p_i,p_{i-1}}.
\end{equation}
Then $[d,h_c]$ is a finite telescoping sum, hence zero.  The homotopies $h_{p,q}$ can all be chosen so that $h_{p,q}^2=0$, and distinct homotopies anti-commute with one another.  Thus, the following defines the chain complex
\[
\CY(\b):=\HH(F(\b))\otimes_\Q \Q[y_c]_{c\in \components(L)},\qquad\qquad d_{\CY(\b)}=d\otimes 1 + \sum_{c\in \components(L)} h_c\otimes y_c
\]
as discussed earlier in the introduction.

Now we discuss gradings.  We want to put a $\Z\times \Z\times \Z$ grading on $\CY$ so that the differential is homogeneous of tridegree $(0,0,1)$.  Observe that the Hochschild cohomology $\HH(B)=\Ext_{R\otimes_\Q R}(R,B)$ of a Soergel bimodule is a bigraded vector space.  If $C$ is a complex of Soergel bimodules, then applying $\HH$ term-by-term gives a complex
\[
\HH(C) = \begin{diagram} \cdots & \rTo^{\HH(d_C)} & \HH(C^k) & \rTo^{\HH(d_C)} &  \HH(C^{k+1})& \rTo^{\HH(d_C)} &  \cdots\end{diagram}.
\]
We regard $\HH(C)$ as being a triply graded chain complex, where the term in tridegree $(i,j,k)$ is $\HH^{i,j}(C^k)$.  An element $c\in \HH^{i,j}(C^k)$ will be said to have \emph{bimodule degree} (or \emph{quantum degree}) $\deg_Q(c)=i$, \emph{Hochschild degree} $\deg_A(c)=j$, and \emph{homological degree} $\deg_T(c)=k$. 

The endomorphism $x_p$ has degree $(2,0,0)$, hence the homotopies $h_{p,q}$ have degree $(2,0,-1)$ when regarded as homogeneous linear endomorphisms of $\HH(C)$. In order for each summand $h_c\otimes y_c$ of the differential of $\CY(F(\b))$ to be homogeneous of degree $(0,0,1)$ we must declare that the $y_i$ have degree $(-2,0,2)$.  Extending multiplicatively, we regard $\Q[y_c]_{c\in \components(L)}$ as a triply graded algebra.  This yields the desired grading on $\CY(\b)$.

\subsection{A local picture for our homology}
\label{subsec:localpicture}
We now give a more category-theoretic picture for our homology, which is necessary for defining our homology at the level of braids, much like Khovanov-Rozansky homology.  In \S \ref{subsec:yifications} we define a category $\YC(\SBim_n)$ of \emph{$y$-ified complexes of Soergel bimodules}.  This category is monoidal and triangulated.  An object of $\YC(\SBim_n)$ is a triple $(C,w,\Delta)$ such that $C\in \KC^b(\SBim_n)$ is a finite complex of Soergel bimodules, $w\in S_n$ is a permutation, and $\Delta$ is a degree $(0,0,1)$ $R\otimes_\Q R\otimes_\Q \Q[\yy]$-linear endomorphism of $C\otimes_\Q \Q[\yy]$ such that
\[
\Delta^2 = \sum_i (x_{w(i)}\otimes 1-1\otimes x_i)\otimes y_i.
\]
Here the action of $x_i\otimes 1$ and $1\otimes x_i$ uses the bimodule structure on $C$.  Since $\Delta^2$ is nonzero, a $y$-ification is a particular kind of matrix factorization or curved complex.

\begin{remark}
The category $\YC(\SBim_n)$ splits into blocks according to the permutation $w$, and there are no nonzero morphisms relating objects in different blocks. 
\end{remark}

There is a triangulated, monoidal forgetful functor $\YC(\SBim_n)\rightarrow \KC^b(\SBim_n)$ sending $(C,w,\Delta)\mapsto C$, and we say that $(C,w,\Delta)$ is a \emph{$y$-ification of $C$ }.  In \S \ref{sec:setup} we prove that Rouquier complexes $F(\b)$ admit unique $y$-ifications $\FY(\b)$ up to canonical equivalence.  Then, in \S \ref{sec:homology} we define the homology $\HY$ of a $y$-ification and show that $\HY(\b):=\HY(\FY(\b))$ is an invariant of the braid closure (Theorem \ref{th:intro invariance}).

It is sometimes useful to view $\HY$ as a representable functor.  To do this, we first embed $\SBim_n$ into the bounded derived category $\DC_n:=D^b(R\gbimod)$ of graded $R$-bimodules.  The category $\DC_n$ is monoidal, and carries grading shift functors of the form $(i,j)$, where $(1,0)$ is the shift of graded bimodules, and $(0,1)$ is the homological shift.  Let $\SBim_n^{\Ext}\subset \DC_n$ denote the full subcategory on the objects $B(i,j)$ where $B\in \SBim_n$.  Then the usual Khovanov-Rozansky homology can be described succinctly as follows: the Rouquier complex can be viewed as living in $\KC^b(\SBim_n^{\Ext})$.  This category carries grading shift functors of the form $(i,j)[k]$ where $(i,j)$ is inherited from $\DC_n$ and $[k]$ is the homological shift in $\KC^b(\cdots)$ (not to be confused with the homological shift in $\DC_n$!).  Then
\[
\HKR(\hat\b) \cong \Hom^{\Z\times \Z\times \Z}_{\KC(\SBim_n^{\Ext})} (\one, F(\b)) :=\bigoplus_{i,j,k\in \Z}\Hom_{\KC(\SBim_n^{\Ext})} (\one, F(\b)(i,j)[k]).
\]

For the $y$-ified homology we can say something similar.  The $y$-ified Rouquier complexes can be thought of as living in the triangulated monoidal category $\YC(\SBim_n^{\Ext})$ of $y$-ifications in $\SBim_n^{\Ext}$.  This category again carries grading shift functors of the form $(i,j)[k]$.  If $\b$ is a pure braid then
\[
\HY(\hat\b)\cong \Hom^{\Z\times \Z\times \Z}_{\YC(\SBim_n^{\Ext})} (\one, \FY(\b)).
\]
If $\b$ is not a pure braid then one must replace $\one$ by an appropriate representing object $\one_w$, where $w\in S_n$ is the permutation represented by $\b$ (\S \ref{subsec:modules}).

With this description, we can now understand the algebra structure on $\bigoplus_{k\geq 0}\HY(\ft_n^k)$ as inherited from the tensor product of morphisms in $\one\rightarrow \FY(\ft_n^k)$.

\subsection{Connection to Hilbert schemes of points in the plane}
\label{subsec:introHilb}
Next, we would like to explain the relation between Theorem \ref{th: intro ft} and the conjectures of \cite{GNR,OR1,OR2,OR3}. Let $\Hilb^n(\C^2)$ denote the Hilbert scheme of $n$ points on the plane.  It has a natural projection (Hilbert-Chow map) to the symmetric product $\Sym^n\C^2$.  The {\em isospectral Hilbert scheme} $X_n$ is defined \cite{haiman2001hilbert} as the reduced fibered product of $\Hilb^n(\C^2)$ and $(\C^2)^n$ over $\Sym^n\C^2$.  Let $J_n\subset \Q[\xx,\yy]$ denote the ideal generated by the anti-symmetric polynomials with respect to the diagonal $S_n$ action. Define the graded algebra:
\[
\AC=\bigoplus_{k\ge 0}J_{n}^k.
\]
Theorem of Haiman (see \S \ref{sec:Hilbert}) states that $\Proj (\AC\otimes_{\Q}\C)=X_n.$

According to Theorem \ref{th: intro ft}, the $\deg_A=0$ component $\HY^0(\ft_n^k)\subset \HY(\ft_n^k)$ is isomorphic to $J_n^k$.  Given a braid $\beta$, we can define 
\[
M_{\beta}:=\bigoplus_{k\ge 0}\HY^0(\beta \FT_n^k)\otimes\Q \C
\]
This is clearly a graded $\AC$-module, so $M_{\beta}$ defines a quasicoherent sheaf $\FC_{\beta}$ on $X_n.$
We expect the cohomology of $\FC_{\beta}$ to be related to the $y$-ified homology of $\beta$.  In particular, we prove the following
\begin{theorem}
\label{th: intro Hilbert}
For all $k>0$ the $y$-ified homology of $\FT_n^k$ (with complex coefficients) is isomorphic to the space of global sections of the vector bundle $\Lambda(T^*)\otimes \OC(k)$ on the isospectral Hilbert scheme $X_n$, where $T$ is the tautological rank $n$ bundle on $X_n$.
\end{theorem}

Similarly, following \cite{GNR,OR1} we define the space $X_{n}(\C^2,\C)$, the isospectral Hilbert scheme supported on $\{y=0\}$  as the preimage of the subspace $\{y_1=\ldots=y_n=0\}$ under the natural projection $X_n\to (\C^2)^n$. 
We expect that $\FC_{\beta}$ is in fact a coherent sheaf, and its restriction to $X_{n}(\C^2,\C)$ agrees with the conjectures of \cite{GNR,OR1,OR2,OR3}.

\begin{remark}
We have natural projection maps $\Hilb^n(\C^2)\rightarrow \Sym^n(\C^2)$, $X_n(\C^2)\rightarrow (\C^2)^n$, and $\Hilb^n(\C^2,\C)\rightarrow \Sym^n(\C\times\{0\})$, $X_n(\C^2,\C)\rightarrow (\C\times\{0\})^n$.  This maps endow sheaves on these schemes with the structure of modules over
\[
\C[\xx,\yy]^{S_n},\qquad \C[\xx,\yy],\qquad \C[\xx]^{S_n}, \qquad \C[\xx],
\]
respectively, where the $S_n$ action on $\C[\xx,\yy]$ permutes both sets of variables simultaneously.

Meanwhile the choice of braid representative $\b$ for $L$ allows us to identify $\C[x_c]_{c\in \components(L)}$ with a quotient of $\C[x_1,\ldots,x_n]$ in which $x_i\mapsto x_c$ if the $i$-th strand in $\b$ is on component $c\in \components(L)$ (similarly for $y_c$).  Thus, to obtain link homology with its module structure ($y$-ified or not), one should work with the isospectral Hilbert schemes.
\end{remark}

Remark that cohomology of sheaves on $X_n$ are $\Q[\xx,\yy]$-modules, while cohomology of sheaves on $X_n(\C^2,\C)$  are only $\Q[\xx]$-modules.

Next we consider another interesting corollary of Theorem \ref{th: intro ft}. It follows from the work of Haiman \cite{haiman2002vanishing} that the bigraded character of $J_{n}^k$ equals the Hall inner product $(\nabla^k p_1^{n},e_n)$, where $e_n$ is the elementary symmetric function and $\nabla$ is the Bergeron-Garsia operator on symmetric functions (see \S \ref{sec:combinatorics}). This yields yet another combinatorial formula for the Poincar\'e series of $\HKR(\FT_n^k)$  which was first conjectured in  \cite{EH} (see also \cite{Wilson}). Here we just write a formula for $k=1$ and the $\deg_A=0$ component, and
give the complete combinatorial description in \S \ref{sec:combinatorics}.

\begin{theorem}
\label{th: intro magic formula}
The Poincar\'e polynomial for $(a=0)$ part of the HOMFLY-PT homology of $\FT_n$ equals
\[
\PC(\FT_n^k)(q,t,a)|_{a=0}=(1-t)^n(\nabla^k p_1^{n},e_n)=\sum_{(e_1,\ldots,e_n)\in (\Z_{\ge 0})^n }q^{\sum e_i}t^{\dinv(e)},
\]
where $\dinv(e)=\sharp\{i<j: e_i=e_j\ \text{or}\ e_i+1=e_j\}$ and we have expressed the Poincar\'e series in terms of the variables \eqref{eq:qat}.
\end{theorem}

\subsection{Comparison with the results of Batson-Seed}
\label{subsec:introcomparison}
As we have mentioned, our construction of $y$-ified homology and some of the results are similar to the work of Batson and Seed \cite{Batson}, although they work in Khovanov homology while we consider HOMFLY-PT homology. For example, our Theorem \ref{th:intro splitting} and Corollary \ref{cor:intro ss} are analogous to \cite[Theorem 1.1]{Batson}. 

However, there are several important differences:
\begin{itemize}
\item As with the undeformed homologies, our approach requires a braid presentation for a link, while they can use arbitrary link diagrams.

\item For most of the paper, we regard $y_i$ as formal variables rather than elements in a field. This allows us to keep all three gradings on the deformed homology and study the $\Q[\yy]$--module structure. We also study the dependence of the specializations $\HY(L,\Q_{\nu})$ 
on the parameters $\nu\in \Q^r$.

\item Our construction is {\em local} and defined for braids. As explained in \cite{Batson}, the proof of the equation $d^2=0$ for links 
requires a certain non-local cancellation of terms. This means that for braids or tangles $d^2\neq 0$ and the local picture of the theory naturally involves curved complexes with $d^2=Z(\xx,\yy)$ for some function $Z$, see section \S \ref{sec:setup}.  This is the key technical novelty of our approach. \footnote{In fact, the authors of \cite{Batson} raised a question of constructing a version of their complex for tangles using curved complexes.} 

\item The parity phenomenon is unique to HOMFLY-PT homology: the Khovanov homology of the full twist and other interesting braids are 
supported both in even and in odd homological degrees. In particular, Theorems \ref{th: intro splitting injective} and  \ref{th: intro splitting ideal}   cannot be applied directly to compute the Khovanov homology of the full twist. 
\end{itemize}

We expect that there exist similar $y$-ified pictures for categories of $\sl_N$ matrix factorizations and $\sl_N$ link homology, where there is no need to restrict to braid presentations of a link.  We leave such investigations for future work.

\subsection{Conventions}
\label{subsec:conventions}
The HOMFLY-PT polynomial of links is the polynomial $P(L)(Q,\a)$ which is uniquely characterized by the skein relation
\[
\a P\left(\:\ig{.6}{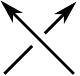}\:\right) - \a\inv P\left(\:\ig{.6}{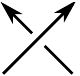}\:\right)= (Q-Q\inv)P\left(\:\ig{.6}{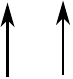}\:\right)
\]
together with $P_{L\sqcup L'}=P_LP_{L'}$, where $L\sqcup L'$ denotes the split union.  These rules force the invariant of the unknot to be $\frac{\a-\a\inv}{Q-Q\inv}$.

\begin{remark}
In a link diagram, a \emph{negative} (or \emph{left-handed} crossing is $\:\ig{.6}{lcrossing}\:$, while a \emph{positive} (or \emph{right-handed}) crossing is $\:\ig{.6}{rcrossing}\:$.
\end{remark}

To define Khovanov-Rozansky homology of a link $L$, first choose a braid representative and consider $\HKR(\b)$ whose definition we reviewed above.

\begin{remark}
We often prefer to write degrees multiplicatively; if $c\in \HKR^{i,j,k}(\b)$, then we say that $c$ has \emph{weight} $\wt(c)=Q^iT^jA^k$.
\end{remark}

To obtain a properly normalized link homology define
\[
(\Hnorm)^{i,j,k}(L) = \HKR^{i',j',k'}(\b),
\]
where
\[
i=i'-e+2n-2r,\qquad j=j'+\frac{1}{2}(e-n+r),\qquad k=k'+\frac{1}{2}(-e-n+r).  
\]
Here $n$ is the number of strands in $\b$, $r$ is the number of components of $L$, and $e=e(\b)$ is the signed number of crossings, i.e.~the image of $\b$ under the group homomorphism $e:\Br_n\rightarrow \Z$ sending each $\sigma_i\mapsto 1$.  It is straightforward to show that $e$ always has the same parity as $n-r$, hence the above is still a $\Z\times \Z\times \Z$ grading.  Let $\PC_\b(Q,A,T)$ denote the Poincar\'e series of $\HKR(\b)$ and $\Pnorm_L(Q,A,T)$ denote the Poincar\'e series of $\Hnorm(L)$, so that
\[
\Pnorm_L:= Q^{-e+2n-2r} A^{(e-n+r)/2}T^{(-e-n+r)/2}\PC(\b).
\]
Then $\Pnorm(L)$ is an honest link invariant, and the relationship with the HOMFLY-PT polynomial is
\[
\Pnorm_L(Q,-\a^2Q^2,-1) = Q^{-2r}\a^r P(L)(Q,\a).
\]
For the reader's convenience, we mention that the superpolynomial in \cite{DGR} is related to the Poincar\'e series $\Pnorm$ by
\[
Q^{2r}\a^{-r}\PC_L|_{A=\a^2T\inv Q^2}
\]
with the following exceptions.
\begin{itemize}
\item  We prefer a cohomological convention for differentials, hence our $T$ is $t_{\DGR}\inv$.
\item Our Poincar\'e series are \emph{unreduced}, meaning that the invariant of the unknot is not 1.
\item What we call the HOMFLY-PT polynomial of $L$ is called the HOMFLY-PT polynomial of the mirror $\bar{L}$ in \cite{DGR}.
\end{itemize}

This last remark deserves some attention.  In most link homology theories, there is a duality between the homology of a link and its mirror.  In particular, over $\Q$, the homology of $L$ and $\bar{L}$ are dual to one another as (appropriately graded) vector spaces.  However, for the HOMFLY-PT homology this symmetry is much more subtle, owing essentially to the fact that the invariant of the unknot is no longer a Frobenius algebra.

Our preferred convention regarding links and their mirrors is determined by either of the equivalent descriptions:
\begin{enumerate}
\item The bigraded $R$-module of homogeneous chain maps $\FT_n\rightarrow R$ in $\KC(\SBim_n)$ up to homotopy is isomorphic to $R$ in degree zero, while the bigraded $R$-module of maps the other direction $R\rightarrow \FT_n$ is much more interesting\footnote{This module is nothing other than the Hochschild degree zero part of the homology of the $(n,n)$-torus link, which was computed in \cite{EH}}.
\item The HOMFLY-PT homology of every positive torus link is supported in even homological degrees, while the homology of negative torus links fails to be parity already for the $(2,-4)$ link.
\item The idempotents obtained as ``infinite torus braids $\FT_n^{\otimes \infty}$'' are unital idempotents in the language of \cite{Hog17b}, rather than counital idempotents.
\end{enumerate}

The conjectural $Q\leftrightarrow TQ\inv$ symmetry of $\PC_L$ is therefore a special case of a symmetry first conjectured in \cite{DGR}.

\begin{remark}To obtain a precisely normalized $y$-ified homology, one regrades $\HY(\b)$ in precisely the same way as $\HKR(\b)$.  Then properly normalized Poincar\'e series of $\HY(L)$ satisfies
\[
\PC^{\text{norm}}_{\HY(L)}(Q,-\a^2Q^2,-1) = \frac{Q^{-2r}\a^r}{(1-Q^{-2})^r} P(L)(Q,\a).
\]
\end{remark}

We will henceforth deal entirely with the un-normalized homology and Poincar\'e series, since carrying around the specific normalization factors can be quite a nuisance.  Further, the most compact and aesthetically pleasing formulae are usually obtained by using the unnormalized Poincar\'e series, expressed in terms of the variables
\begin{equation}\label{eq:qat}
q=Q^2 \qquad \qquad t=T^2Q^{-2}\qquad \qquad a=AQ^{-2}.
\end{equation}
For comparison, we have the invariants of the unknot:
\[
\PC(U) =\frac{1+a}{1-q} \qquad\qquad \Pnorm(U)=\frac{1+\a^2 T\inv}{1-Q^2}
\]
and the positive trefoil:
\[
\frac{\PC(T(2,3))}{\PC(U)} =q+t+a \qquad \qquad \qquad \frac{\Pnorm(T(2,3))}{\Pnorm(U)}= Q^{-2}\a^2 +Q^2T^{-2}\a^2 + T^{-3}\a^4.
\]

All complexes representing knot homology are defined over $\Q$, so we use rational coefficients almost everywhere in the paper. 
By the universal coefficients theorem, one can replace $\Q$ by any field of characteristic zero. In sections \S \ref{subsec:introHilb} and \S \ref{sec:Hilbert} we work over $\C$ to match the Hilbert scheme literature.

\subsection{Organization of the paper}

In Section \ref{sec:setup} we discuss general properties of matrix factorizations and define the notion of $y$-ification of an abstract complex.  

In Section \ref{sec:setup} we study  $y$-ifications of the Rouquier complex of Soergel bimodules.  In Theorem \ref{thm:rouquier} we prove that given a braid $\beta$, the corresponding Rouquier complex  $F(\beta)$ admits a canonical $y$-ification for the permutation associated with $\beta$. 

In Section \ref{sec:homology} we define the Hochschild cohomology for $y$-ified complexes and prove that the Hochschild cohomology of $y$-ified Rouquier complexes is invariant under Markov move (Proposition \ref{prop:markov}) and hence defines a topological link invariant $\HY(L)$. Theorem \ref{thm:rouquier} and Proposition \ref{prop:markov} imply Theorem \ref{th:intro invariance}. We  prove that $\HY(L)$ has a well-defined module structure over the $y$-ified homology of the unlink. We also define a family of {\em homologies with coefficients in $\Q_{\nu}$} $\HY(L,\Q_{\nu})$ for $\nu\in \Q^r$ (if $L$ has $r$ components)   by specializing $y_c$ to $\nu_c$ in the corresponding complex. In Theorem \ref{thm:deformationSS} we construct a spectral sequence with $E_1$ page isomorphic to the Khovanov-Rozansky HOMFLY-PT homology $\HKR(L)$ and (with one grading collapsed) and $E_{\infty}$ page isomorphic to $\HY(L,\Q_{\nu})$. This is a direct generalization of the Batson-Seed spectral sequence \cite{Batson}.

The next section \ref{sec:splitting} contains the new computational tools for HOMFLY-PT homology. We define a map $\Psi$ relating the $y$-ified homology of a link and of the corresponding split link. We use it to prove Theorem \ref{th:intro splitting} and Corollary \ref{cor:intro ss} as Corollary \ref{cor:linksplittingproperty}. 
Furthermore, we prove Theorem \ref{th: intro parity is flat} as Theorem \ref{th:flat vs collapsible}, Theorem \ref{th: intro splitting injective} as Theorem \ref{thm:splittableInjectivity}  and Corollary \ref{th: intro splitting ideal} as Proposition \ref{prop:purebraidProps}. 

In Section \ref{sec:ideals} we study the ideals corresponding to some pure braids. We prove Proposition \ref{prop: intro JM} as Proposition \ref{prop:jucys murphy}.  Section \S \ref{sec:ftideals} proves that the vector space of maps $\FY(\FT_n^k)\rightarrow \one$ is rank 1, generated by the splitting map, then uses this to give a complete description of the the full twist ideal.

In the following section \ref{sec:Hilbert} we give a more detailed exposition of the algebraic geometry of the Hilbert scheme of points and polygraph rings. We connect it to the ideals appearing in the previous section and prove Theorem \ref{th: intro ft} (see Corollary \ref{cor: full twist from AJ}) and Theorem \ref{th: intro Hilbert} (see Corollary \ref{cor: full twist from isospectral}).  By using the results of Haiman connecting the isospectral Hilbert scheme to combinatorics, we prove Theorem \ref{th: full magic formula} generalizing Theorem \ref{th: intro magic formula}.

In Section \ref{sec:symmetry} we discuss a conjectural symmetry of HOMFLY-PT homology. 
Finally, in the appendix \ref{sec:appendix} we discuss more homological algebra of $y$-ified complexes. 

\section*{Acknowledgments}

The authors would like to thank Erik Carlsson, Ciprian Manolescu, Andrei Negu\cb{t}, Jacob Rasmussen, and Paul Wedrich for the useful discussions.
The work of E.~G.~ was partially supported by the NSF grant DMS-1700814, Russian Academic Excellence Project 5-100 and the grant RSF-16-11-10018.  The work of M.~H.~ was supported by NSF grant DMS 1702274.

\section{Rouquier complexes and $y$-ification}
\label{sec:setup}
In this section we introduce the notion of a $y$-ification of a complex of Soergel bimodules, and we prove that Rouquier complexes admit unique $y$-ifications up to equivalence.

As a matter of notation, if $\AC$ is a $\Z$-linear category, we let $\KC(\AC)$ denote the homotopy category of complexes.  We prefer the cohomological convention for differentials, hence the differential of a complex $C$ will map $C^k\rightarrow C^{k+1}$ for all $k$.  We let $[1]$ denote the ``left shift'' of complexes, defined by $C[1]^k:=C^{k+1}$.  By convention the shift $[1]$ negates the differential $d_{C[1]}=-d_C$.

We let $\KC^b(\AC)$, $\KC^+(\AC)$, and $\KC^-(\AC)$ denote the categories of complexes which are bounded, respectively bounded from below, respectively bounded from above.


\subsection{The Soergel category and Rouquier complexes}
\label{sub:soergel rouquier}
Fix an integer $n\geq 1$.  Let $R=\Q[x_1,\ldots,x_n]$, and $R^e=R\otimes_\Q R$. 
We will identify
\[
R^e\cong \Q[x_1,\ldots,x_n,x_1',\ldots,x_n'],
\]
via an isomorphism that sends $x_i\otimes 1\mapsto x_i$ and $1\otimes x_i\mapsto x_i'$.  We regard $R$ and $R^e$ as graded algebras via $\deg(x_i)=\deg(x_i')=2$.  We use bold symbols to denote entire sets of variables, so the above polynomial rings might also be denoted $R=\Q[\xx]$ and $R^e\cong \Q[\xx,\xx']$.  The categories of graded $R^e$-modules and graded $(R,R)$-bimodules will henceforth be identified with one another.  We let $(1)$ be the shift of graded (bi)modules, so $M(1)^k = M^{k+1}$.

The symmetric group $S_n$ acts on $R=\Q[x_1,\ldots,x_n]$ by permuting indices.   For each subgroup $G\subset S_n$, let $R^G$ denote the subalgebra of polynomials $f\in R$ such that $g(f)=f$ for all $g\in G$.  If $G=\{1,s_i\}$, where $s_i=(i,i+1)$ is the simple transposition, then we also write $R^i = R^G$.  Let $B_i$ denote the graded $(R,R)$-bimodule
\[
B_i = R\otimes_{R^i} R(1).
\]
Let $\SBim_n$ denote the smallest full subcategory of $(R,R)\gbimod$ containing the identity bimodule $R$, the bimodules $B_1,\ldots,B_{n-1}$, and closed under isomorphism, tensor product $\otimes_R$, taking direct sums, direct summands, and grading shifts.  Objects of $\SBim_n$ are called \emph{Soergel bimodules}.

The category $\SBim_n$ is monoidal with respect to the tensor product $\otimes_R$, with identity the trivial bimodule $\one=R$.

\begin{notation}
Henceforth, we will denote $\otimes_R$ simply by $\otimes$.
\end{notation}

We have a bimodule map $b:B_i\rightarrow R(1)$ sending $1\otimes 1\mapsto 1$, and a bimodule map $b^\ast: R(-1)\rightarrow B_i$ sending $1\mapsto x_i\otimes 1- 1\otimes x_{i+1}$.  These will be referred to as the \emph{dot maps}.

The braid group $\Br_n$ is generated by the elementary braids (Artin generators) $\sigma_1,\ldots,\sigma_{n-1}$ modulo the usual braid relations
\[
\sigma_i\sigma_{i+1}\sigma_i=\sigma_{i+1}\sigma_{i}\sigma_{i+1},\ \sigma_i\sigma_j=\sigma_j\sigma_i\ \text{for}\ |i-j|\ge 2.
\]
Associated to $\sigma_i^\pm$ we have the \emph{Rouquier complexes}
\begin{subequations}
\begin{equation}
F(\sigma_i) \ \ = \ \
\begin{tikzpicture}[baseline=-.5em]
\tikzstyle{every node}=[font=\small]
\node (a) at (0,0) {$\underline{B_1}$};
\node (b) at (2,0) {$R(1)$};
\path[->,>=stealth',shorten >=1pt,auto,node distance=1.8cm]
(a) edge node[above] {$b$}		(b);
\end{tikzpicture}
\end{equation}
\begin{equation}
F(\sigma_i\inv) \ \ = \ \ \begin{tikzpicture}[baseline=-.5em]
\tikzstyle{every node}=[font=\small]
\node (a) at (0,0) {$R(-1)$};
\node (b) at (2,0) {$\underline{B_1}$};
\path[->,>=stealth',shorten >=1pt,auto,node distance=1.8cm]
(a) edge node[above] {$b^\ast$}		(b);
\end{tikzpicture}
\end{equation}
\end{subequations}

Associated to each braid word $\bbeta=(\sigma_{i_1}^{\e_1},\ldots,\sigma_{i_r}^{\e_r})$ with $\e_{i_j}\in \{1,-1\}$, let
\[
F(\bbeta) := F(\sigma_{i_1}^{\e_1}) \otimes \cdots \otimes  F(\sigma_{i_r}^{\e_r}).
\]
In \cite{Rouquier} it is proven that $F(\bbeta)$ depends only on the braid $\b=\sigma_{i_1}^{\e_1}\cdots\sigma_{i_r}^{\e_r}\in \Br_n$, and not on the chosen representation as a product of generators, up to canonical isomorphism in $\KC^b(\SBim_n)$.  

\subsection{Motivating examples}
As a motivating example for the abstraction which follows, let us consider the most important examples of $y$-ifications.  Let $\yy$ denote a set of variables $y_1,\ldots,y_n$.  We place a bigrading on $\Q[\xx,\xx',\yy]$ by declaring $\deg(x_i)=\deg(x_i')=(2,0)$ and $\deg(y_i)=(-2,2)$.  Written multiplicatively, this is $\wt(x_i)=\wt(x_i')=Q^2$ and $\wt(y_i)=T^2Q^{-2}$.  Given a (bigraded $\Q$-vector space $V$), we will let $V[\yy]:=V\otimes_\Q \Q[y_1,\ldots,y_n]$.  Consider the following diagrams:
\begin{subequations}
\begin{equation}\label{eq:firstycrossing+}
\FY(\sigma_1) \ \ := \ \ \begin{tikzpicture}[baseline=-.5em]
\tikzstyle{every node}=[font=\small]
\node (a) at (0,0) {$\underline{B_1}[\yy]$};
\node (b) at (5,0) {$R[\yy](1)$};
\path[->,>=stealth',shorten >=1pt,auto,node distance=1.8cm,
  thick]
([yshift=3pt] a.east) edge node[above] {$b\otimes 1$}		([yshift=3pt] b.west)
([yshift=-2pt] b.west) edge node[below] {$b^{\ast}\otimes (y_1-y_2)$}		([yshift=-2pt] a.east);
\end{tikzpicture}
\end{equation}
\begin{equation}\label{eq:firstycrossing-}
\FY(\sigma_1\inv) \ \ := \ \ 
\begin{tikzpicture}[baseline=-.5em]
\tikzstyle{every node}=[font=\small]
\node (a) at (0,0) {$R[\yy](-1)$};
\node (b) at (5,0) {$\underline{B_1}[\yy]$};
\path[->,>=stealth',shorten >=1pt,auto,node distance=1.8cm,
  thick]
([yshift=3pt] a.east) edge node[above] {$b^\ast\otimes 1$}		([yshift=3pt] b.west)
([yshift=-2pt] b.west) edge node[below] {$b\otimes (y_1-y_2)$}		([yshift=-2pt] a.east);
\end{tikzpicture}.
\end{equation}
\end{subequations}
We say that the above diagrams define $y$-ifications of the Rouquier complexes $F(\sigma_1^\pm)$.  Let us discuss how to interpret, for example, the first diagram.  We may extract a bigraded vector space from this diagram, of the form
\[
V = (B_1\oplus R(1)[-1])\otimes_\Q \Q[y_1,y_2].
\]
The homological shift $[-1]$ corresponds to the fact that the $R$ term is to the right of $B_1$ (and $[1]$ is the leftward shift).  Since $R$ and $B_1$ are graded $R^e$-modules, the vector space $V$ has the structure of a bigraded $R^e[\yy]$-module.  The arrows in the diagram \eqref{eq:firstycrossing+} determine a $R^e[\yy]$-linear map $\Delta:V\rightarrow V[1]$.  However, $(V,\Delta)$ is not a chain complex, since $\Delta^2\neq 0$.  On the contrary, direct computation shows that $b_1\circ b_1^\ast=x_1-x_2'$ acting on $R$, and $b_1^\ast\circ b_1=x_1-x_2'$ acting on $B_1$, hence
\[
\Delta^2= -(x_1-x_2')(y_1-y_2)  = (x_2-x_1')y_1 + (x_1-x_2')y_2 
\]
where we have used the fact that
\[
x_1 - x_2' = -(x_2 - x_1')
\]
as endomorphisms of $B_1$ and $R$.  Thus, the diagram \eqref{def:abstractYification} determines a sort of chain complex in which the differential squares to a nonzero, but central element of the category of bigraded $R^e[\yy]$-modules.  Such objects are called \emph{curved complexes}.  We now digress to introduce some abstract theory for discussing such objects.

\subsection{Graded categories}
\label{subsec:gradedCats}
Let $\Gamma$ be an abelian group, and let $\CC$ be a category with invertible functors $\Sigma_\gamma:\CC\rightarrow \CC$, $\gamma\in \Gamma$, such that $\Sigma_0 = \Id_\CC$ and $\Sigma_\gamma\circ \Sigma_{\gamma'}=\Sigma_{\gamma+\gamma'}$ for all $\gamma,\gamma'\in \Gamma$.  Then for objects $C,D\in \CC$ we will let
\[
\Hom_{\CC}^{\Gamma}(C,D):=\bigoplus_{\gamma\in\Gamma}\Hom_{\CC}(C,\Sigma_\gamma D).
\]
In this case we say that $\CC$ has a \emph{strict $\Gamma$-action}, or that $\CC$ is $\Gamma$-graded.  For us, $\Gamma$ will typically be $\Z$ or $\Z\times \Z$ with the $\Gamma$ action on $\CC$ determined by appropriate grading shift functors $\Sigma_i=(i)$ or $\Sigma_{i,j}=(i,j)$.

Let $S$ be a $\Gamma$-graded ring.  That is to say, $S=\bigoplus_{\gamma\in \Gamma} S^\gamma$, and ring multiplication restricts to a morphism $S^\gamma\otimes S^{\gamma'}\rightarrow S^{\gamma+\gamma'}$, and let $\CC$ be a $\Gamma$-graded category.  We say that $\CC$ is $S$-linear if $\CC$ is equipped with a strict $\Gamma$-action, each graded hom space $\Hom^\Gamma_{\CC}(C,D)$ is a graded $S$-module, and composition of morphisms is $S$-bilinear.  Equivalently, there is a map of graded algebras from $S$ in to the graded algebra of natural transformations of $\Id_\CC$:
\[
S\rightarrow \bigoplus_{\gamma\in \Gamma}\Hom(\Id_{\CC}, \Sigma_\gamma).
\]
In particular, each element $r\in S^\gamma$ corresponds to a morphism $\Sigma_\gamma(C)\rightarrow C$ for each $C\in \CC$.

\begin{example}
If $S=\bigoplus_{i\in \Z}S^i$ is a commutative graded ring in the usual sense, then the category of graded $S$-modules is $S$-linear with $\Sigma_i=(i)$.
\end{example}

If $\CC$ is any category, we let $\seq{\CC}$ denote the category whose objects are sequences $(C^i)_{i\in \Z}$ with $C^i\in \CC$.  A morphism $C^\ast \rightarrow D^\ast$ in $\seq{\CC}$ is by definition a sequence of morphisms $(f^i)_{i\in \Z}$ with $f^i\in \Hom_{\CC}(C^i,D^i)$.  If $\CC$ has a strict action of $\Gamma$, then $\seq{\CC}$ has a strict action of $\Gamma\times \Z$, via functors $\Sigma_{\gamma,j}$ defined by
\[
(\Sigma_{\gamma,j} C)^k = \Sigma_\gamma C^{k+j} .
\]
In other words, $\Sigma_{\gamma,0}$ is inherited from the action of $\Sigma_\gamma$ on $\CC$, while $\Sigma_{0,j}$ is the shift of sequences.  We also write $C[j]:=\Sigma_{0,j}(C)$.   In this language, a complex (i.e.~ an object of $\KC(\CC)$) is a pair $(C,d)$ where $C\in \seq{\CC}$ and $d:C\rightarrow C[1]$ is a morphism in $\seq{\CC}$ such that $d[1]\circ d =0$.   By convention, the homological shift negates the differential of complexes: $(C,d)[1] = (C[1],-d[1])$.

If $\CC$ is a $\k$-linear category for some commutative ring $\k$ and $C,D\in \KC(\CC)$ are complexes, then we let $\Homc^{\Gamma,\Z}(C,D)$ denote the chain complex with chain groups
\[
\Homc^{\Gamma,k}_{\CC}(C,D):=\prod_{j-i=k}\Hom^{\Gamma}(C^i,D^j)
\]
and differential given by the usual super-commutator $f\mapsto d_C\circ f - (-1)^k f\circ d_D$.  Note that $\Homc^{\Gamma}_\CC(C,D)$ is a complex of $\Gamma$-graded $\k$-modules, hence is $\Gamma\times \Z$-graded.  The homotopy category $\KC(\CC)$ has a strict action of $\Gamma\times \Z$, and $\Hom_{\KC(\CC)}^{\Gamma\times \Z}(C,D)$ is simply the homology of $\Homc^{\Gamma\times \Z}_{\KC(\CC)}(C,D)$.

\subsection{Graded curved complexes}
Let $\Gamma$ be an abelian group, $S$ a $\Gamma\times \Z$-graded ring and $\AC$ an $S$-linear category, as in the previous section. 

Let $z\in S$ be a homogeneous element with degree $\deg(Z)=(0,2)\in \Gamma\times \Z$.  A \emph{(graded) $Z$-factorization in $\AC$} is a pair $(V,\Delta)$ where $V\in \AC$ and $\Delta:V\rightarrow V[1]$ is a morphism in $\AC$ such that $\Delta[1]\circ \Delta  = Z$.  We will also refer to $(V,\Delta)$ as a \emph{(graded) curved complex} in $\AC$ with \emph{curvature} $Z$, and \emph{connection} $\Delta$.  By abuse, we usually refer to $\Delta$ as the differential.  The terminology comes from the notion of curved dg algebras and modules  \cite{Posit93, Posit11}.

If $(V,\Delta)$ and $(V',\Delta')$ are $Z$-factorizations in $\AC$, then the graded morphism space $\Hom_{\AC}^{\Gamma\times \Z}(V,V')$ has the structure of a chain complex, where the differential sends
\[
f\mapsto \Delta\circ f - (-1)^k f\circ \Delta=:[\Delta,f]
\]
whenever $f$ is a morphism $f:V\rightarrow \Sigma_{\gamma,k}(V')$.  To verify that this defines a differential, compute that
\begin{eqnarray*}
[\Delta,[\Delta,f]] &=& \Delta^2\circ f +(-1)^k \Delta\circ f\circ \Delta -(-1)^k\Delta\circ f\circ \Delta - f\circ \Delta^2\\
&=& Z\circ f - f\circ Z\\
&=& Z(\Id\circ f - f\circ \Id)\\
&=&0.
\end{eqnarray*}
where in the second to last step we have used the fact that $\AC$ is an $S$-linear category (hence $Z$ corresponds to an element in the center of $\AC$).  The complex of homs will also be denoted as $\underline{\Hom}_{\Fac(\AC,Z)}^{\Gamma\times \Z}(V,V')$.

Note that, in particular, a degee zero cycle $f\in \Hom_{\AC}^{\Gamma\times \Z}(V,V')$ is simply a morphism $f:V\rightarrow V'$ in $\AC$ such that $\Delta\circ f = f\circ \Delta$, and such a cycle is zero in homology if there exists a morphism $h:V\rightarrow V'[-1]$ such that $\Delta\circ h+ h\circ \Delta = f$.

\begin{definition}
Retain notation as above.  Let $\Fac(\AC,Z)$ denote the homotopy category of $Z$-factorizations in $\AC$.  The space of morphisms $(V,\Delta)\rightarrow (V',\Delta')$ in $\Fac(\AC,Z)$ is by definition the degree zero homology of $\Hom^{\Gamma\times \Z}(V,V')$.  We let $[1]$ denote the endofunctor $\Fac(\AC,Z)\rightarrow \Fac(\AC,Z)$ such that
\[
(V,\Delta)[1] = (V[1],-\Delta[1]).
\]
\end{definition}

\begin{example}
The category $\Fac(\seq{\CC},0)$ is equivalent to the homotopy category of complexes $\KC(\CC)$.  In this case the group $\Gamma$ plays no role, and we could assume $\Gamma=0$.
\end{example}

Isomorphism in $\Fac(\AC,Z)$ is called \emph{homotopy equivalence}, and is denoted by $\simeq$.  We say that $(V,\Delta)\in \Fac(\AC,Z)$ is \emph{contractible} if $(V,\Delta)\simeq 0$.  Equivalently, $(V,\Delta)$ is contractible if the identity $\Id_V$ is null-homologous in $\underline{\Hom}_{\Fac(\AC,Z)}^{\Gamma\times\Z}(V)$.

One can define mapping cones in $\Fac(\AC,Z)$ by direct analogy with the case of chain complexes.  Let $(M,\Delta_M)$ and $(N,\Delta_N)$ be $Z$-factorizations in $\AC$ and let $f:M\rightarrow N$ be a morphism.    The \emph{mapping cone} of $f$ is
\[
\Cone(f) = \left(M[1]\oplus N , \smMatrix{-\Delta_{M} & 0 \\ f & \Delta_{N}}\right),
\]
which is a $Z$-factorization.  This definition of mapping cones gives $\Fac(\AC,Z)$ the structure of a triangulated category with suspension $[1]$. We will not need (or prove) this fact.  However, we will need the following consequence, which is proven in the appendix, \S \ref{subsec:cones}.

\begin{lemma}\label{lemma:coneDetectsHomotopyEquiv}
A morphism $f:(M,\Delta_M)\rightarrow (N,\Delta_N)$ is an isomorphism in $\Fac(\AC,Z)$ if and only if $\Cone(f)\simeq 0$.
\end{lemma}

In the next section we define the notion of a $y$-ification, which is a particular kind of curved complex.

\subsection{$y$-ifications}
\label{subsec:yifications}
Let $\CC$ be a graded $R^e$-linear category (for example, we could take $\CC=\SBim_n$ or $\CC=R^e\gmod$).  

\begin{definition}
Given any sequence $C\in \seq{\CC}$, we have the ``formal tensor product''
\[
C[\yy]:=C\otimes_\Q[y_1,\ldots,y_n]\cong \bigoplus_{k_1,\ldots,k_n\geq 0} C(2\ell)[-2\ell].
\]
where $\ell=k_1+\cdots+k_n$.  If $C$ is bounded below, then this direct sum is finite in each degree, hence is a defines a (bigraded) $\Q[\yy]$-module in $\seq{\CC}$.
\end{definition}


\begin{definition}\label{def:abstractYification}
Retain notation as above.  Given $w\in S_n$, let $Z_w:=\sum_i (x_{w(i)}-x_i')y_i\in \Q[\xx,\xx',\yy]$.  A \emph{$y$-ification} of $C\in \KC^b(\CC)$ is a triple $\CB=(C,w,\Delta)$ such that $(C[\yy],\Delta)$ is a curved complex in $\seq{\CC}$ with curvature $Z_w$ and $\Delta=d_C$ modulo the ideal $(\yy)$.   Let $\YC(\CC)$ denote the homotopy category of $y$-ifications.  An object in $\YC(\CC)$ is a $y$-ification $(C,w,\Delta)$, and morphisms in $\YC(\CC)$ are morphisms of curved complexes (up to homotopies).
\end{definition}

The following is obvious after unpacking the definitions.
\begin{proposition}
Let $\CC$ and $\DC$ be graded $R^e$-linear categories, and  $\FC:\CC\rightarrow \DC$ be a functor compatible with the shifts $(1)$ and the $R^e$-action on morphism spaces.  Then applying $\FC$ term-by-term induces a functor $\YC(\CC)\rightarrow \YC(\DC)$.\qed
\end{proposition}



\begin{remark}
A complex $C$ might have many inequivalent $y$-ifications, or it might have none.  Further, the permutation $w$ is not necessarily uniquely determined by $C$, in the sense that $C$ might have $y$-ifications $(C,w,\Delta)$ and $(C,w',\Delta')$ associated to different permutations $w\neq w'$. 
\end{remark}

Let $(C, w, \Delta)$ be a $y$-ification of $C\in \KC^b(R^e\gmod)$.   Given $\kk=(k_1,\ldots,k_n)\in \Z_{\geq 0}^n$, we let $|\kk|=k_1+\cdots+k_n$ and $\yy^\kk = y_1^{k_1}\cdots y_n^{k_n}$.  We may express $\Delta$ in terms of its \emph{components}:
\begin{equation}
\Delta = \sum_{\kk} \Delta_{\kk} \otimes \yy^{\kk}
\end{equation}
where $\Delta_\kk \in \End^{\Z\times \Z}_{R^e}(C)$ are endomorphisms of degree $(2|\kk|,1-2|\kk|)$.  By definition, $\Delta=d_C\otimes 1$ modulo $y_1,\ldots,y_n$, which is equivalent to $\Delta_0 = d_C$.

If $\kk = (0,\ldots,0,1,0,\ldots,0)$ (with 1 in position $i$) then we will set $h_i:=\Delta_\kk$, so that
\[
\Delta = d_C\otimes 1 + (h_1\otimes y_1+\cdots h_n\otimes y_n) + (\text{higher}).
\]
 The equation $\Delta^2 = \sum_i (x_{w(i)}-x_i')y_i$ implies that
\[
d_C\circ h_i + h_i\circ d_C = (x_{w(i)}-x_i')\otimes y_i
\]
so a $C$ admits a $y$-ification with permutation $w$ only if left multiplication by $x_{w(i)}$ is homotopic to right multiplication by $x_i$ as endomorphisms of $C$ (for all $i$).

The components $\Delta_{\kk}$ with $|\kk|=2$ are homotopies which realize the fact that $h_ih_j+h_jh_i$ and $h_i^2$ are null-homotopic in $\End^{4,-2}(C)$.  If $|\kk|\geq 3$, the components $\Delta_\kk$ are interpeted as certain higher homotopies.

\begin{definition}\label{def:strict}
A $y$-ification $(C,w,\Delta)$ is \emph{strict} if $\Delta_\kk=0$ for $|\kk|\geq 2$.
\end{definition}

\begin{proposition}
Given $C\in \KC^b(R^e\gmod)$ with differential $d$, the structure of a strict $y$-ification $(C[\yy],w,\Delta)$ is equivalent to a choice of elements $h_i\in \End^{2,-1}_{R^e}(C)$ such that
\begin{subequations}
\begin{equation}\label{eq:homotopy}
d h_i + h_id = x_{w(i)}\otimes 1 - 1\otimes x_i.
\end{equation}
\begin{equation}
h_ih_j+h_jh_i = 0.\label{eq:strictH}
\end{equation}
\begin{equation}
h_i^2=0\label{eq:strictH2}
\end{equation}
\end{subequations}
\end{proposition}
\begin{proof}
If $\Delta_\kk=0$ for $|\kk|\geq 2$, then $\Delta$ may be written in terms of components as
\[
\Delta = d_C\otimes 1 + \sum_i h_i\otimes y_i
\]
for some elements $h_i\in \End^{2,-1}_{R^e}(C)$.  It is straightforward to verify that the equation $\Delta^2=0$ is equivalent to \eqref{eq:homotopy}, \eqref{eq:strictH}, and \eqref{eq:strictH2}. 
\end{proof}

\begin{remark}\label{rmk:balanced}
If $B\in \SBim_n$ is a Soergel bimodule, then $f(x_1,\ldots,x_n)-f(x_1',\ldots,x_n')$ acts by zero on $B$ for any symmetric polynomial $f$.  In particular, $x_1+\cdots+x_n=x_1'+\cdots+x_n'$.  Thus, for a complex $C\in \KC^b(\SBim_n)$ it makes sense to ask for $y$-ifications $(C,w,\Delta)$ whose curved differentials $\Delta\in \End(C\otimes_\Q\Q[y_1,\ldots,y_n])$ satisfy the stronger condition that $\Delta=d_C\otimes 1$ modulo the ideal of differences $y_1-y_2,\ldots,y_{n-1}-y_n$.  On the level of strict $y$-ifications this amounts to imposing
\begin{equation}
h_1+\cdots+h_n = 0\label{eq:strictH3}
\end{equation}
In fact all of the $y$-ifications in this paper satisfy this extra condition.
\end{remark}

\subsection{An example}
\label{subsec:FT2}
The Rouquier complex $F(\sigma_1^2)$ is homotopy equivalent to the complex
\[
C \ \ = \ \ \begin{tikzpicture}[baseline=-.5em]
\tikzstyle{every node}=[font=\small]
\node (a) at (0,0) {$B_1(-1)$};
\node (b) at (3,0) {$B_1(1)$};
\node (c) at (6,0) {$R(2)$};
\path[->,>=stealth',shorten >=1pt,auto,node distance=1.8cm]
(a) edge node[above] {$x_1-x_1'$}		(b)
(b) edge node[above] {$b$}(c);
\end{tikzpicture}
\]
(we leave it as an exercise to verify this fact, see \cite{Kh07}).  Let us construct a $y$-ification of this complex.  First, consider $C\otimes_\Q\Q[y_1-y_2]$ with a curved differential $\Delta'$ expressed in the following diagram:
\[
C[y_1-y_2] \ = \ \left(
\begin{tikzpicture}[baseline=-8em]
\tikzstyle{every node}=[font=\small]
\node (a1) at (0,0) {$\underline{B_1}(-1)$};
\node (b1) at (3,0) {$B_1(1)$};
\node (c1) at (6,0) {$R(2)$};
\node (a2) at (0,-2) {$B_1(1)$};
\node (b2) at (3,-2) {$B_1(3)$};
\node (c2) at (6,-2) {$R(4)$};
\node (a3) at (0,-4) {$B_1(3)$};
\node (b3) at (3,-4) {$B_1(5)$};
\node (c3) at (6,-4) {$R(6)$};
\node (a4) at (0,-6) {$\vdots$};
\node (b4) at (3,-6) {$\vdots$};
\node (c4) at (6,-6) {$\vdots$};
\path[->,>=stealth',shorten >=1pt,auto,node distance=1.8cm]
(a1) edge node[above] {$x_1-x_1'$}		(b1)
(b1) edge node[above] {$b$}(c1)
(a2) edge node[above] {$x_1-x_1'$}		(b2)
(b2) edge node[above] {$b$}(c2)
(a3) edge node[above] {$x_1-x_1'$}		(b3)
(b3) edge node[above] {$b$}(c3)
(b1) edge node[above] {$\Id$} (a2)
(b2) edge node[above] {$\Id$} (a3)
(b3) edge node[above] {$\Id$} (a4);
\end{tikzpicture}
\right)
\]
The top row of the diagram in parentheses denotes $C\otimes 1$, the second row denotes $C\otimes (y_1-y_2)$, and so on.  Multiplication by $y_1-y_2$ is visualized as the endomorphism of the above which identifies the first row with the second, the second with the third, and so on.  Because $y_1-y_2$ carries homological degree 2, the bimodules in the first column appear in homological degrees $0,2,4,\ldots,$. By inspection, we have
\[
(\Delta')^2 = (x_1-x_1')(y_1-y_2)
\]
Now we fatten this up to a $y$-ification $C\otimes_\Q \Q[y_1,y_2]=C\otimes_\Q \Q[y_1-y_2]\otimes_\Q \Q[y_1]$, with the curved differential $\Delta$ induced from $\Delta'$ (i.e.~$\Delta= \Delta'\otimes 1 $).  The curvature is
\[
(\Delta)^2 = (x_1-x_1')(y_1-y_2)
\]
since $x_1-x_1'=-(x_2-x_2')$ as endomorphisms of any Soergel bimodule in $\SBim_2$.  To faithfully draw a picture of the full $y$-ification would be unreasonably difficult on two dimensional paper, so instead we ``collapse the $y$-directions'' and indicate the $y$-ification simply by the diagram:
%
%
%
\begin{equation}
C[y_1,y_2] \ := \ \begin{tikzpicture}[baseline=-.5em]
\tikzstyle{every node}=[font=\small]
\node (a) at (-1,0) {$\underline{B_1}[y_1,y_2](-1)$};
\node (b) at (4,0) {${B}_1[y_1,y_2](1)$};
\node (c) at (8,0) {${R}[y_1,y_2](2)$};
\path[->,>=stealth',shorten >=1pt,auto,node distance=1.8cm]
([yshift=3pt] a.east) edge node[above] {$(x_1-x_1')\otimes 1$}		([yshift=3pt] b.west)
([yshift=-2pt] b.west) edge node[below] {$\Id\otimes (y_1-y_2)$}		([yshift=-2pt] a.east)
(b) edge node[above] {$b\otimes 1$} (c);
\end{tikzpicture}.
\end{equation}
This $y$-ification is strict, since the curved differential is linear in $y_1,y_2$.  To recover the homotopies $h_1,h_2$, one could ``differentiate'' the curved differential with respect to $y_i$.

\subsection{Tensor products}
\label{subsec:tensorproduct}
Note that in order to define the category of $y$-ifications $\YC(\CC)$, we only assume that $\CC$ is a graded $R^e$-linear category.  In particular we have not assumed that $\CC$ has any monoidal structure.  However, if $\CC$ is monoidal in a way which is compatible with the $R^e$-structure, then $\YC(\CC)$ has the structure of a monoidal category.  This extra compatibility is discussed next.  The reader who wishes to skip ahead to Definition \ref{def:tensor} may take $\CC$ to be the category of graded $(R,R)$-bimodules.

Suppose that $\CC$ has the structure of a graded monoidal category with grading shift $(1)$.  We will also assume we are given isomorphisms
\[
\one(1)\otimes C\cong C(1)\otimes C\otimes \one(1),
\]
natural in $C\in \CC$.\footnote{One should also impose additional coherence relations, which we do not discuss here.}

Assume that we are given a map of graded algebras $\rho(R)\rightarrow \End_{\CC}^\Z(\one)$.  This makes $\CC$ into a graded $R^e$-linear category, where the action of $x_i$ on an object $C$ is via
\[
C \cong \one \otimes C \buildrel\rho(x_i)\otimes\Id_C\over  \rightarrow \one(2)\otimes C\cong C(2),
\]
and the action of $x_i'$ on $C$ is defined similarly using $\Id_C\otimes \rho(x_i)$.  To define the monoidal structure we only need the observation that on any tensor product of objects $C\otimes D$ in $\CC$ we have $x_i'\otimes \Id_D =\Id_C\otimes x_i$ as morphisms
\[
C\otimes D\rightarrow C(2)\otimes D\cong C\otimes D(2).
\]

\begin{definition}\label{def:tensor}
Let $\CC$ be as above, and let $\CB=(C,v\Delta_C)$ and $\DB=(D,w,\Delta_D)$ be objects of $\YC(\CC)$.  Then the \emph{tensor product $y$-ification} is $\CB\otimes \DB=(C\otimes D, vw, \Delta_{C\otimes D})$, where $\Delta_{C\otimes D}$ is defined in terms of components as
\[
(\Delta_{C\otimes D})_{\kk} = (\Delta_C)_{w(\kk)}\otimes \Id_D+\Id_C\otimes (\Delta_D)_{\kk} \in \Endc^{\Z}(C\otimes D)
\]
where $\kk=(k_1,\ldots,k_n)$ and $w(\kk) = (k_{w(1)},\ldots,k_{w(n)})$.
\end{definition}

Let us interpret the above definition.  Choose a permutation $w\in S_n$.  Let $D\in \CC$ be given, and consider $D[\yy]\in \seq{\CC}$.  Define a map of bigraded algebras $\sigma:\Q[\xx,\yy]\otimes \Q[\xx,\yy]\rightarrow \End_{\seq{\CC}}^{\Z\times \Z}(D[\yy])$ where $\sigma(x_i\otimes 1)$ and $\sigma(1\otimes x_i)$ are the actions of $x_i$ and $x_i'$ defined already, $\sigma(1\otimes y_i)$ is multiplication by $y_i$, and $\sigma(y_i\otimes 1)$ is multiplication by $y_{w\inv(i)}$.   We regard this algebra map as giving $D$ the structure of a graded $(\Q[\xx,\yy],\Q[\xx,\yy])$-bimodule in $\seq{\CC}$.

Then the tensor product of $y$-ifications is simply $C[\yy]\otimes_{\one[\yy]} D[\yy]$ with its tensor product curved differential $\Delta_C\otimes \Id + \Id\otimes \Delta_D$.  We leave it as an exercise to verify that one recovers the above definition.


\begin{lemma}\label{lemma:tensor}
Retain notation as in Definition \ref{def:tensor}.  The tensor product $(\CB\otimes \DB,vw,\Delta_{C\otimes D})$ is a $y$-ification.
\end{lemma}
\begin{proof}
%
First form the tensor product in $\seq{\CC}$:
\[
C[\yy]\otimes_{\one} D[\yy].
\]
The subscript `$\one$' indicates that we are performing the tensor product in $\seq{\CC}$ and we are not yet identifying the $y$'s on the left factor with $y$'s on the right.

Consider $\Delta_C\otimes \Id_{D[\yy]} + \Id_{C[\yy]}\otimes \Delta_D$, interpreted via the Koszul sign rule.  This is a morphism in $\seq{\CC}$:
\[
C[\yy]\otimes_{\one} D[\yy]\rightarrow (C[\yy]\otimes_{\one} D[\yy])[1].
\]
The sign rule ensures that the two summands anti-commute, hence
\[
(\Delta_C\otimes \Id_{D[\yy]} + \Id_{C[\yy]}\otimes \Delta_D)^2 = \sum_{i}(x_{v(i)}-x_i')y_i \otimes \Id + \sum_{i}\Id\otimes (x_{i}-x_{w\inv(i)}')y_{w\inv(i)},
\]
The action of $x_i'$ on the first tensor factor coincides with the action of $x_i$ on the second factor, hence half the terms above cancel upon identifying the $y_i$ on the left factor with $y_{w\inv(i)}$ on the second factor.  This is exactly the identification that occurs when forming $C[\yy]\otimes_{\one[\yy]} D[\yy]$.  Collecting the surviving terms gives
\[
(\Delta_C\otimes \Id_{D[\yy]} + \Id_{C[\yy]}\otimes \Delta_D)^2 = \sum_{i}(x_{vw(i)}-x_i')y_i.
\]
This shows that $\CB\otimes \DB$ is a $y$-ification of $C\otimes D$, with permutation $vw$.
\end{proof}

\begin{example}
Let's describe explicitly the tensor product of $y$-ified Rouquier complexes $\FY(\sigma_1)\otimes \FY(\sigma_1)$ defined by the diagram \eqref{eq:firstycrossing+}.  First, we label the strands with the numbers $1,2$, and form the tensor product
\[
 \left(\begin{tikzpicture}[baseline=-.5em]
\tikzstyle{every node}=[font=\small]
\node (a) at (0,0) {$\underline{B_1}[\yy]$};
\node (b) at (5,0) {$R[\yy](1)$};
\path[->,>=stealth',shorten >=1pt,auto,node distance=1.8cm,
  thick]
([yshift=3pt] a.east) edge node[above] {$b\otimes 1$}		([yshift=3pt] b.west)
([yshift=-2pt] b.west) edge node[below] {$-b^{\ast}\otimes (y_1-y_2)$}		([yshift=-2pt] a.east);
\end{tikzpicture}\right)
\otimes
 \left(\begin{tikzpicture}[baseline=-.5em]
\tikzstyle{every node}=[font=\small]
\node (a) at (0,0) {$\underline{B_1}[\yy]$};
\node (b) at (5,0) {$R[\yy](1)$};
\path[->,>=stealth',shorten >=1pt,auto,node distance=1.8cm,
  thick]
([yshift=3pt] a.east) edge node[above] {$b\otimes 1$}		([yshift=3pt] b.west)
([yshift=-2pt] b.west) edge node[below] {$b^{\ast}\otimes (y_1-y_2)$}		([yshift=-2pt] a.east);
\end{tikzpicture}\right)
\]
The sign on the differential on the first factor is coming from reordering of components.  Expanding the tensor product yields
\[
 \begin{tikzpicture}
\tikzstyle{every node}=[font=\small]
\node (aa) at (0,0) {$\underline{B_1\otimes B_1}[\yy]$};
\node (ba) at (5,0) {$B_1[\yy](1)$};
\node (ab) at (0,-5) {$B_1[\yy](1)$};
\node (bb) at (5,-5) {$R[\yy](2)$};
\path[->,>=stealth',shorten >=1pt,auto,node distance=1.8cm,
  thick]
([yshift=3pt] aa.east) edge node[above] {$b\otimes \Id\otimes 1$}		([yshift=3pt] ba.west)
([yshift=-3pt] ba.west) edge node[below] {$b^{\ast}\otimes \Id\otimes (y_1-y_2)$}		([yshift=-3pt] aa.east)
([xshift=-3pt] aa.south) edge node[left=3pt,yshift=11pt] {$-\Id \otimes b\otimes 1$} ([xshift=-3pt] ab.north)
([xshift=3pt] ab.north) edge node[right=3pt,yshift=11pt] {$\Id \otimes b^{\ast}\otimes (y_1-y_2)$} ([xshift=3pt] aa.south)
([xshift=3pt] ba.south) edge node[right=3pt,yshift=-11pt] {$b\otimes 1$} ([xshift=3pt] bb.north)
([xshift=-3pt] bb.north) edge node[left=3pt,yshift=-11pt] {$-b^{\ast}\otimes (y_1-y_2)$} ([xshift=-3pt] ba.south)
([yshift=-3pt] ab.east) edge node[below] {$b\otimes 1$}		([yshift=-3pt] bb.west)
([yshift=3pt] bb.west) edge node[above] {$b^{\ast}\otimes (y_1-y_2)$}		([yshift=3pt] ab.east);
\end{tikzpicture}
\]
Taking the direct sum along diagonals yields the following:
\[
 \begin{tikzpicture}
\tikzstyle{every node}=[font=\small]
\node (a) at (0,0) {$\underline{B_1\otimes B_1}[\yy]$};
\node at (7,.4) {$B_1[\yy](1)$};
\node at (7,-.4) {$B_1[\yy](1)$};
\node (b) at (7,0) {$\oplus$};
\node (c) at (12,0) {$R[\yy](2)$};
\path[->,>=stealth',shorten >=1pt,auto,node distance=1.8cm,
  thick]
([yshift=3pt] a.east) edge node[above] {$\sqmatrix{b\otimes \Id\\ -\Id\otimes b}$}		([yshift=3pt,xshift=-10pt] b.west)
([yshift=-3pt,xshift=-10pt] b.west) edge node[below] {$\sqmatrix{b^\ast\otimes \Id & \Id\otimes b^\ast}(y_1-y_2)$}		([yshift=-3pt] a.east)
([yshift=3pt,xshift=10pt] b.east) edge node[above] {$\sqmatrix{b & b}$}		([yshift=3pt] c.west)
([yshift=-3pt] c.west) edge node[below] {$\sqmatrix{-b^\ast \\ b^\ast}(y_1-y_2)$}		([yshift=-3pt,xshift=10pt] b.east);
\end{tikzpicture}.
\]
\end{example}

\begin{remark} It is perhaps better to regard the $y$-variables in a $y$-ification $(D,w,\Delta)$ as being associated to the strands of $w$.  Here a ``strand'' is interpreted in the usual string diagram notation for permutations, and can also be regarded as a pair $(i,j)$ with $i=w(j)$.  Fix a set $\Lab$ with $n$ elements, and label the strands of $w$ with elements of $\Lab$ (bijectively).  If $c\in \Lab$ is the label on the strand $(i,j)$, then we let $y_c$ denote the endomorphism of $D[\yy]$ which is left multiplication by $y_i$ (equivalently right multiplication by $y_j$).  Then the tensor product of $y$-ifications $(C,v,\Delta_C)\otimes (D,w,\Delta_D)$ is defined only when the labels on $v$ and $w$ are compatible, meaning they induce a well defined labeling on the strands of $vw$.  The tensor product of $y$-ifications with compatibly labeled permutations is now essentially the obvious tensor product, identifying $y$-variables with the same label.
\end{remark}

\subsection{Strict $y$-ifications of Rouquier complexes}
Consider the strict $y$-ifications depicted below:
\begin{subequations}
\begin{equation}\label{eq:ycrossing+}
\FY(\sigma_i) \ \ := \ \ \begin{tikzpicture}[baseline=-.5em]
\tikzstyle{every node}=[font=\small]
\node (a) at (0,0) {$\underline{B_i}[\yy]$};
\node (b) at (5,0) {$R[\yy](1)$};
\path[->,>=stealth',shorten >=1pt,auto,node distance=1.8cm,
  thick]
([yshift=3pt] a.east) edge node[above] {$b\otimes 1$}		([yshift=3pt] b.west)
([yshift=-2pt] b.west) edge node[below] {$b^{\ast}\otimes (y_i-y_{i+1})$}		([yshift=-2pt] a.east);
\end{tikzpicture}
\end{equation}
\begin{equation}\label{eq:ycrossing-}
\FY(\sigma_i\inv) \ \ := \ \ 
\begin{tikzpicture}[baseline=-.5em]
\tikzstyle{every node}=[font=\small]
\node (a) at (0,0) {$R[\yy](-1)$};
\node (b) at (5,0) {$\underline{B_i}[\yy]$};
\path[->,>=stealth',shorten >=1pt,auto,node distance=1.8cm,
  thick]
([yshift=3pt] a.east) edge node[above] {$b^\ast\otimes 1$}		([yshift=3pt] b.west)
([yshift=-2pt] b.west) edge node[below] {$b\otimes (y_i-y_{i+1})$}		([yshift=-2pt] a.east);
\end{tikzpicture}.
\end{equation}
\end{subequations}
Direct computation verifies that $\FY(\sigma_i^\pm)$ are $y$-ifications of $F(\sigma_i^\pm)$.  Tensoring these together (using the definition in \S \ref{subsec:tensorproduct}) defines $y$-ifications of arbitrary Rouquier complexes.

\begin{definition}
For each braid word $\bbeta = \sigma_{i_1}^{\e_1}\cdots \sigma_{i_r}^{\e_r}$ with $\e_{i_j}\in \{1,-1\}$, let
\[
\FY(\bbeta) := \FY(\sigma_{i_1}^{\e_1}) \otimes \cdots \otimes  \FY(\sigma_{i_r}^{\e_r}).
\]
\end{definition}

The remainder of this section is dedicated to proving the following.
\begin{theorem}\label{thm:rouquier}
The $y$-ification $\FY(\bbeta)$ depends only on the braid $\b$ up to canonical isomorphism in $\YC(\SBim_n)$.  In fact each Rouquier complex $F(\b)$ admits a unique $y$-ification up to homotopy equivalence.
\end{theorem}

We prove Theorem \ref{thm:rouquier} in \S \ref{subsec:uniqueness}, following a slight technical detour.


  \subsection{Contracting $y$-ifications}
\label{subsec:contractible}
In this section we wish to prove the following result.

\begin{lemma}\label{lemma:contractible}
If $C\in \KC^b(R^e\gmod)$ is a contractible complex then any $y$-ification of $C$ is contractible.
\end{lemma}
This, combined with Lemma \ref{lemma:coneDetectsHomotopyEquiv} gives us a method for proving that a morphism of $y$-ifications is a homotopy equivalence.  This will be crucial in the proof of Theorem \ref{thm:rouquier} below.
\begin{proof}
Let $h\in\End^{0,1}(C)$ be an element with $[d,h]=\Id_C$.  We may also assume that $h^2=0$.  Let
\[
\diff:=\Delta-d_C \in \End^1(C[\yy])
\]
denote the higher degree components of $\Delta$.  Let $(\yy)\subset \Q[\yy]$ denote the ideal generated by the $y_i$.  Observe that $C[\yy]$ is complete with respect to $(\yy)$ in the sense that any formal series of the form
\[
\sum_{\kk} f_{\kk} y^{\kk}
\]
gives rise to a well-defined endomorphism of $C[\yy]$, where $f_{\kk}\in \Endc(C)$ are any elements with
\[
\deg(f_{\kk}) = t^{-|\kk|} \deg(f_0).
\]
This is because $C$ is assumed bounded, hence all but finitely many of the $f_{\kk}$ must be zero for degree reasons.

Now, let $\Phi:=\Id+h\diff \in \End^0(C[\yy])$.  By the above discussion, $\Phi$ is invertible with inverse given by expanding the expression $(\Id+h\diff)\inv$ into a formal power series, since $\diff$ is zero modulo the maximal ideal $(\yy)$.

Now, define
\[
H:=(1+h \diff)\inv h = h(1+\diff h)\inv.
\]
We wish to show that $(d+\diff)H + H(d+\diff)=\Id$.  Multiplying on the left by $(1+h\diff)$ and on the right by $(1+\diff h)$, this is equivalent to
\[
(1+h\diff) (d+\diff)h + h(d+\diff)(1+\diff h) = (1+h\diff)(1+\diff h).
\]
Now, compute: the first term on the left is
\[
dh+h\diff d h + \diff h + h\diff^2 h
\]
The second term on the left is
\[
hd+ h d \diff h + h \diff + h\diff^2 h.
\]
Adding them gives
\[
dh+hd+ h(\diff d + d\diff)h + \diff h + h \diff+ 2 h\diff^2 h. 
\]
Now, $\diff d + d\diff = -\diff^2$, and $dh+hd=1$.  Thus, the above equals $1+\diff h + h\diff+h\diff^2 h$, as claimed.
\end{proof}

\subsection{Obstruction theoretic construction of morphisms}
\label{subsec:chainMapExtension}
In this section we consider a technical point which will be used in our proof of Theorem \ref{thm:rouquier} below.  Let $C,D\in \KC^b(R^e\gmod)$ be given.  Any homogeneous element $f\in \Homc_{R^e[\yy]}(C[\yy],D[\yy])$ with $\deg(f)=(a,b)$ may be written in terms of its components:
\[
f = \sum_{\kk\in\Z_{\geq0}^n} f_{\kk}\otimes y^\kk
\]
where $f\in \Homc_{R^e}(C,D)$ are linear maps of degree $\deg(f_{\kk})=(a+2|\kk|,b-2|\kk|)$.

Assume that $C[\yy]$ and $D[\yy]$ have the structures of $y$-ifications, associated with the same permutation $w$.  If $f:C[\yy]\rightarrow D[\yy]$ is a morphism in $\Fac(Z_w)$, then $f_0:C\rightarrow D$ is a chain map.  Unpacking the definitions, we see that $\Cone(f:C[\yy]\rightarrow D[\yy])$ is a $y$-ification of $\Cone(f_0:C\rightarrow D)$.  

Expressing $\Delta_D\circ f - f\circ \Delta_C=0$ in terms of components gives
\[
0=\sum_{\ii+\jj=\kk}(\Delta_{D,\ii}\circ f_\jj - f_\jj\circ \Delta_{C,\ii})\otimes \yy^{\ii+\jj}.
\]
We obtain the following system of equations:
\begin{subequations}
\begin{equation}
\label{eq:f0}[d,f_0] = 0.
\end{equation}
\begin{equation}
\label{eq:f1}[d,f_\kk] = -\sum_{\ii,\jj} \Big(\Delta_{\ii}f_\jj - f_{\ii}\Delta_{\jj}\Big).
\end{equation}
\end{subequations}
where the sum on the right-hand side of \eqref{eq:f1} is over nonzero sequences $\ii,\jj$ with $\ii+\jj=\kk$.

For each $\kk\in \Z_{\geq0}^n$, let $z(\kk)\in\Hom^{2|\kk|,-2|\kk|}(C,D)$ denote the right-hand side of \eqref{eq:f1}.  It is easy to see that if \eqref{eq:f1} holds for all $\lseq<\kk$, then $z(\kk)$ is a cycle.  Here we say $\lseq\leq \kk$ if the coordinates of $\kk-\lseq$ are non-negative, and $\lseq<\kk$ if $\lseq\geq \kk$ and $\lseq\neq \kk$. Thus, the construction of $f$ may be accomplished by induction, as follows:
\begin{enumerate}
\item Choose a chain map $f_0:C\rightarrow D$ in the usual sense, for $C,D\in \KC^b(\SBim_n)$.
\item Fix $\kk$ and assume by induction that we have defined $f_{\lseq}$ for all $\lseq<\kk$.
\item If $z(\kk)$ is a boundary, then we define $f_{\kk}$ by negating a chosen homotopy for $z(\kk)$.
\item Repeat steps (2) and (3).
\end{enumerate}

Since we assume $C$ and $D$ are bounded complexes, the $f_{\kk}$ are zero for all but finitely many $\kk$, and the above process eventually terminates.  

\subsection{Uniqueness of $y$-ifications of Rouquier complexes}
\label{subsec:uniqueness}

\begin{proposition}\label{prop:invertible}
Let $F,F'\in \KC^b(R^e\gmod)$ be invertible complexes and $\phi_0:F\rightarrow F'$ a homotopy equivalence.  Then $\phi_0$ extends (uniquely up to homotopy) to a homotopy equivalence $\phi:F[\yy]\simeq F'[\yy]$ for any $y$-ifications $(F[\yy],w,\Delta)$, $(F'[\yy],w,\Delta')$.
\end{proposition}
\begin{proof}
Since $F\simeq F'$ are invertible, the complex of homs satisfies
\[
\Homc(F,F')_{R^e}\simeq \Homc_{R^e}(\one,F'\otimes F\inv)\simeq \Endc_{R^e}(\one).
\]
In particular any chain map $\Homc(F,F')$ of nonzero homological degree is null-homotopic.  Thus, $\phi_0$ extends to a chain map $\phi:F[\yy]\rightarrow F'[\yy]$ by the procedure discussed in the previous subsection  \S \ref{subsec:chainMapExtension}.  Then $\Cone(\phi)$ is contractible by Lemma \ref{lemma:contractible}, since it is a $y$-ification of the contractible complex $\Cone(\phi_0)$.  It follows that $\phi$ is a homotopy equivalence of $y$-ifications.
\end{proof}

The main theorem of this section is an immediate corollary.

\begin{proof}[Proof of Theorem \ref{thm:rouquier}]
The Rouquier complexes are well-defined up to canonical isomorphism in $\KC^b(R^e\gmod)$.  Given the uniqueness of lifts in Proposition \ref{prop:invertible}, it follows that the $y$-ified Rouquier complexes $\FY(\b)$ are well-defined up to canonical isomorphism in $\KC(\Fac(Z_w))$.
\end{proof}

\section{The deformed homology theory}
\label{sec:homology}

In this section we define a homology theory for bigraded $Z_w$-factorizations, and we prove that the homology of a $y$-ified Rouquier complex is a well-defined invariant of the oriented link $L=\hat\b$, up to isomorphism and overall shift of triply graded modules over $\Q[x_c,y_c]_{c\in \pi_0(L)}$.

A permutation $w\in S_n$ determines an equivalence relation on $\{1,2,\ldots,n\}$, which is the transitive closure of the relation $j\sim i$ if $j=w(i)$.  Equivalence classes in $\{1,\ldots,n\}$ are also called \emph{cycles} of $w$.  We let $\Q[\xx,\yy]_w$ denote the quotient of $\Q[x_1,\ldots,x_n,y_1,\ldots,y_n]$ by the ideal generated by $x_i-x_j$ and $y_i-y_j$ for all $i\sim j$.

\begin{remark}\label{rmk:conjugation}
If $v,w\in S_n$, then we have a bijection between cycles of $w$ and cycles of $vwv\inv$ defined by
\[
\{1,\ldots,n\}_{\sim_{w}} \rightarrow \{1,\ldots,n\}_{\sim_{vwv\inv}},\qquad [i]\mapsto [v(i)].
\]
This induces an isomorphism $\Q[\xx,\yy]_w \cong \Q[\xx,\yy]_{vwv\inv}$.
\end{remark}

If $\b\in \Br_n$ is a braid with underlying permutation $w$, then we also let $\Q[\xx,\yy]_\b:=\Q[\xx,\yy]_w$.  Just as above, any expression of the form $\b'=\gamma\b\gamma\inv$ induces an isomorphism $\Q[\xx,\yy]_\b\cong \Q[\xx,\yy]_{\b'}$.

If $L$ is a link, then we let $\Q[\xx,\yy]_L$ denote the polynomial ring in variables $x_c$, $y_c$, indexed by components $c\in \pi_0(L)$.  If $\b\in \Br_n$ is a braid representative of $\b$, then we have a canonical isomorphism $\Q[\xx,\yy]_\b\cong \Q[\xx,\yy]_L$.


\subsection{The (Hochschild co)homology of a $y$-ification}
\label{subsec:cohomology}
If $B$ is any graded $R^e$-module, then we have the Hochschild cohomology of $B$, defined by
\[
\HH^{i,j}(B):= \Ext_{R^e\gmod}^j(R,B(i)).
\]
To compute these groups, choose a projective resolution $P\rightarrow R$ of $R$, as a graded $R^e$-module.  Then $\HH^{i,j}(B)$ is the $(i,j)$ homology group of $\Homc_{R^e}(P,B)$.  We also let $\HH(B)$ denote the bigraded vector space $\bigoplus_{i,j}\HH^{i,j}(B)$.


Now, let $\CB=(C,w,\Delta)\in \YC(R^e\gmod)$ be a $y$-ification of $C\in \KC^b(R^e\gmod)$.  We may write $\Delta$ in terms of its coordinates as in \S \ref{subsec:yifications}.  First, consider the complex
\[
\HH(C) = \begin{diagram} \cdots & \rTo^{\HH(d_C)} & \HH(C^k) & \rTo^{\HH(d_C)} &  \HH(C^{k+1})& \rTo^{\HH(d_C)} &  \cdots\end{diagram}.
\]
We regard $\HH(C)$ as being a triply graded chain complex, where the term in tridegree $(i,j,k)$ is $\HH^{i,j}(C^k)$.  An element $c\in \HH^{i,j}(M^k)$ will be said to have \emph{bimodule degree} $\deg_Q(c)=i$, \emph{Hochschild degree} $\deg_A(c)=j$, and \emph{homological degree} $\deg_T(c)=k$.  Note that the homology of $\HH(C)$, denoted $\HHH(C)$ in \cite{Kh07}, yields the Khovanov-Rozansky homology of $\b$ when $C=F(\b)$ is a Rouquier complex.  We prefer the notation $\HKR(C)$ over $\HHH(C)$ for purely aesthetic reasons.

Now, regard $\Q[y_1,\ldots,y_n]$ as a triply graded $\Q$-algebra by declaring $\deg(y_i)=(-2,0,2)$, and form the triply graded vector space
\[
\HH(\CB):=\HH(C)\otimes_\Q \Q[y_1,\ldots,y_n] 
\]
This vector space has a degree $(0,0,1)$ endomorphism of the form
\begin{equation}\label{eq:deltaPrime}
\Delta':=\sum_{k_1,\ldots,k_n} \HH(\Delta_{k_1,\ldots,k_n})\otimes y_1^{k_1}\cdots y_n^{k_n}.
\end{equation}
Since $\HH$ is a functor, it follows that 
\begin{equation}\label{eq:HHDsquared}
(\Delta')^2=\sum_{i=1}^n \HH(x_{w(i)}-x_i)\otimes y_i,
\end{equation}
since $\HH(x_i)=\HH(x_i')$ as endomorphisms of $\HH(C)$. 
\begin{remark}\label{rmk:HH(CB)}
We may think of $\HH(\CB)$ as a $y$-ification of $\HH(C)$.  To be precise, an $R$-module may be regarded as an $(R,R)$-bimodule on which the left and right actions coincide.  Then $\HH$ is an $R^e$-linear functor which takes (graded) $R^e$-modules to (bi)graded $R$-modules.  The assignment $\CB\mapsto \HH(\CB)$ is nothing other than the functor $\HH:R^e\mod^\Z\rightarrow R\mod^{\Z\times \Z}$, extended to the categories of $y$-ifications $\YC(R^e\mod^\Z)\rightarrow \YC(R\mod^{\Z\times \Z})$.  Here the superscripts $\Z$ and $\Z\times \Z$ are indicating the groups with respect to which the corresponding modules are graded.  
\end{remark}

After reindexing, \eqref{eq:HHDsquared} becomes
\[
(\Delta')^2=\sum_{i=1}^n \HH(x_i)\otimes (y_{w\inv(i)} - y_i),
\]
hence we arrive at the following.
\begin{lemma}
$(\Delta')^2=0$ modulo $(y_1-y_{w(1)},\ldots,y_n-y_{w(n)})$. \qed
\end{lemma}

\begin{definition}\label{def:CY}
Retain notation as above.  In particular $\CB=(C,w,\Delta)\in \YC(R^e\gmod)$ is a $y$-ification.  Let $\CY(C,w,\Delta)$ denote the complex
\[
\CY(\CB)= \HH(C)\otimes_\Q \Q[y_1,\ldots,y_n] / (y_1-y_{w(1)},\ldots, y_n-y_{w(n)}).
\]
with differential induced by \eqref{eq:deltaPrime}.  Let $\HY(\CB):=H(\CY(\CB))$ denote the homology.  We refer to $\HY(\CB)$ as the \emph{homology} of $\CB$.
\end{definition}

 \begin{example}[$n=1$]\label{ex:unknot}
Let $M=\Q[x,y]$, which is a $y$-ification of $R=\Q[x]$, with zero differential.  The Hochschild cohomology of $R$ is isomorphic to the super-polynomial ring $\Q[x,\theta]$ with $\wt(\theta)=AQ^{-2}$, as can be computed from the \emph{Koszul resolution} of $R$:
 \[
0\rightarrow  \Q[x,x'](-2)\buildrel x-x'\over\longrightarrow \Q[x,x'] \rightarrow \Q[x]\rightarrow 0.
 \]
To be more specific, $\HH(R)$ is isomorphic to the homology of the complex of homs
\[
\underline{\Hom}_{\Q[x,x']}^{\Z\times \Z}\Big(\Q[x,x'](-2)\buildrel x-x'\over\longrightarrow \Q[x,x'] ,\Q[x]\Big) \cong \underline{\Q[x]}\buildrel 0\over \longrightarrow  \Q[x](-2).
\]
This homology is isomorphic to $\Q[x]\oplus \Q[x](2,-1)$, which can also be written as $\Q[x,\theta]$ where $\deg(\theta)=(-2,1,0)$.  Thus, $\HH(R[y])= \HH(R)[y]\cong \Q[x,y,\theta]$.  Note that the Poincar\'e series of this homology is
\[
 \frac{1+a}{(1-q)(1-t)}
 \]
 in terms of the variables \eqref{eq:qat}.
 \end{example}
 
 \begin{remark}\label{rmk:knots}
 If $K$ is a knot, then $\HY(K)\cong \HKR(K)\otimes_\Q \Q[y_1]$ by Remark \ref{rmk:balanced}.
 \end{remark}
 
\begin{example}\label{ex:Hopf}
Let us compute our invariant for the positive Hopf link.  The Hopf link is the closure of $\sigma_1^2$, whose $y$-ification was studied in  \S \ref{subsec:FT2}.  We will focus on the Hochschild degree zero part.  First, apply the functor $\Hom_{R^e}(R,-)$ to each term.  This functor sends $R\mapsto R$ and $B_1\mapsto R(-1)$, since $\Hom(R,B_1)$ is isomorphic to $R$, generated by the dot map.  The resulting complex is
\[
\begin{tikzpicture}[baseline=-.5em]
\tikzstyle{every node}=[font=\small]
\node (a) at (0,0) {$\underline{R[\yy](-2)}$};
\node (b) at (3,0) {$R[\yy]$};
\node (c) at (6,0) {$R[\yy](2)$};
\path[->,>=stealth',shorten >=1pt,auto,node distance=1.8cm,
  thick]
([yshift=3pt] a.east) edge node[above] {$0$}		([yshift=3pt] b.west)
([yshift=-2pt] b.west) edge node[below] {$y_1-y_2$}		([yshift=-2pt] a.east)
(b) edge node[above] {$x_1-x_2$} (c);
\end{tikzpicture}
\]
The left-most rightward arrow is zero because $x_1-x_1'$ acts by zero on $R$.  The right-most rightward arrow is $x_1-x_2$ because the composition of dot maps $R\rightarrow B(1)\rightarrow R(2)$ is $x_1-x_2$.

The homology of this complex is generated over $R[y_1,y_2]$ by two elements $\a$, $\b$, where $\a$ and $\b$ are represented by 1 inside the left-most and right-most copies of $R$, respectively.  The left-most copy of $R$ appears in homological degree zero, hence $\wt(\a)=Q^2=q$.  The right-most copy of $R$ is in homological degree 2, hence $\wt(\b)=T^2Q^{-2}=t$. As a bigraded $\Q[\xx,\yy]$-module the homology is isomorphic to
\[
\Q[x_1,x_2,y_1,y_2]\cdot \{\a,\b\} \ \ \ \ \ \text{modulo } \ \ (y_1-y_2)\a = (x_1-x_2)\b.
\]
As a module over $\Q[\xx,\yy]$ the above has two generators (with weights $q,t$)  and one relation (with weight $qt$), and no syzygies.  Thus, the Poincar\'e series of the positive Hopf link at $a=0$ is
\[
\PC_{\operatorname{Hopf}}|_{a=0} = \frac{1}{(1-q)^2(1-t)^2}(q+t-qt).
\]
For the higher Hochschild degree components, see \S \ref{subsec:2strands}.
\end{example}

\subsection{Markov invariance of the $y$-ified homology}
\label{subsec:invariance}

Let us abbreviate by writing $\CY(\b)$ for $\CY(\FY(\b))$, and similarly for $\HY(\b):=H(\CY(\b))$.  We wish to prove that $\HY(\b)$ depends only on the link $L=\hat\b$, up to isomorphism and overall grading shift.

\begin{proposition}\label{prop:markov}
The $y$-ified homology is unchanged by Markov moves, up to isomorphism and overall shift: for all braids $\b,\gamma\in \Br_n$ we have 
\begin{enumerate}
\item $\CY(\b\gamma) \simeq \CY(\gamma\b)$,
\item $\CY(\b\sigma_n) \simeq \CY(\b)(1,0)[-1]$,
\item $\CY(\b\sigma_n\inv) \simeq \CY(\b)(3,1)[0]$,
\end{enumerate}
These equivalences are equivariant with respect to the action of $\Q[\yy]_{L}$ in each case.  For (1), the statement of equivariance involves the isomorphism $\Q[\yy]_{\b\gamma}\cong \Q[\yy]_{\gamma\b}$ from Remark \ref{rmk:conjugation}.
\end{proposition}
Statement (1) of Proposition \ref{prop:markov} follows from the following more general statement.

\begin{lemma}\label{lemma:traceproperty}
Let $\CB=(C,v,\Delta_C)$ and $\DB=(D,w,\Delta_D)$ be $y$-ifications of $C,D\in R^e\gmod$.  Then  $\CY(\CB\otimes \DB)\cong \CY(\DB\otimes \CB)$.  This isomorphism is equivariant with respect to the action of $\Q[\yy]_{vw}$, using the ring isomorphism $\Q[\yy]_{vw}\cong \Q[\yy]_{wv}$ from Remark \ref{rmk:conjugation}.
\end{lemma}
\begin{proof}
Recall that $\CB=C[y_1,\ldots,y_n]$ and $\DB=D[y_1,\ldots,y_n]$.  The tensor product $\CB\otimes \DB$ may be identified with
\[
\CB\otimes \DB\cong (C\otimes_R D)\otimes_\Q \Q[y_1,\ldots,y_n].
\]
The differential takes the form
\[
\sum_{k_1,\ldots,k_n}(\Delta_{k_{w(1)},\ldots,k_{w(n)}}\otimes \Id_D+\Id_C\otimes \Delta_{k_1,\ldots,k_n})\otimes y_1^{k_1}\cdots y_n^{k_n}.
\]
Modulo the ideal generated by $y_i-y_{vw(i)}$ for all $i$, this differential equals (after splitting up the sum and reindexing independently)
\[
\sum_{k_1,\ldots,k_n}(\Delta_{k_1,\ldots,k_n}\otimes \Id_D+\Id_C\otimes \Delta_{k_{v(1)},\ldots,k_{v(n)}})\otimes (y_1')^{k_1}\cdots (y_n')^{k_n},
\]
where $y_i'=y_{w\inv(i)}$.  Taking $\HH$ and using the fact that $\HH(C\otimes_R D)\cong \HH(D\otimes_R C)$ (naturally in $C$ and $D$), shows that $\CY(\CB\otimes \DB)\cong \CY(\CB\otimes \DB)$, as claimed.   This isomorphism is $\Q[y_1,\ldots,y_n]$-equivariant, provided that we twist the action on $\DB\otimes \CB$ by $w$.
\end{proof}

For the other Markov move (statements (2) and (3) of Proposition \ref{prop:markov}) we first recall a result on the Hochschild homology of Soergel bimodules, whose proof can be found in \cite{Kh07}.

\begin{lemma}\label{lemma:HHmarkov}
For every $M\in R^e_n\gmod$ we have
\begin{subequations}
\begin{equation}
\HH(M\sqcup\one_1) \cong \HH(M)[x_{n+1}](0,0)\oplus \HH(M)[x_{n+1}](2,-1)
\end{equation}
\begin{equation}
\HH((M\sqcup\one_1)\otimes B_n) \cong \HH(M)[x_{n+1}](-1,0)\oplus \HH(M)[x_{n+1}](3,-1),
\end{equation}
\end{subequations}
where $(i,j)$ denotes the shift of bigraded vector spaces, $V(i,j)^{k,l}=V^{i+k,j+l}$.  These isomorphisms are natural in $M$.  Furthermore, if $b:B\rightarrow R(1)$ and $b^\ast: R(-1)\rightarrow B$ are the dot maps, then in terms of these decompositions we have
\[
\HH(\Id_{M\sqcup\one_1}\otimes b)  \ \ =\ \ 
\begin{tikzpicture}[baseline=-.3em]
\tikzstyle{every node}=[font=\small]
\node (a) at (0,.7) {$\HH(M)[x_{n+1}](-1,0)$};
\node at (0,0) {$\oplus$};
\node (b) at (0,-.7) {$\HH(M)[x_{n+1}](3,1)$};
\node (c) at (5,.7) {$\HH(M)[x_{n+1}](1,0)$};
\node at (5,0) {$\oplus$};
\node (d) at (5,-.7) {$\HH(M)[x_{n+1}](3,1)$};
\path[->,>=stealth',shorten >=1pt,auto,node distance=1.8cm,
  thick]
(a) edge node {$x_n-x_{n+1}$} (c)
(b) edge node {$\Id$} (d);
\end{tikzpicture}
\]
and
\[
\HH(\Id_{M\sqcup\one_1}\otimes b^\ast)\  \ = \  \ 
\begin{tikzpicture}[baseline=-.3em]
\tikzstyle{every node}=[font=\small]
\node (a) at (0,.7) {$\HH(M)[x_{n+1}](-1,0)$};
\node at (0,0) {$\oplus$};
\node (b) at (0,-.7) {$\HH(M)[x_{n+1}](1,1)$};
\node (c) at (5,.7) {$\HH(M)[x_{n+1}](-1,0)$};
\node at (5,0) {$\oplus$};
\node (d) at (5,-.7) {$\HH(M)[x_{n+1}](3,1)$};
\path[->,>=stealth',shorten >=1pt,auto,node distance=1.8cm,
  thick]
(b) edge node {$x_n-x_{n+1}$} (d)
(a) edge node {$\Id$} (c);
\end{tikzpicture}
\]\qed
\end{lemma}

\begin{proof}[Proof of Proposition \ref{prop:markov}]
Statement (1) of the proposition is an immediate consequence of Lemma \ref{lemma:traceproperty}.  Statements (2) and (3) follow by a careful application of Lemma \ref{lemma:HHmarkov}.  We include the details only for (2), since (3) is similar.

Let $\CB=(C,w,\Delta)$ be a $y$-ification of $C\in K^b(\R^e_n\gmod)$.   We will show that $\CY((\CB\sqcup \one_1)\otimes \FY(\sigma_n))$ is homotopy equivalent to $\CY(\CB)$ up to overall shift.  Let $\yy$ denote the set of variables $y_1,\ldots,y_n$.  We may expand $\FY(\sigma_n)$ according to \eqref{eq:ycrossing+}, obtaining
\[
\Big((C\sqcup \one_1)\otimes_R F(\sigma_n)\Big)[\yy,y_{n+1}] \ \ \cong \ \ \hskip3in
\]
\[
\hskip1.1in \left(
\begin{tikzpicture}[baseline=-.5em]
\tikzstyle{every node}=[font=\small]
\node (a) at (0,0) {$((C\sqcup \one_1)\otimes_R B_n)[\yy,y_{n+1}]$};
\node (b) at (7,0) {$((C\sqcup \one_1)\otimes_R R)[\yy,y_{n+1}]$};
\path[->,>=stealth',shorten >=1pt,auto,node distance=1.8cm,
  thick]
([yshift=3pt] a.east) edge node[above] {$(\Id\otimes b)\otimes 1$}		([yshift=3pt] b.west)
([yshift=-2pt] b.west) edge node[below] {$(\Id\otimes b^{\ast})\otimes (y_n-y_{n+1})$}		([yshift=-2pt] a.east);
\end{tikzpicture}
\right)
\]
Throughout the remainder of the the proof we adopt the convention that a diagram enclosed in parentheses denotes the complex obtained by taking the direct sum of all the enclosed objects, with differential given by the sum of all visible arrows together with the differentials internal to each summand.

 Observe that $n$ and $n+1$ are in the same cycle of $(w\otimes 1)s_n$, hence $y_n$ and $y_{n+1}$ end up being identified in forming $\CY((\CB\sqcup \one_1)\FY(\sigma_n))$, hence the leftward arrow above becomes zero.  Now we apply $\HH$, set $y_n=y_{n+1}$, and simplify according to Lemma \ref{lemma:HHmarkov}:
\[
\HH\Big((C\sqcup \one_1)\otimes_R F(\sigma_n)\Big)[\yy,y_{n+1}]/(y_n-y_{n+1}) \ \ \cong \ \  \hskip2.5in
\]
\[
\hskip1.5in \left(
\begin{tikzpicture}[baseline=-.3em]
\tikzstyle{every node}=[font=\small]
\node (a) at (0,1) {$\HH(C)[x_{n+1},\yy](-1,0)[0]$};
\node at (0,0) {$\oplus$};
\node (b) at (0,-1) {$\HH(C)[x_{n+1},\yy](3,1)[0]$};
\node (c) at (6,1) {$\HH(C)[x_{n+1},\yy](1,0)[-1]$};
\node at (6,0) {$\oplus$};
\node (d) at (6,-1) {$\HH(C)[x_{n+1},\yy](3,1)[-1]$};
\path[->,>=stealth',shorten >=1pt,auto,node distance=1.8cm,
  thick]
(a) edge node {$x_n-x_{n+1}$} (c)
(b) edge node {$\Id$} (d);
\end{tikzpicture}
\right)
\]
The bottom row represents a contractible summand.  Expanding the other summand according to $\Q[x_{n+1}]\cong \Q\oplus \Q(-2)\oplus \Q(-4)\oplus \cdots $ gives
\[
\left(
\begin{tikzpicture}[baseline=-7em]
\tikzstyle{every node}=[font=\small]
\node (aa) at (0,0) {$\HH(C)[\yy](-1,0)[0]$};
\node (ab) at (4,0) {$\HH(C)[\yy](1,0)[-1]$};
\node (ba) at (0,-3) {$\HH(C)[\yy](-3,0)[0]$};
\node (bb) at (4,-3) {$\HH(C)[\yy](-1,0)[-1]$};
\node (ca) at (0,-6) {$\cdots $};
\node (cb) at (4,-6) {$\cdots$};
\path[->,>=stealth',shorten >=1pt,auto,node distance=1.8cm,
  thick]
(aa) edge node {$x_n$} (ab)
(aa) edge node {$-\Id$} (bb)
(ba) edge node {$x_n$} (bb)
(ba) edge node {$-\Id$} (cb);
\end{tikzpicture}
\right).
\]
After a Gaussian elimination (\S \ref{subsec:gauss}), this object is homotopy equivalent to $\HH(C)[\yy](1,0)[-1]$.  We have thus shown that
\[
\HH((C\sqcup \one_1)\otimes_R F(\sigma_n))[\yy,y_{n+1}]/(y_n-y_{n+1}) \simeq \HH(C)[\yy](1,0)[-1].
\]
This is an equivalence of $y$-ifications in $R_n[\yy]\gmod$.  In particular, taking the further quotient by $(y_1-y_{w(1)},\ldots,y_n-y_{w(n)})$ yields
\[
\CY((\CB\sqcup\one_1)\otimes \FY(\sigma_n)) \simeq \CY(\CB)(1,0)[-1],
\]
as claimed.  This proves (2); statement (3) is similar.
\end{proof}

\subsection{Homology with coefficients}
\label{subsec:homology with coefficients}


\begin{definition}\label{def:homologyWithCoeffs}
Let $\CB=(C,w,\Delta)\in \YC(R^e\gmod)$ be a $y$-ification.  Let $M$ be a module over $\Q[\yy]_w$.  We define the homology of $\CB$ \emph{ with coefficients in $M$} to be
\[
\HY(\CB;M):=H(\CY(\CB;M)),\qquad\qquad \CY(\CB;M):=\CY(\CB)\otimes_{\Q[\yy]_w}M.
\]
If $\CB=\FY(\b)$, then we also write $\CY(\b;M)=\CY(\FY(\b);M)$ and $\HY(\b;M):=H(\CY(\b;M))$. 
\end{definition}

Proposition \ref{prop:markov} states that $\CY(\b)$ depends only on the link $\hat\b$ up to overall shift and homotopy equivalence of differential trigraded $\Q[\yy]_L$.  This gives the following as an immediate consequence.
\begin{proposition}
The homology with coefficients $\HY(\b;M)$ depends only on the link $L=\hat\b$ up to isomorphism.\qed
\end{proposition}

\begin{example}\label{ex:Q0}
Let $\Q_0$ denote the module $M=\Q$ with each $y_c$ acting by zero, and let $\CB=(C,w,\Delta)$ be a $y$-ification. Then $\CY(\CB;\Q_0) = \HH(C)$.
\end{example}

\begin{example}\label{ex:Qnu}
Given a choice of scalars $\nu_c\in \Q$, indexed by cycles of $w$, let $\Q_\nu$ denote the module $M=\Q$ with each $y_c$ acting by $\nu_c$.  Then the differential on $\CY(\CB;\Q_\nu)$ is a sum:
\[
\sum_{k_1,\ldots,k_n} \HH(\Delta_{k_1,\ldots,k_n})\otimes \nu_1^{k_1}\cdots \nu_n^{k_n},
\]
where we write $\nu_i = \nu_c$ whenever $i$ is in cycle $c$.  Each $\Delta_{k_1,\ldots,k_n}$ has tridegree $(2\ell,0,1-2\ell)$, where $\ell=k_1+\cdots+k_n$.  Since the $\nu_i$ are rational numbers, their tridegree is $(0,0,0)$.  Thus  the differential on $\CY(\CB;\Q_\nu)$ is in general $\deg_Q$ non-decreasing, and is $(\deg_A,\deg_Q+\deg_T)$ preserving.  In particular $\HY(\CB;\Q_\nu)$ is bigraded via $(\deg_Q,\deg_A,\deg_T)$, and filtered with respect to $\deg_Q$.
\end{example}

\begin{remark}\label{rmk:familyOfHomology}
If $L=\hat\b$ is a link and $\pi_0(L)$ is its set of components, then we will regard $\HY(L;\Q_\nu)=\HY(\b;\Q_\nu)$ as a family of homology theories, parametrized by points $\nu\in \Q^{\pi_0(L)}$.  The fiber over zero, i.e.~$\HY(L;\Q_0)$, is the usual triply graded Khovanov-Rozansky homology, and is the unique fiber in which the trigrading does not collapse.
\end{remark}

Given remarks in Example \ref{ex:Qnu} we have a spectral sequence which computes $\HY(L;Q_\nu)$.

\begin{theorem}\label{thm:deformationSS}
There is a spectral sequence with $E_1$ page isomorphic to Khovanov-Rozansky homology $\HKR(F(\b))$ and $E_\infty$ page the associated graded of $\HY(\b;\Q_\nu)$ with respect to the $\deg_Q$-filtration.\qed
\end{theorem}
\begin{remark}\label{rmk:SSgradings}
The pages of this spectral sequence are triply graded, via $(\deg_Q,\deg_A,\deg_T)$.  The differential of the $r$-th page has degree $(2r,0,1-2r)$, hence is zero for $r\gg 0$ since, e.~g.~, the $E_1$ page $\HKR(F(\b))$ is supported in finitely many homological degrees.  In fact $d_r=0$ for $r>\ell/2$, where $\ell$ is the homological width of $\HKR(\b)$. 
\end{remark}

\begin{definition}\label{def:parity}
A complex $C\in \KC^b(R^e\gmod)$ is called \emph{parity} if $\HKR(C)$ is supported only in even or odd homological degrees.  A link $L$ is called \emph{parity} if the Rouquier complex $F(\b)$ is parity for some (hence every) braid representative $\b$. 
\end{definition}

\begin{corollary}\label{cor:parityAndFlatness}
Let $\CB=(C,w,\Delta)$ be a $y$-ification of $C\in \KC^b(R^e\gmod)$.  If $C$ is parity then $\HKR(C)\cong \HY(\CB;\Q_\nu)$ as bigraded vector spaces for all $\nu\in \Q^{\cyc{w}}$, after collapsing the trigrading to a bigrading as above.\qed
\end{corollary}

\begin{remark}
The isomorphism in Corollary \ref{cor:parityAndFlatness} should be regarded as an isomorphism of the associated graded spaces with respect to the $\deg_Q$ filtration.  In particular, it is not generally an isomorphism of $\Q[\xx]_w$-modules.  See the example below.
\end{remark}

\begin{example}
Recall the complex which computes the Hochschild degree zero part of the homology of the positive Hopf link:
\[
\begin{tikzpicture}[baseline=-.5em]
\tikzstyle{every node}=[font=\small]
\node (a) at (0,0) {$\underline{R[\yy](-2)}$};
\node (b) at (3,0) {$R[\yy]$};
\node (c) at (6,0) {$R[\yy](2)$};
\path[->,>=stealth',shorten >=1pt,auto,node distance=1.8cm,
  thick]
(b) edge node[above] {$y_1-y_2$}	(a)
(b) edge node[above] {$x_1-x_2$} (c);
\end{tikzpicture}
\]
If we specialize $y_i=\nu_i$ with $\nu_1=\nu_2\in \Q$, then the homology is isomorphic to $\Q[x_1,x_2]/(x_1-x_2)$, generated by a class $\b$ in tridegree $(-2,0,2)$, direct sum a copy of $\Q[x_1,x_2]$, generated by a class $\a$ in tridegree $(2,0,0)$.  After collapsing the tridegree $(\deg_Q,\deg_T,\deg_A)$ to $(\deg_A,\deg_Q+\deg_T)$, we find that
\[
\HKR(\sigma_1^2)=\HY(\sigma_1^2;\Q_{\nu})\cong \beta\Q[x_1,x_2]/(x_1-x_2) \oplus \Q[x_1,x_2]\alpha
\]
where $\deg(\alpha)=(0,2)$.

On the other hand, if we specialize to numbers $y_i=\nu_i$ with $\nu_1\neq \nu_2\in \Q$, then the first two terms of the complex cancel one another, and we see that
\[
\HY(\sigma_1^2;\Q_\nu)\cong \Q[x_1,x_2].
\]
As a bigraded vector space, this is isomorphic to $\HKR(\sigma_1^2)$, as claimed.
\end{example}

In \S \ref{subsec: negative} we consider examples in which the spectral sequence from Theorem \ref{thm:deformationSS} does not collapse immediately.

\subsection{Modules over the unlink}
\label{subsec:modules}
In what follows we will be interested in studying the homology of an $r$-component link $\HY(L)$ as a module over the homology of the $r$-component unlink $\Q[x_1,\ldots,x_r,y_1,\ldots,y_n,\theta_1,\ldots,\theta_n]$.  In order to bring the odd variables $\theta_i$ onto equal footing with the even variables, it is necessary to pass to derived categories.

\begin{definition}
Let $\DC_n:=D^b(R^e_n\gmod)$ denote the derived category of graded $(R_n,R_n)$-bimodules.  For each complex
\[
B = \begin{diagram}[small] \cdots &\rTo & B^{k} &\rTo^{d} & B^{k+1} &\rTo & \cdots\end{diagram}
\]
and each $i,j\in \Z$, let $B(i,j)$ denote the complex
\[
B(i,j) = \begin{diagram}[small] \cdots &\rTo & B^{k+j}(i) &\rTo & B^{k+j+1}(i) &\rTo & \cdots\end{diagram}.
\]
In other words, $B(i,j)^k = B^{k+j}(i)$.  
\end{definition}

\begin{remark}
We regard $\DC_n$ as a monoidal category with derived tensor product $\buildrel L\over{\otimes}_R$ and identity $\one=R$.  Soergel bimodules are free as left (and right) $R$-modules, hence the derived tensor product of Soergel bimodules coincides with the ordinary tensor product.  It follows that $\SBim_n$ is a full monoidal subcategory of $\derived_n$, hence $\KC(\SBim_n)$ and $\YC(\SBim_n)$ may be regarded as full monoidal subcategories of $\KC(\DC_n)$ and $\YC(\DC_n)$. 
\end{remark}

Since $\DC_n$ has grading shift functors $(i,j)$, we have the bigraded hom spaces $\Hom_{\DC_n}^{\Z\times \Z}(M,N)$ as in \S \ref{subsec:gradedCats}.
\begin{lemma}\label{lemma:HHasHom}
We have $\Hom_{\DC_n}^{\Z\times \Z}(\one, M)\cong \HH(M)$ for all graded $(R,R)$-bimodules $M$. In particular,
\[
\End_{\DC_n}^{\Z\times \Z}(\one)=\HH(R) = \Q[x_1,\ldots,x_n]\otimes_\Q \Lambda[\theta_1,\ldots,\theta_n].
\]
\end{lemma}
\begin{proof}
Given a graded $(R,R)$-bimodule $M$, the Hoschchild cohomology of $M$ is the direct sum of $\Ext^j_{R^e}(R,M(i))$ over all $i,j$.  But $\Ext$s are nothing other than homs in the derived category $\DC_n$ of graded $R^e$-modules.  This proves the first statement, given that $R=\one$ is the identity bimodule. 

For the second, we note that $\End_{\DC_n}^{\Z\times \Z}(R,R)$ can be computed by choosing an $R^e$-free resolution $K^\bullet\rightarrow R$ and using
\[
\End_{\DC_n}^{\Z\times \Z}(R)\cong H(\underline{\End}_{\DC_n}^{\Z\times \Z}(K,K)) \cong H(\underline{\End}_{\DC_n}^{\Z\times \Z}(K,R)),
\]
where the underlines mean complex of homs.  The algebra structure is most evident from the expression in the middle, but computations at the level of graded vector spaces are more accessible from the expression on the right.  In any case, for $K$ we have the well known Koszul resolution of $R$:
\[
K:= (\Q[x_1,x_1'](-2)\buildrel x_1-x_1'\over \longrightarrow \underline{\Q[x_1,x_1']})\otimes_\Q \cdots \otimes_\Q (\Q[x_n,x_n'](-2)\buildrel x_n-x_n'\over \longrightarrow \underline{\Q[x_n,x_n']}).
\]
Taking $\underline{\Hom}_{\DC_n}^{\Z\times\Z}(K,R)$ yields
\[
\underline{\Hom}_{\DC_n}^{\Z\times\Z}(K,R)  = (\underline{\Q[x_1]}\buildrel 0\over \rightarrow \Q[x_1](2))\otimes_\Q \cdots \otimes_\Q  (\underline{\Q[x_n]}\buildrel 0\over \rightarrow \Q[x_n](2)),
\]
whose homology (easy to compute, given that the differential is zero) is $\Q[x_1,\ldots,x_n]\otimes_\Q \L[\theta_1,\ldots,\theta_n]$, in which the $x_i$ have bidegree $(2,0)$ and the $\theta_i$ have bidegree $(-2,1)$.  The algebra structure is now characterized by the fact that $\one$, being the monoidal identity in a graded monoidal category $\DC_n$, has a graded commutative endomorphism ring.  Alternatively, one may regard $\theta_i$ as the endomorphism of Koszul complexes which is the identity on all but one factor, and on the remaining factor is
\[
\begin{diagram}
0 &\rTo & \Q[x_i,x_i'](-2) &\rTo^{x_i-x_i'} & \underline{\Q[x_i,x_i']}\\
\dTo &&\dTo^{\Id} &&\dTo \\
 \Q[x_i,x_i'](-4) &\rTo^{x-x'} & {\Q[x_i,x_i'](-2)} &\rTo & \underline{0}
\end{diagram}.
\] 
From this description it is clear that the composition $\theta_i(-2,1)\circ \theta: R\rightarrow R(-4,2)$ is zero.  Given the Koszul sign rule for defining the tensor product of morphsims with nonzero degree, it follows also that the $\theta_i$ anti-commute with one another, as claimed.
\end{proof}

\begin{definition}
A class $c\in \HH^{i,j}(M)$ will be said to have $\wt(c)=Q^iA^j$.  Henceforth, we let $\ttheta$ denote the set of variables $\theta_1,\ldots,\theta_n$.  These are odd variables with $\wt(\theta_i)=AQ^{-2}$.  Here ``odd'' means that the $\theta_i$ square to zero and anti-commute with other $\theta_j$'s.  Additionally, $\Q[\xx,\ttheta]$ will denote the \emph{super-polynomial ring}
\[
\Q[\xx,\ttheta]:=\Q[x_1,\ldots,x_n,\theta_1,\ldots,\theta_n]:=\Q[x_1,\ldots,x_n]\otimes_\Q \Lambda[\theta_1,\ldots,\theta_n].
\]
Similarly for $\Q[\xx,\yy,\ttheta]$.
\end{definition}

Now, given any complex $C^\bullet\in \KC(\SBim_n)$, we may include $\SBim_n$ as a full subcategory in $\DC_n$ and then regard $C^{\bullet}$ as an object in $\KC(\DC_n)$ (a complex of complexes).  It is important that we are forgetting the triangulated structure, and regarding $\DC_n$ only as an $R^e$-linear category with grading shift functors $(i,j)$.  However, $\KC(\DC_n)$ is a triangulated category via the usual mapping cone construction and the usual homological shift $[k]$ in $\KC$.
Unpacking definitions, we have
\[
\HKR(C) = \Hom_{\KC(\DC_n)}^{\Z\times \Z\times \Z}(R,C) := \bigoplus_{i,j,k\in \Z}\Hom_{\KC(\DC_n)}(R,C(i,j))[k].
\]

Let us discuss $y$-ifications in this context.  This makes sense since $\DC_n$ is linear with respect to $R^e$.  The category of $y$-ifications $\YC(\DC_n)$ is monoidal via the construction in \S \ref{subsec:tensorproduct}, and contains $\YC(\SBim_n)$ as a full monoidal subcategory; it has grading shift functors of the form $(i,j)[k]$ where $(i,j)$ is inherited from $\DC_n$ and $[k]$ is the homological shift in $\YC(\cdots)$.

\begin{lemma}
We have $\End_{\YC(\DC_n)}^{\Z\times \Z\times \Z}(\one_{\YC})\cong \Q[\xx,\yy,\ttheta]$.  Moreover, if $\CB=(C,w,\Delta)\in \YC(\DC_n)$ is a $y$-ification associated to the trivial permutation $w=1\in S_n$, then
\[
\Hom_{\YC(\DC_n)}^{\Z\times \Z\times \Z}(\one_{\YC},\CB) \cong \HY(\CB).
\]
\end{lemma}
Finding a similar representation of the $y$-ified homology when $w\neq 1$ is slightly more involved, and we return to this after the proof.
\begin{proof}
First, recall a special case of tensor-hom adjunction, which states that of $M$ is an $S$ module and $N$ is an $S[\yy]$-module, then
\[
\Hom_{S[\yy]}(M[\yy],N)\cong \Hom_S(M,N).
\]
An analogue of this fact gives us the the first in the following sequence of isomorphisms:
\begin{eqnarray*}
\underline{\Hom}^{\Z\times\Z\times \Z}_{\YC(\DC_n)}(R[\yy],C[\yy]) &\cong& \underline{\Hom}^{\Z\times\Z\times \Z}_{\seq{\DC_n}}(R,C[\yy])\\
& \cong & \underline{\Hom}^{\Z\times\Z\times \Z}_{\seq{\DC_n}}(R,C)\otimes_\Q \Q[\yy]\\
& \cong & \HH(C)\otimes \Q[\yy]\\
& = & \CY(C).
\end{eqnarray*}
The second isomorphism uses the fact that we only consider homogeneous morphisms, hence even though $\Q[\yy]$ is infinite dimensional as a vector space, it is finite dimensional in each degree, hence tensoring with $\Q[\yy]$ commutes with $\Hom(R,-)$ (the reader should contrast this with the case of ungraded vector spaces: if $V$ is countable dimensional then $\Hom(V,W\otimes_\Q \Q^\omega)$ has uncountable dimension, while $\Hom(V,W)\otimes_\Q  \Q^\w$ has countable dimension).  The third isomorphism follows from Lemma \ref{lemma:HHasHom}, and the last is by definition.  It is an exercise to follow the differential through this chain of isomorphisms.  Taking homology completes the proof.
\end{proof}

It follows that every $y$-ification $\CB=(C,w,\Delta)$ can be thought of as a bimodule over $\Q[\xx,\yy,\ttheta]$, where the left (resp.~right) action is induced from the action on the first factor of $\one_\YC\otimes \CB\cong\CB$ (resp.~second factor of $\CB\otimes\one_\YC\cong \CB$).   

\begin{proposition}\label{prop:HYmoduleStruct}
The $y$-ified homology $\HY(L)$ is a well-defined triply graded module over $\Q[\xx,\yy,\ttheta]_L:=\Q[x_c,y_c,\theta_c]_{c\in \pi_0(L)}$ up to isomorphism and overall shift.
\end{proposition}

This essentially reduces to a detailed computation in $\DC_2$, which we take up in the next section.  But first we discuss how to represent the $y$-ified homology as a hom space, when $w\neq 1$.

\begin{definition}
Let $w\in S_n$ be given.  Let $\pi_0(w)$ denote the set of cycles of $w$, which may be identified with equivalence classes in $\{1,\ldots,n\}$.   For each pair $i\neq j\in \{1,\ldots,n\}$, let $\KB_{i,j}(\one)\in \YC(\SBim_n)$ denote the \emph{Koszul $y$-ification}
\[
\KB_{i,j}(\one) \ := \ 
\left(
\begin{tikzpicture}[baseline=-.5em]
\tikzstyle{every node}=[font=\small]
\node (a) at (0,0) {$\Q[\xx,\yy](-2)$};
\node (b) at (4,0) {$\underline{\Q[\xx,\yy]}$};
\path[->,>=stealth',shorten >=1pt,auto,node distance=1.8cm,
  thick]
([yshift=3pt] a.east) edge node[above] {$x_{i}-x_{j}$}		([yshift=3pt] b.west)
([yshift=-2pt] b.west) edge node[below] {$-(y_i-y_{j})$}		([yshift=-2pt] a.east);
\end{tikzpicture}\right).
\]
If $w\in S_n$ is a permutation and $w=\tau_{i_1,j_1}\cdots\tau_{i_r,j_r}$ is a minimal length expression of $w$ as a product of transpositions, then we let $\KB_w(\one):=\KB_{i_1,j_1}(\one)\otimes\cdots\otimes \KB_{i_r,j_r}(\one)$ be the tensor product $y$-ification.
\end{definition}
Note that $\KB_w(\one)$ is a $y$-ification of the Koszul complex (in $\KC(\SBim_n)$) associated to the sequence $x_i-x_j$, indexed by pairs $(i,j)$ with $i=w(j)$.

\begin{proposition}
If $\CB=(C,w,\Delta)\in \YC(\DC_n)$ is a strict $y$-ification, then
\[
\Hom_{\YC(\DC_n)}^{\Z\times \Z\times \Z}(\KB_w(\one), \CB) \cong \HY(\CB).
\]
\end{proposition}
We omit the proof, since we will not need this result elsewhere.

\subsection{The odd variables} 
\label{subsec:oddvariables}

\begin{proposition}\label{prop:odd}
Inside $\DC_2$ we have
\begin{subequations}
\begin{equation}\label{eq:derivedHomRR}
\End^{\Z\times \Z}_{\DC_2}(R)\cong \Q[x_1,x_2,\theta_1,\theta_2]=:A.
\end{equation}
\begin{equation}\label{eq:derivedHomBR}
\Hom^{\Z\times \Z}_{\DC_2}(B,R)\cong A\{b,\omega\} / (x_1-x_2)\omega = (\theta_1-\theta_2)b
\end{equation}
\begin{equation}\label{eq:derivedHomRB}
\Hom^{\Z\times \Z}_{\DC_2}(R,B)\cong A\{b^\ast,\omega^\ast\} / (x_1-x_2)\omega^\ast = (\theta_1-\theta_2)b^\ast
\end{equation}
\end{subequations}
where the degrees are $\deg(\theta_i)=(-2,1)$, $\deg(b)=\deg(b^\ast)=(1,0)$, and $\deg(\omega)=\deg(\omega^\ast)=(-3,1)$.  The isomorphisms \eqref{eq:derivedHomBR} and \eqref{eq:derivedHomRB} are isomorphisms of bigraded $A$-modules.  These morphisms satisfy the following relations
\begin{subequations}
\begin{equation}
b\circ \omega^\ast = \theta_1-\theta_2 = \omega\circ b^\ast, \qquad\qquad b\circ b^\ast = x_1-x_2 \in \End^{\Z\times \Z}(R)
\end{equation}
\begin{equation}
b^\ast\circ \omega = \theta_1\otimes \Id - \Id\otimes\theta_2  = -(\theta_2\otimes \Id - \Id\otimes \theta_1) = \omega^{\ast}\circ b \in \End^{\Z\times \Z}_{\DC_2}(B)
\end{equation}
\end{subequations}
In fact, post composition with $b$ identifies $\Hom^{\Z\times \Z}_{\DC_2}(R,B)$ with the ideal in $A$ generated by $b\circ b^\ast = x_1-x_2$ and $b\circ \omega^\ast= \theta_1-\theta_2$.
\end{proposition}
The proof of this proposition occupies the remainder of this subsection.

\begin{definition}
In this subsection, let $\d_i:=x_i-x_i'\in R^e$.  Let $\xi_i$ and $\phi_i$ denote odd variables ($i=1,2$) with $\wt(\xi_i)=\wt(\phi_1)=A\inv Q^{2}$ and $\wt(\phi_2)=A\inv Q^{4}$.

Let $K_R=R^e[\xi_1,\xi_2]$ with $R^e$-linear differential determined by $d(\xi_i)=\d_i$ together with the graded Leibniz rule.  In other words,
\[
d \ : \ (1,\xi_1,\xi_2,\xi_1\xi_2) \ \mapsto (0, \d_1, \d_2, \d_1 \xi_2 - \d_2\xi_1)
\]
Let $K_B=R^e[\phi_1,\phi_2]$ with $R^e$-linear differential determined by $d(\phi_1)=\d_1+\d_2$, $d(\phi_2)=x_1\d_2 + x_2'\d_1$ together with the graded Leibniz rule.
\end{definition}

Then $K_R$ is isomorphic to $R$, and $K_B$ is isomorphic to $B_1(-1)$, in $\DC_2$.  In fact, $K_R$ and $K_B$ are resolutions of $R$ and $B_1(-1)$ by free graded $R^e$-modules.  In particular $x_i-x_i'$ are homotopic to zero as endomorphisms of $K_R$ and $x_1+x_2-x_1'-x_2'$ and $x_1x_2-x_1'x_2'$ are homotopic to zero as endomorphisms of $K_B$.   The complexes $K_R$ and $K_B$ are special cases of Koszul complexes, hence come with some naturally defined endomorphisms given by ``contracting'' with the generating odd variables.  We describe these next. 

For $R$ we have endomorphisms $\iota_{\xi_i}:R\rightarrow R(-2,1)$, $i=1,2$ defined on the level of resolutions $R\cong K_R$ by
\[
\iota_{\xi_1} \ : \ (1,\xi_1,\xi_2,\xi_1\xi_2)\mapsto (0, 1,0,\xi_2),
\]
\[
\iota_{\xi_2} \ : \ (1,\xi_1,\xi_2,\xi_1\xi_2)\mapsto (0, 0, 1, - \xi_1).
\]

For $B$ we have endomorphisms $\iota_{\phi_1}:B\rightarrow B(-2,1)$ and $\iota_{\phi_2}:B\rightarrow B(-4,1)$ defined on the level of resolutions $B\cong K_B(1)$ by
\[
\iota_{\phi_1} \ :  \ (1,\phi_1,\phi_2,\phi_1\phi_2) \ \mapsto \ (0,1,0,\phi_2),
\]
\[
\iota_{\phi_2} \ :  \ (1,\phi_1,\phi_2,\phi_1\phi_2) \ \mapsto \ (0,0,1,-\phi_1).
\]

\begin{definition}
Let $\theta_i\in \End_{\DC_2}^{\Z\times \Z}(R)$ denote the endomorphisms of degree $(-2,1)$ induced by $\iota_{\xi_i}$.
\end{definition}

Since $K_R$ is a projective resolution of the identity bimodule, we have $K_R\otimes_R K_B\simeq K_B\simeq K_B\otimes_R K_R$.  These equivalences give $K_B$ the structure of a bimodule over $\End^{\Z\times \Z}_{\DC_2}(K_R)$, with actions denoted by $\gamma \otimes \kappa$ and $\kappa \otimes \gamma$, for $\gamma\in \Endg(K_R)$ and $\kappa\in \Endg(K_B)$.  A precise understanding of these actions requires the construction of explicit homotopy equivalences $K_R\otimes_R K_B\simeq K_B\simeq K_B\otimes_R K_R$.  We focus mainly on the left action, since the right action is similar.

We have a chain map $\pi:K_R\otimes_R K_B \rightarrow K_B$ sending
\[
f(\xx,\xx')\otimes g(\xx,\xx')\mapsto f(\xx,\xx)g(\xx,\xx'), \qquad \xi_i\otimes 1\mapsto 0, \qquad 1\otimes \phi_i\mapsto \phi_i,
\]
extended multiplicatively to all of $K_R\otimes_R K_B$.  The inverse homotopy equivalence $\sigma:K_B\rightarrow K_R\otimes_R K_B$ is defined by
\[
x_i\mapsto x_i\otimes 1,\qquad x_i'\mapsto 1\otimes x_i',\qquad \phi_1\mapsto 1\otimes \phi_1 + \xi_1\otimes 1+ \xi_2\otimes 1,  \qquad \phi_2\mapsto 1\otimes \phi_2 + x_1\xi_2\otimes 1 +x_2'\xi_1\otimes 1.
\]

The composition $K_B\rightarrow K_R\otimes_R K_B\rightarrow K_B$ is the identity, and the composition $K_R\otimes_R K_B\rightarrow K_B \rightarrow K_R\otimes_R K_B$ is homotopic to the identity, with the homotopy defined by
\[
h(1\otimes 1) = h(x_i\otimes 1)=h(1\otimes x_i')=h(\xi_i\otimes 1) = h(1\otimes \phi_i)=0, \qquad h(x_i'\otimes 1) = h(1\otimes x_i) = \xi_i\otimes 1,
\]
extended via the graded Leibniz rule to all of $K_R\otimes_R K_B$.  We leave it as an exercise to verify this; in fact, $(\pi,\sigma,h)$ are the data of a \emph{strong deformation retract} $K_R\otimes_R K_B \rightarrow K_B$.

Using the above homotopy equivalences, we obtain an explicit action of $\End(K_R)$ on $K_B$.  The compositions $K_B\simeq K_R\otimes_R K_B \buildrel \theta_i\otimes \Id \over\longrightarrow  K_R\otimes_R K_B(-2,1) \simeq K_B$ satisfy
\[
\theta_1 \otimes \Id_{K_B} \ :\  (1,\phi_1,\phi_2,\phi_1\phi_2) \ \mapsto \ (0, 1, x_2, \phi_2 - x_2\phi_1),
\]
\[
\theta_2 \otimes \Id_{K_B} \ :\  (1,\phi_1,\phi_2,\phi_1\phi_2) \ \mapsto \ (0, 1, x_1, \phi_2 - x_1\phi_1),
\]
extended by $R^e$-linearity.

Similar formulae produce inverse equivalences $K_B\simeq K_B\otimes_R K_R$, with corresponding right action of $\End(K_R)$ determined by
\[
\Id_{K_B}\otimes \theta_1 \ :\ (1,\phi_1,\phi_2,\phi_1\phi_2) \ \mapsto \ (0, 1, x_2', \phi_2 - x_2'\phi_1),
\]
\[
\Id_{K_B}\otimes \theta_2\  :\ (1,\phi_1,\phi_2,\phi_1\phi_2) \ \mapsto \ (0, 1, x_1', \phi_2 - x_1'\phi_1).
\]
We have proven the following.

\begin{lemma}\label{lemma:leftrightderivedaction}
The left and right module $\End^{\Z\times \Z}_{\DC_2}(R)$ structures on $\End^{\Z\times \Z}_{\DC_2}(B)$ satisfy:
\[
\theta \otimes \Id_{B} = \iota_{\phi_1} -x_2 \iota_{\phi_2}  \qquad\qquad \theta_2 \otimes \Id_{K_B} = \iota_{\phi_1} -x_1 \iota_{\phi_2}
\]
\[
\Id_{B}\otimes\theta_1 = \iota_{\phi_1}-x_2' \iota_{\phi_2}\qquad\qquad \Id_{B} \otimes \theta_2  = \iota_{\phi_1}-x_1' \iota_{\phi_2}.
\]\qed
\end{lemma}
From this lemma, it is clear that the left and right actions of $\theta_1+\theta_2$ on $B$ are equal in $\DC_2$.

The dot maps $b:B\rightarrow R(1)$ and $b^\ast:R(-1)\rightarrow B$ are represented on the level of resolutions as the chain maps $K_B\rightarrow K_R$ and $K_R\rightarrow K_B(2)$ defined by
\[
b : (1,\phi_1,\phi_2,\phi_1\phi_2) \mapsto (1,\xi_1+\xi_2,x_1\xi_2+x_2'\xi_1, (x_1-x_2')\xi_1\xi_2),
\]
\[
b^\ast : (1,\xi_1,\xi_2,\xi_1\xi_2) \mapsto (x_1-x_2',x_1\phi_1-\phi_2, -x_2'\phi_1+\phi_2, \phi_1\phi_2).
\]

Finally, we define the special morphisms  $\w :B\rightarrow R(-3,1)$ and $\w^\ast :R\rightarrow B(-3,1)$ on the level of resolutions by
\[
\w \ : \ (1,\phi_1,\phi_2,\phi_1\phi_2) \ \mapsto \ (0,0, -1,\xi_1+\xi_2),
\]
\[
\w^\ast \ : \ (1,\xi_1,\xi_2,\xi_1\xi_2)\ \mapsto \ (0, 1, -1, \phi_1).
\]
Given these explicit formulae, it is straightforward to check all of the relations in Proposition \ref{prop:odd}.  Now, Lemma \ref{lemma:HHmarkov} implies that, as bigraded $\Q[x_1,x_2]$-modules, $\HH(B)=\Hom^{\Z\times \Z}_{\DC_2}(R,B)$ is free of graded rank $(Q+AQ^{-3})(1+AQ^{-2})$.  A dimension count shows that post-composing with $b$ embeds $\HH(B)$ into $\HH(R)$, with image as asserted.  We leave the remaining details as exercises.\qed

\begin{corollary}
Regard the complexes $F(\sigma_1^{\pm})$ as objects in $\KC^b(\SBim_2)\subset \KC^b(\DC_2)$.  Then $\theta_1\otimes \Id - \Id\otimes \theta_2$ and $\theta_2\otimes \Id- \Id\otimes \theta_1$ are null-homotopic as morphisms $F(\sigma_1^{\pm})\rightarrow F(\sigma_1^\pm)(-2,1)$.
\end{corollary}
\begin{proof}
The following diagrams define the required homotopies.
\begin{equation}\label{eq:oddhomotopies}
\begin{diagram}
B & \rTo^{b} & R(1,0)\\
&\ldTo^{\w^\ast}&\\
B(-2,1) & \rTo^{b} & R(-1,1)
\end{diagram}
\qquad\qquad
\begin{diagram}
R(-1,0) & \rTo^{b^{\ast}} & B(0,0)\\
&\ldTo^{\w}&\\
R(-3,1) & \rTo^{b^{\ast}} & B(-2,1).
\end{diagram}
\end{equation}
\end{proof}

\begin{proof}[Proof of Proposition \ref{prop:HYmoduleStruct}]
Consider the triply graded dg algebra $E$ of $R^e$-linear endomorphisms of the Rouquier complex $F(\sigma_1^{\pm})\in \KC^b(\DC_2)$.  Inside $E$ we have homotopies $h_{x}$ and $h_{\theta}$ with degrees $(2,0,-1)$ and $(-2,1,-1)$, satisfying
\[
d(h_x)= x_1-x_2' = -(x_2-x_1'),\qquad\qquad d(h_\theta)= \theta_1-\theta_2' = -(\theta_2-\theta_1').
\]
These homotopies $h_\theta$ are defined by \eqref{eq:oddhomotopies} and, $h_x$ is given by
\[
\begin{diagram}
B & \rTo^{b} & R(1)\\
&\ldTo^{b^\ast}&\\
B(2) & \rTo^{b} & R(3)
\end{diagram}
\qquad\qquad
\begin{diagram}
R(-1) & \rTo^{b^{\ast}} & B\\
&\ldTo^{b}&\\
R(1) & \rTo^{b^{\ast}} & B(2).
\end{diagram}.
\]
These homotopies square to zero and anti-commute with one another.  Now, let $\b$ is any braid, and let $D$ be the link diagram associated to $\hat\b$.  If $p,q$ are two points which are on opposite sides of a crossing, then we have homotopies $h^{(x)}_{p,q}$ and $h^{(\theta)}_{p,q}$ of $F(\b)$, given by the appropriate homotopy $h_x$ or $h_\theta$ acting on the appropriate tensor factor $F(\sigma_i^{\pm})$.

These homotopies anti-commute with one-another and square to zero.  By construction, the strict $y$-ification $\CY(\b)$ can be described as
\[
\CY(\b):=\HH(F(\b))\otimes_\Q \Q[\yy]_L
\]
with differential $\Delta = d_C + \Delta_1$, where $\Delta_1$ is linear in $\yy$ and is a sum of terms $h_{p,q}\otimes y_c$ where $c$ is a component of the link $\hat\b$ and $p,q$ run over a set of pairs of points on $c$, lying on opposite sides of some crossing.  From the aforementioned anti-commuting property, it follows that $h^{(x)}_{p,q}$ and $h^{(\theta)}_{p,q}$ also define homotopies on the $y$-ification $\CY(\b)$, which slide $x$ and $\theta$ past crossings.
\end{proof}

\section{Link splitting and flatness}
\label{sec:splitting}
Throughout this section we abbreviate by writing
\[
\yring{L}:=\Q[y_c]_{c\in \pi_0(L)},\qquad \Q[\xx,\yy,\boldsymbol{\theta}]_L:=\Q[x_c,y_c,\theta_c]_{c\in \pi_0(L)}.
\]

\subsection{Split union}

\begin{definition}\label{def:split}
If $L,L'\subset \R^3$ are links, their \emph{split union} is the link $L\sqcup L'$ in which $L$ and $L'$ are separated by a 2-sphere. Given any link $L$, we let $\split(L)$ denote the split union of the components of $L$.  A link is \emph{split} if it can be written as the split union of two nontrivial links, and a link is \emph{totally split} if it is the split union of a finite collection of knots, i.e.~$L=\split(L)$.
\end{definition}

Let $\pi_0(L)$ denote the set of components of $L$, we assume that we are given a natural bijection between $\pi_0(L)$ and $\pi_0(\split(L))$.  More generally, if two links $L$ and $L'$ are related by a sequence of crossing changes which  only involves crossing changes between different components, we can identify the components of $L$ with the ones of $L'$, and define a bijection between $\pi_0(L)$ and $\pi_0(L')$.   This bijection defines a canonical isomorphism $\yring{L}\simeq \yring{\split(L)}$.

The following describes the behavior of $\HY(L)$ with respect to split union.

\begin{proposition}\label{prop:splittensor}
If $L\sqcup L'$ is the split union of $L$ and $L'$, then
\[
\HY(L\sqcup L')\cong\HY(L)\otimes_\Q \HY(L')\qquad\text{and}\qquad \HY(L\sqcup L';\Q_{\nu,\nu'})\cong\HY(L;\Q_\nu)\otimes_\Q \HY(L';\Q_{\nu'})
\]
for all $\nu\in \Q^{\pi_0(L)}$, $\nu'\in \Q^{\pi_0(L')}$.
\end{proposition} 
\begin{proof}
We may represent $L\sqcup L'$ by a ``split union'' of braids $\b\sqcup \b'$ where $\hat\b=L$ and $\hat\b'=L'$. On the level of complexes it is clear that $\CY(\b\sqcup \b')\cong \CY(\b)\otimes_\Q\CY(\b')$.  The K\"unneth theorem (for complexes over $\Q$) gives us an isomorphism in homology $\HY(\b\sqcup \b)\cong \HY(\b)\otimes_\Q\HY(\b')$, which can be checked is equivariant with respect to $\Q[\xx,\yy,\ttheta]_L$.  This proves the first statement.

For the second, we need only observe that $\CY(\b\sqcup \b';\Q_{\nu,\nu'})\cong \CY(\b;\Q_\nu)\otimes_\Q\CY(\b';\Q_{\nu'})$ and then use the K\"unneth theorem over $\Q$.
\end{proof}



\subsection{Link splitting maps}
\label{subsec:linksplitting}
An astute reader may have noticed the similarity between our definitions of $\FY(\sigma_i)$ and $\FY(\sigma_i\inv)$.  Specifically, if $y_i-y_{i+1}$ is invertible, then the two complexes are isomorphic.  In fact, the following diagram defines a morphism $\FY(\sigma_i)\rightarrow \FY(\sigma_i\inv)$:
\begin{equation}\label{eq:crossingchange}
\begin{tikzpicture}[baseline=-.5em]
\tikzstyle{every node}=[font=\small]
\node (a) at (0,0) {$\underline{B_i}[\yy]$};
\node (b) at (4.5,0) {$R[\yy](1)$};
\node (c) at (0,-3) {$R[\yy](-1)$};
\node (d) at (4.5,-3) {$\underline{B_i}[\yy]$};
\path[->,>=stealth',shorten >=1pt,auto,node distance=1.8cm,
  thick]
([yshift=3pt] a.east) edge node[above] {$b\otimes 1$}		([yshift=3pt] b.west)
([yshift=-2pt] b.west) edge node[below] {$-b^{\ast}\otimes (y_i-y_{i+1})$}		([yshift=-2pt] a.east)
([yshift=3pt] c.east) edge node[above] {$b^\ast \otimes 1$}		([yshift=3pt] d.west)
([yshift=-2pt] d.west) edge node[below] {$-b \otimes (y_i-y_{i+1})$}		([yshift=-2pt] c.east)
(a) edge node[xshift=1cm,yshift=-.6cm] {$\Id\otimes 1$} (d);
\draw[frontline,->,>=stealth',shorten >=1pt,auto,node distance=1.8cm,thick]
(b) to node[xshift=-3.7cm,yshift=0cm] {$-\Id\otimes (y_i-y_{i+1})$} (c);
\end{tikzpicture}
\end{equation}
After a Gaussian elimination (\S \ref{subsec:gauss}), the cone of this morphism is
\begin{equation}\label{eq:splitmapcone}
\Cone(\psi_i) \ \ \simeq  \ \ \begin{tikzpicture}[baseline=-.5em]
\tikzstyle{every node}=[font=\small]
\node (a) at (0,0) {$R[\yy](-1)$};
\node (b) at (4,0) {$\underline{R}[\yy](1)$};
\path[->,>=stealth',shorten >=1pt,auto,node distance=1.8cm,
  thick]
([yshift=3pt] a.east) edge node[above] {$(x_i-x_{i+1})\otimes 1$}		([yshift=3pt] b.west)
([yshift=-2pt] b.west) edge node[below] {$-\Id\otimes (y_i-y_{i+1})$}		([yshift=-2pt] a.east);
\end{tikzpicture}
\end{equation}
If we invert $(y_i-y_{i+1})$, then this complex becomes contractible, hence $\psi_i:\FY(\sigma_i)\rightarrow \FY(\sigma_i\inv)$ becomes an equivalence.  We also define another morphism $\FY(\sigma_i\inv)\rightarrow \FY(\sigma_i)$ of tridegree $(-2,0,2)$:

\begin{equation}\label{eq:crossingchange2}
\begin{tikzpicture}[baseline=-.5em]
\tikzstyle{every node}=[font=\small]
\node (a) at (0,0) {$\underline{B_i}[\yy]$};
\node (b) at (4.5,0) {$R[\yy](1)$};
\node (c) at (0,-3) {$R[\yy](-1)$};
\node (d) at (4.5,-3) {$\underline{B_i}[\yy]$};
\path[->,>=stealth',shorten >=1pt,auto,node distance=1.8cm,
  thick]
([yshift=3pt] a.east) edge node[above] {$b\otimes 1$}		([yshift=3pt] b.west)
([yshift=-2pt] b.west) edge node[below] {$-b^{\ast}\otimes (y_i-y_{i+1})$}		([yshift=-2pt] a.east)
([yshift=3pt] c.east) edge node[above] {$b^\ast \otimes 1$}		([yshift=3pt] d.west)
([yshift=-2pt] d.west) edge node[below] {$-b \otimes (y_i-y_{i+1})$}		([yshift=-2pt] c.east)
(d) edge node[xshift=3cm,yshift=0.2cm] {$\Id\otimes (y_i-y_{i+1})$} (a);
\draw[frontline,<-,>=stealth',shorten >=1pt,auto,node distance=1.8cm,thick]
(b) to node[xshift=-2.0cm,yshift=0cm] {$-\Id\otimes 1$} (c);
\end{tikzpicture}
\end{equation}

This construction yields the following.

\begin{proposition}\label{prop:yequivalence}  Let $\b^+$ and $\b^-$ be braids which differ in a single crossing:
\[
\b^\pm = \b'\sigma_i^\pm\b''
\]
for some braids $\b',\b''\in \Br_n$, and some $i$.  There exists a chain map $\psi^{+}_{\b',i,\b''}:\FY(\b^+)\rightarrow \FY(\b^-)$ and a chain map $\psi^{-}_{\b',i,\b''}:\FY(\b^-)\rightarrow \FY(\b^+)(-2)[2]$ such that
\[
\psi^{+}_{\b',i,\b''}\psi^{-}_{\b',i,\b''}=\psi^{-}_{\b',i,\b''}\psi^{+}_{\b',i,\b''}=\Id\otimes (y_i-y_{i+1}).
\]
 If the two strands involved in the $\sigma_i^\pm$ lie on components $c\neq c'$ then $\psi^{\pm}_{\b',i,\b''}$ become homotopy equivalences after inverting $(y_c-y_c')$. \qed
\end{proposition}

These splitting maps $\psi^{\pm}$ allow us to discuss the homology of $L$ ``in reference'' to the the homology of $\split(L)$.  Surprisingly, in many cases the $y$-ified homology of $L$ embeds in the $y$-ified homology of $\split(L)$ (see Theorem \ref{thm:splittableInjectivity}).  How this works for the positive Hopf link is illustrated in \S \ref{subsec:FT2ideal}.  We return to this point in \S \ref{sec:ideals}.

We now arrive at our link splitting property.

\begin{corollary}\label{cor:linksplittingproperty}
Let $L=L_1\cup\cdots \cup L_r$ be a link with components $L_i$, and let $M$ be a module over $\yring{L}$ such that $y_i-y_1$ acts invertibly on $M$ for all $i\neq 1$.  Then
\[
\HY(L;M)\cong \HY(L_1;M)\otimes_\Q \HY(L_2\cup\cdots\cup L_r;M)
\]
as (bigraded or triply graded) $\Q[\xx,\yy,\boldsymbol{\theta}]_L$-modules.
\end{corollary}
\begin{proof}
By switching crossings between $L_1$ and the other components, we may unlink $L_1$ from the rest of the link $L_2\cup\cdots\cup L_r$.  This gives us a link splitting map $\CY(L)\rightarrow \CY(L_1\sqcup (L_2\cup\cdots\cup L_r))$, where $\sqcup$ here denotes split union.  
\end{proof}
\begin{remark}
For instance, we could take $M=\yring{L}[(y_c-y_d)\inv]_{c\neq d}$.  Localizing in this way involves no collapse of degrees, since we are inverting homogeneous elements.  Thus, for this choice of $M$, $\HY(L;M)$ is a triply graded $\Q[\xx,\yy,\boldsymbol{\theta}]$-module.  For other choices of $M$ the triple grading on $\CY(\b;M)$ collapses to a bigrading upon taking homology, as in Example \ref{ex:Qnu}.
\end{remark}
\begin{remark}
Combining this with the spectral sequence from Theorem \ref{thm:deformationSS}, we see that there is a spectral sequence with $E_2$ page $\HKR(L)$ and $E_\infty$ page $\HKR(\split(L))$ with the collapsed gradings $(\deg_A,\deg_Q+\deg_T)$.
\end{remark}

Corollary \ref{cor:linksplittingproperty} allows us to reinterpret the homology with coefficients $\HY(L,\Q_{\nu})$ in a more intuitive way.

\begin{corollary}
A point $\nu\in \Q^{\pi_0(L)}$ defines a set partition 
$\pi_0(L)=\sqcup_{a} \Pi_a$ such that $\nu_c=\nu_{c'}$ if and only if $c,c'$ belong to the same block $\Pi_a$ for some $a$.
Define $L(\Pi_a)=\cup_{c\in \Pi_a}L_c$. Then
\[
\HY(L,\Q_{\nu})\cong \bigotimes_{a}\HKR(L(\Pi_a))
\]
with the collapsed gradings. 
\end{corollary}

\begin{proof}
By construction, $c$ and $c'$ are in different blocks then $\nu_c-\nu_{c'}$ is invertible and by Proposition \ref{prop:yequivalence}
we can arbitrarily change the crossings between $L_c$ and $L_{c'}$. This implies
\[
\HY(L,\Q_{\nu})\cong \HY(\sqcup_{a}L(\Pi_a),\Q_{\nu})=\bigotimes_{a}\HY(L(\Pi_a),\Q_{\nu}).
\]
Since $\nu_c=\nu_{c'}$ for all $c,c'\in \Pi_a$, one has $\HY(L(\Pi_A),\Q_{\nu})\cong \HKR(L(\Pi_a)).$
\end{proof}

\begin{example}
\label{ex:HY sheaf}
Suppose that $L$ has 3 components. The space $\Q^3$ can be stratified as follows: there are three planes $\{\nu_1=\nu_2\}$,
$\{\nu_1=\nu_3\}$ and $\{\nu_2=\nu_3\}$ which intersect in a line $\{\nu_1=\nu_2=\nu_3\}$. If $\nu$ is on the latter line,
$\HY(L,\Q_{\nu})\cong \HKR(L)$. If $\nu$ belongs to the plane $\{\nu_1=\nu_2\}$ (outside of the line)
then $\HY(L,\Q_{\nu})\cong \HKR(L_1\cup L_2)\otimes \HKR(L_3)$. If $\nu$ is generic (outside of the planes) then
\[
\HY(L,\Q_{\nu})\cong \HKR(L_1)\otimes \HKR(L_2)\otimes \HKR(L_3)\cong \HY(\split(L),\Q_{\nu}).
\]
\end{example}

\begin{definition}
A link $L$ is said to be \emph{parity} if $\HKR(L)$ is supported in only even homological degrees.  Similarly, a braid $\b$ is \emph{parity} if $L=\hat\b$ is parity.
\end{definition}

For example positive torus links are parity by work of Elias, Hogancamp, and Mellit \cite{EH,Hog17a,Mellit}.

\begin{corollary}\label{cor:paritysplitting}
If $\HKR(L)$ is parity, then $\HKR(L)\cong \HKR(\split(L))$ as bigraded vector spaces.  In particular
\[
\PC_L(q,t,a)|_{t=1} = \PC_{\split(L)}(q,t,a)_{t=1}.
\]
\end{corollary}
\begin{proof}
Choose a braid representative $\b\in \Br_n$ of $L$.  Let $\nu_1,\ldots,\nu_n\in \Q$ be such that $\nu_i-\nu_j$ is zero if $i$ and $j$ are on the same component of $L$ and invertible otherwise (here we are identifying the integers $1,\ldots,n$ with the points along the bottom of the braid diagram).  Assume that $L$ is parity.  Then the spectral sequence from Theorem \ref{thm:deformationSS} collapses at the $E_1$ page.  Indeed, the differential $d_r$ has homological degree $1-2r$, hence is zero for $r\geq 1$ by the parity condition. Thus $E_1=E_2=\cdots=E_\infty$.  The $E_1$ page is $\HKR(L)$, and the $E_\infty$ page is the associated graded of $\HY(L;\Q_\nu)$ with respect to the $\deg_Q$-filtration which, if we forget the $\Q[\xx]$-action and think of $\HY(L;\Q_\nu)$ only as a bigraded vector space, is $E_\infty\cong \HY(L;\Q_{\nu})$.  This shows that $\HKR(L)\cong \HY(L;\Q_\nu)$ as bigraded vector spaces.  On the other hand, $\HY(L;\Q_\nu)\cong \HY(\split(L);\Q_\nu)$ by Corollary \ref{cor:linksplittingproperty}, which is isomorphic to $\HKR(\split(L))$ by the same argument.  This completes the proof of the first statement.

For the final statement, recall that $q=Q^2$, and $t=T^2Q^{-2}$.  Thus collapsing the trigrading $(\deg_Q,\deg_A,\deg_T)$ to a bigrading $(\deg_A,\deg_Q+\deg_T)$ corresponds to setting $t=1$ in the Poincar\'e series.
\end{proof}

\begin{remark}
This isomorphism does not respect the $\Q[\xx]$-action.  For instance if $L=\widehat{\sigma_1^2}$ is the positive Hopf link and $R=\Q[x_1,x_2]$, then the Hochschild degree zero part
\[
\HKR^0(L)\cong R(-2)\oplus R/(x_1-x_2)(2)[-2]
\]
has $R$-torsion, but $\HKR^0(\split(L))=R$ has no $R$-torsion.  Nonetheless, $\HKR^0(L)\cong R$ as vector spaces, graded via $\deg_Q+\deg_T$.
\end{remark}

In what follows we will need an explicit map between the $y$-ified complexes for $L$ and for $\split(L)$.  

\begin{proposition}
\label{prop: crossing change maps}
Suppose that a link $L'$ is obtained from $L$ by a sequence of crossing changes between different components which involves $p_{c,c'}$ (resp. $n_{c,c'}$) positive to negative (resp. negative to positive) crossing changes between the components $L_c$ and $L_{c'}$. Let $p=\sum_{c,c'}p_{c,c'}$ and $n=\sum_{c,c'}n_{c,c'}$. Then there is a pair of chain maps 
$\Psi_{L\to L'}:\CY(L)\to \CY(L'), \Psi_{L'\to L}:\CY(L')\to \CY(L)$ of tridegrees $(-2n,0,2n)$ and $(-2p,0,2p)$ such that
\[
\Psi_{L\to L'}\Psi_{L'\to L}=\prod_{c,c'}(y_c-y_{c'})^{p_{c,c'}+n_{c,c'}}\cdot \Id_{\CY(L')},
\]
\[
\Psi_{L'\to L}\Psi_{L\to L'}=\prod_{c,c'}(y_c-y_{c'})^{p_{c,c'}+n_{c,c'}}\cdot \Id_{\CY(L)}.
\]
\end{proposition}

\begin{proof}
We can define $\Psi_{L\to L'}$ as a composition of maps $\psi^{\pm}$ from Proposition \ref{prop:yequivalence} associated to each crossing change. The desired equation immediately follows  from Proposition \ref{prop:yequivalence}. 
\end{proof}

\begin{remark}
The map $\Psi_{L\to L'}$ actually depends on the sequence of crossing changes and not only on $L$ and $L'$.
\end{remark}

\begin{corollary}
\label{cor:spliiting map}
Given a link $L$ and a sequence of crossing changes between different components which transforms it to $\split(L)$, one can define a map
\[
\Psi_{L\to \split(L)}:\CY(L)\to \CY(\split(L)).
\]
We will call it a {\em splitting map} and simply denote by $\Psi$ if a link $L$ and a sequence of crossing changes is clear from the context.
\end{corollary}

\subsection{Full twist on two strands}
\label{subsec:FT2ideal}
Consider the following diagram:
\begin{equation}\label{eq:FT2map}
\begin{tikzpicture}[baseline=-.5em]
\tikzstyle{every node}=[font=\small]
\node (e) at (-1,3) {$R[y_1,y_2](-2)$};
\node (f) at (9,3) {$R[y_1,y_2](2)$};
\node (a) at (-1,0) {$\underline{B}_1[y_1,y_2](-1)$};
\node (b) at (4,0) {${B}_1[y_1,y_2](1)$};
\node (c) at (9,0) {${R}[y_1,y_2](2)$};
\node (d) at (4,-3) {$R[y_1,y_2]$};
\path[->,>=stealth',shorten >=1pt,auto,node distance=1.8cm]
([yshift=3pt] a.east) edge node[above] {$(x_1-x_1')\otimes 1$}		([yshift=3pt] b.west)
([yshift=-2pt] b.west) edge node[below] {$\Id\otimes (y_1-y_2)$}		([yshift=-2pt] a.east)
(b) edge node[above] {$b\otimes 1$} (c)
(a) edge node {$b^\ast\otimes 1$} (d)
(c) edge node {$\Id\otimes(y_1-y_2)$} (d)
(e) edge node {$b\otimes 1$} (a)
(f) edge node {$\Id\otimes 1$}  (c) ;
\end{tikzpicture}.
\end{equation}

The middle row is $\FY(\sigma_1^2)$, up to equivalence (see \S \ref{subsec:FT2}), and the lower half of the diagram is an illustration of the canonical map $\psi:\FY(\sigma_1^2)\rightarrow R[\yy]$.  The two arrows in the top half of the diagram describe morphisms $\alpha_1:R[\yy](-2)\rightarrow \FT_2$ and $\alpha_2:R[\yy](2)[-2]\rightarrow \FT_2$.

Taking homology, we see that $\a_1$ and $\a_2$ are the generators of $\HY(\sigma_1^2)$, and composing with $\psi$, their images are $x_1-x_2$ and $-(y_1-y_2)$, respectively.  We conclude that the $A$-degree zero part of $\HY(\sigma_1^2)$ is isomorphic to the ideal in $R[y_1,y_2]$ generated by $\alpha:=x_1-x_2$ and $\beta:=y_1-y_2$.  As a module over $R[y_1,y_2]$, this ideal is generated by $\alpha,\beta$, modulo $(y_1-y_2)\alpha = (x_1-x_2)\beta$.  This agrees with the computation in Example \ref{ex:Hopf}.  For the higher Hochschild degree parts, see \S \ref{subsec:2strands}.

Observe that after inverting $y_1-y_2$ we have $\alpha = \frac{x_1-x_2}{y_1-y_2}\beta$, and the $y$-ified homology collapses to a copy of $R=\Q[x_1,x_2]$, generated by $\beta$.

\subsection{Flatness over $y$}
\label{subsec:flatness}

\begin{definition}
We will say that a braid $\b\in \Br_n$ (or the link $L$ it represents) is \emph{$y$-flat} if $\HY(\b)$ is free as a module over  $\yring{L}$.
\end{definition}

Below we give several examples of $y$-flat links and criteria for $y$-flatness of links.  To establish these, we recall some classical facts from commutative algebra and include their proofs for completeness. Consider the algebra $A=\Q[y_c]_{,\in \pi_0(L)}$. 
Let $C$ be a finite complex of finitely generated free graded $A$-modules.  Then, for an appropriate graded subspace $C_0\subset C$, we have
\[
C = C_0\otimes_\Q A.
\]
The differential on $C$ is $A$-linear, hence induces a differential on $C_0$, uniquely characterized by
\[
d_C = d_{C_0}\otimes 1 + (\text{higher degree terms}).
\]
Clearly, in every grading the corresponding subcomplex of $C$ is finite-dimensional, and so is its homology.   For every point $\nu\in \Q^{\pi_0(L)}$ one can define the specialization $C_{\nu}=C\otimes_A \Q_\nu$ of $C$ at $\nu$. In particular, $C_0$ is the specialization of $C$ at $0$. 

\begin{lemma}
\label{lem: zero and flatness}
Suppose that $H^*(C)\simeq H^*(C_0)\otimes A$ as graded vector spaces. Then $H^*(C)$ is a free $A$-module.
\end{lemma}

\begin{proof}
We can write the differential on $C$ as $d=d_0+d_{>0}$, where $d_{>0}$ consists of all terms of positive $y$-degree.
By Gaussian elimination one can replace $(C_0,d_0)$ by a homotopy equivalent complex $C'_0$ with zero differential.  Thus, $C_0'=H(C_0)$.  Since $d_0$ consists of the $y$-degree zero terms of $d$, this sequence of Gaussian eliminations is the reduction modulo $A_{>0}$ of an $A$-equivariant homotopy equivalence (also a sequence of Gaussian eliminations) $C\buildrel\sim \over\rightarrow C'$.  By construction, $C'\simeq H^*(C_0)\otimes A$ with some $A$-equivariant differential $d'$, hence
\[
H(C)\cong H(C')\cong H(H(C_0)\otimes A,d').
\]
But by hypothesis, $H(C)\cong H(C_0)\otimes A$.  This forces $d'=0$ (otherwise the total rank of $H(C)$ would be smaller than expected; recall that $H(C)$ is finite dimensional in each degree).  We conclude that $C\simeq H(C)\cong H(C_0)\otimes A$, which is free.
\end{proof}

We  say that $\nu$ is {\em generic} if $\nu_c\neq \nu_{c'}$ for all $c\neq c'$.

\begin{lemma}
\label{lem: zero, generic point and flatness}
Suppose that $H^*(C_0)\simeq H^*(C_{\nu})$ as graded vector spaces for  all generic $\nu$. Then $H^*(C)$ is a free  $A$-module.
\end{lemma}

\begin{proof}
Let us replace $C$ by a homotopy equivalent complex $C'$ as in the proof of Lemma \ref{lem: zero and flatness}, then
$C'\simeq H^*(C,d_0)\simeq H^*(C_0)\otimes A$.
Therefore $C'_{\nu}\simeq H^*(C_0)$. On the other hand, $C_{\nu}$ and $C'_{\nu}$ are homotopy equivalent for all $\nu$,
so for  generic $\nu$ the chain groups of $C'$ are isomorphic to its homology, and the differential in $C'_{\nu}$ vanishes.
On the other hand, the differentials in $C'$ are given by certain matrices with entries in $A$, so if these matrices vanish for all generic $\nu$ then they vanish everywhere. 
\end{proof}

\begin{lemma}
\label{lem:injective}
Suppose that $M$ and $N$ are complexes of free graded $A$--modules,  $f:M\to N$ is a chain map and $H^*(M)$ is free over $A$. Let $f^*_{\nu}:H^*(M,\Q_{\nu})\to H^*(N,\Q_{\nu})$  be the specialization of $f$ in homology with coefficients $\Q_{\nu}$. If $f_{\nu}$ is injective  for all generic $\nu$  then $f^{*}:H^*(M)\to H^*(N)$ is injective.
\end{lemma}

\begin{proof}
As in the proof of Lemma \ref{lem: zero and flatness}, we can replace $M$ by a homotopy equivalent complex $M'$ of free $A$-modules with zero differential. Suppose that $f^*(\alpha)=0$ for some $\alpha\in M'\cong H^*(M)$, then $f(\alpha)=d(\beta)$ fore some $\beta\in N$.
Let $\alpha_{\nu}$ and $\beta_{\nu}$ denote the specializations of $\alpha$ and $\beta$ at the point $\nu$, then
$f_{\nu}(\alpha_{\nu})=d_{\nu}(\beta_{\nu})$ and $f^*_{\nu}(\alpha_{\nu})=0$ for all $\nu$. For generic $\nu$, the map $f^*_{\nu}$ is injective, so $\alpha_{\nu}=0$, therefore $\alpha=0$. 
\end{proof}

\begin{lemma}
\label{lem:split is flat}
Any totally split link $L=\split(L)$ is $y$-flat.  In particular, knots are $y$-flat.
\end{lemma}

\begin{proof}
If $L=L_1\sqcup \cdots \sqcup L_r$ is a totally split link, then Proposition \ref{prop:splittensor} says
\[
\HY(L)\cong \HY(L_1)\otimes_{\Q}\cdots \otimes_{\Q}\HY(L_r).
\]
Each $L_i$ is a knot, hence is $y_i$-flat (Remark \ref{rmk:knots}), hence $L$ is $\yy$-flat.
\end{proof}

\begin{theorem}
\label{th:flat vs collapsible}
The following are equivalent:
\begin{enumerate}
\item $L$ is $y$-flat.
\item $\HKR(L)\cong \HKR(\split(L))$ upon collapse of the trigrading to the $(\deg_A, \deg_Q+\deg_T)$ bigrading.
\end{enumerate}
In particular, if $L$ is parity, then $L$ is $y$-flat.
\end{theorem}
\begin{proof}
Assume first that $\HKR(L)\cong \HKR(\split(L))$ after collapse of the trigrading.  
Then for generic $\nu$ we have 
\[
\HY(L,\Q_\nu)\cong \HY(\split(L),\Q_\nu)\cong \HKR(\split(L))\cong \HKR(L). 
\]
The first isomorphism is from the link splitting property (Corollary \ref{cor:linksplittingproperty}), the middle isomorphism holds since $\HY(\split(L))$ is $y$-flat, and the last isomorphism holds by hypothesis.  By Lemma \ref{lem: zero, generic point and flatness} $\HY(L)$ is a free $\xyring{L}$-module.

Conversely, if $\HY(L)$ is free then $\HY(L,\Q_\nu)\cong \HKR(L)$ for all $\nu$, after collapse of the trigrading, so if $\nu$ is generic then
\[
\HKR(\split(L))\cong \HY(\split(L),\Q_\nu)\cong \HY(L,\Q_\nu)\cong \HKR(L).
\]



\end{proof}

As an important corollary we have the following.
\begin{theorem}\label{thm:splittableInjectivity}
Let $L$ be arbitrary link.
\begin{enumerate}
\item[(a)] There is an injective map $f:\HY(\split(L))\to \HY(L)$ which is homogeneous but possibly of nonzero degree.
\item[(b)] Let $L'$ be obtained from $L$ by a sequence of crossing changes between different components.  If $L$ is $y$-flat (for instance, if $L$ is parity) then any splitting map $\Psi: \HY(L)\to \HY(L')$ defined as in Corollary \ref{cor:spliiting map} is injective. 
\end{enumerate}
\end{theorem}

\begin{proof}
(a) follows from (b), given that $\split(L)$ is $y$-flat (we can use Proposition \ref{prop: crossing change maps} and define $f=\Psi_{\split(L)\to L}$.

(b) Suppose $\HY(L)$ is a free $\Q[\yy]$-module, and let $\Psi:\CY(L)\rightarrow \CY(L')$ be a splitting map defined as in Proposition \ref{prop: crossing change maps}.  By Corollary \ref{cor:linksplittingproperty} $\Psi_{\nu}:\HY(L,\Q_{\nu})\to \HY(L',\Q_{\nu})$ is an isomorphism for generic $\nu$. By Lemma \ref{lem:injective} $\HY(L,\Q)\to \HY(L')$ is injective. 
\end{proof}

For instance, if $\b=(\sigma_{n-1}\cdots \sigma_2\sigma_1)^m$, then $\widehat\b$ is the $(n,m)$ torus link.  Results of Elias, Hogancamp, and Mellit compute the homology of arbitrary positive torus links.   In particular, positive torus links are parity.

\begin{corollary}
If $m,n\geq 1$ are coprime, then splitting map identifies $\HY(T(nd,md))$ with a $\Q[x_1,\ldots,x_d,y_1,\ldots,y_d,\theta_1,\ldots,\theta_d]$-submodule of $\HY(T(n,m))^{\otimes d}$.
\end{corollary}
As a special case, we see that the $y$-ified homology of $T(n,nk)$ is isomorphic to an ideal in $\Q[x_1,\ldots,x_n,y_1,\ldots,y_n,\theta_1,\ldots,\theta_n]$.  Identifying this ideal exactly is the subject of the next section.

\section{The ideal associated to a braid}
\label{sec:ideals}
In this section we compute the $y$-ified homologies of several pure braids, and identify them as ideals in $\Q[\xx,\yy,\ttheta]$.  Throughout this section let $\DC_n=D^b(R^e_n\gmod)$ be the bounded derived category of graded $(R,R)$-bimodules, as in \S\ref{subsec:oddvariables}.  Let $\one\in \SBim_n\subset \DC_n$ be the monoidal identity, and let $\one_{\YC}=\one[\yy]$ denote the monoidal identity of $\YC(\DC_n)$.


Let $A_{ij}=\sigma_{j-1}\cdots\sigma_{i+1}\sigma_i^2\sigma_{i+1}\inv\cdots \sigma_{j-1}\inv$ denote Artin's generator of the pure braid group.   Graphically, we have
 \[
 A_{ij} \ =\  
 \begin{minipage}{1.5in}
\labellist
\small
\pinlabel $i$ at 15 -7
\pinlabel $j$ at 68 -7
\endlabellist
\includegraphics[scale=.8]{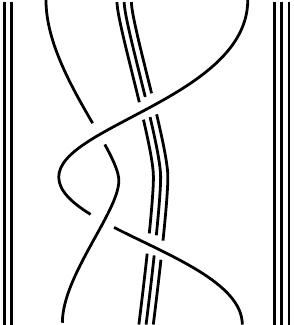}
\end{minipage}
 \]
 \vskip7pt
Since $A_{ij}$ can be split by a single positive to negative crossing change, we have a splitting map $\psi_{ij}:\FY(A_{ij})\rightarrow \one$.

\begin{definition}\label{def:splittingMap}
If $\b\in \Br_n$ is a product of the braids $A_{ij}$, let $\psi_\b:\FY(\b)\rightarrow \one$ denote the corresponding tensor product in $\YC(\SBim_n)$ of the maps $\psi_{ij}$.  Let $\IY(\b)\subset \Q[\xx,\yy,\ttheta]$ denote the image of the induced map in homology $\HY(\b)\rightarrow \HY(\one)= \Q[\xx,\yy,\ttheta]$.
\end{definition}

Since $\one_\YC$ is the monoidal identity, every square of the form
\[
\begin{diagram}
\CB\otimes \one_\YC &\rTo^{\Id\otimes g} & \CB\otimes \one_\YC(i,j)[k]\\
\dTo^{f\otimes \Id}&&\dTo^{f\otimes \Id}\\
\DB\otimes \one_\YC &\rTo^{\Id\otimes g} & \DB\otimes \one_\YC(i,j)[k]
\end{diagram}
\]
commutes, where $f:\CB\rightarrow \DB$ is any morphism of $y$-ifications and $g\in \Q[\xx,\yy,\theta]$ is any element, thought of as an endomorphism of $\one_\YC$.  Since taking homology is a functor it follows that the induced map in homology $\HY(\CB)\rightarrow \HY(\DB)$ is linear with respect to $\Q[\xx,\yy,\ttheta]$.  

It follows that $\IY(\b)$ is an ideal.  If $\HY(\b)$ is parity, then $\HY(\b)\cong \IY(\b)$, and each is free as a $\Q[\yy]$-module by Corollary \ref{cor:parityAndFlatness} and Theorem \ref{thm:splittableInjectivity}.  The goal of this section is to identify these ideals in the case when $\b=\FT_n^k$ and $\b=\LC_n$.

\begin{example}
As special cases we have \emph{Jucys-Murphy braids} $\LC_n=\sigma_{n-1}\cdots \sigma_2\sigma_1^2\sigma_2\cdots \sigma_{n-1}$, so-named because of a relation between these braids and the $q$-analogues of Jucys-Murphy elements in Hecke algebras.  These are pure braids, and can be written as products of the generators $A_{ij}$:
\[
\LC_n = A_{1,n}A_{2,n}\cdots A_{n-1,n}.
\]
Also, the full twist braid $\FT_n = \LC_2\cdots \LC_n=\prod_{i<j}A_{i,j}$ can be written as a positive product of the $A_{i,j}$.
\end{example}

\subsection{Computations on two strands}
\label{subsec:2strands}
Let $\FT_2=\sigma_1^2$ be the full-twist on two strands.  
We claim that the $y$-ified Rouquier complex for $\FT_2^k$ ($k>0$) is homotopy equivalent to the following:
\[
\FY(\FT_2^k) \ \simeq \ 
\begin{tikzpicture}[baseline=-3em]
\tikzstyle{every node}=[font=\small]
\node (a) at (0,0) {$\underline{B}[\yy](1-2k)$};
\node (b) at (4,0) {$B[\yy](3-2k)$};
\node (c) at (7,0) {$\cdots$};
\node (d) at (.5,-2) {$\cdots$};
\node (e) at (4,-2) {$B[\yy](-3+2k)$};
\node (f) at (8,-2) {$B[\yy](-1+2k)$};
\node (g) at (11,-2) {$R[\yy](2k)$};
\path[->,>=stealth',shorten >=1pt,auto,node distance=1.8cm]
([yshift=3pt] a.east) edge node[above] {$x_1-x_1'$}		([yshift=3pt] b.west)
([yshift=-2pt] b.west) edge node[below] {$y_1-y_2$}		([yshift=-2pt] a.east)
(b) edge node[above] {$x_1-x_2'$} (c)
([yshift=3pt] d.east) edge node[above] {$x_1-x_1'$}		([yshift=3pt] e.west)
([yshift=-2pt] e.west) edge node[below] {$y_1-y_2$}		([yshift=-2pt] d.east)
(e) edge node[above] {$x_1-x_2'$} (f)
(f) edge node[above] {$b$} (g);
\end{tikzpicture}
\]
where $B=B_1$.  Based on the computation of $\FT_2^k$ in \cite{Kh07}, the above is indeed a $y$-ification of the Rouquier complex for $\FT_2^{k}$, hence uniqueness of $y$-ifications proves our claim.  Any morphism $\FT_2^{\otimes k} \rightarrow \one$ can be expressed diagrammatically as 
\[
\begin{tikzpicture}[baseline=-3em]
\tikzstyle{every node}=[font=\small]
\node (a) at (0,0) {$B'$};
\node (b) at (2,0) {$B'$};
\node (c) at (4,0) {$B'$};
\node (d) at (6,0) {$\cdots$};
\node (e) at (8,0) {$B'$};
\node (f) at (10,0) {$B'$};
\node (g) at (12,0) {$R'$};
\node (m) at (0,-3) {$R'$};
\path[->,>=stealth',shorten >=1pt,auto,node distance=1.8cm]
([yshift=3pt] a.east) edge node[above] {$x_1-x_1'$}		([yshift=3pt] b.west)
([yshift=-2pt] b.west) edge node[below] {$y_1-y_2$}		([yshift=-2pt] a.east)
(b) edge node[above] {$x_1-x_2'$} (c)
([yshift=3pt] c.east) edge node[above] {$x_1-x_1'$}		([yshift=3pt] d.west)
([yshift=-2pt] d.west) edge node[below] {$y_1-y_2$}		([yshift=-2pt] c.east)
(d) edge node[above] {$x_1-x_2'$} (e)
([yshift=3pt] e.east) edge node[above] {$x_1-x_1'$}		([yshift=3pt] f.west)
([yshift=-2pt] f.west) edge node[below] {$y_1-y_2$}		([yshift=-2pt] e.east)
(f) edge node[above] {$b$} (g)
(a) edge node {} (m)
(c) edge node {} (m)
(e) edge node {} (m)
(g) edge node {} (m);
\end{tikzpicture}
\]
where we have omitted all degree shifts and abbreviated $B':=B[\yy]$, $R':=R[\yy]$.  It is an easy exercise to verify that there is a unique morphism of $y$-ifications (up to scaling by $\Q$) in which the left-most arrow is $(x_1-x_2)^{k-1}b$, the next arrow to the right is $(x_1-x_2)^{k-2}(y_1-y_2)b$, continuing in this way until the second to last arrow on the right $(y_1-y_2)^{k-2}b$, and finally that the right-most arrow is multiplication by $(y_1-y_2)^{k}$.  Thus, this must be splitting morphism up to unit scalar.  Taking $\HH$ yields the diagram:
\[
\begin{tikzpicture}[baseline=-3em]
\tikzstyle{every node}=[font=\small]
\node (a) at (0,0) {$\HH(B)[\yy]$};
\node (b) at (3,0) {$\HH(B)[\yy]$};
\node (c) at (5,-.3) {};
\node (d) at (6,0) {$\cdots$};
\node (e) at (9,0) {$\HH(B)[\yy]$};
\node (f) at (12,0) {$\HH(B)[\yy]$};
\node (g) at (15,0) {$\Q[\xx,\yy,\ttheta]$};
\node (m) at (0,-3) {$\Q[\xx,\yy,\ttheta]$};
\path[->,>=stealth',shorten >=1pt,auto,node distance=1.8cm]
(b) edge node[above] {$y_1-y_2$}		(a)
(b) edge node[above] {$x_1-x_2$} (d)
(d) edge node[above] {$x_1-x_2$} (e)
(f) edge node[above] {$y_1-y_2$}		(e)
(f) edge node[above] {$b$} (g)
(a) edge node {} (m)
(e) edge node {} (m)
(g) edge node {} (m);
\end{tikzpicture}
\]

Recall Proposition \ref{prop:odd}, which states that the map in homology induced by $b:B(-1)\rightarrow R$ identifies $\HH(B(-1))=\Hom^{\Z\times \Z}_{\DC_2}(R,B)(-1)$ with the ideal in $\Q[x_1,x_2,\theta_1,\theta_2]$ generated by $x_1-x_2$ and $\theta_1-\theta_2$.   Every element $c$ in the above complex (top row of the diagram) is a cycle if and only the homological degree of $c$ is even.  Given the above description of the link splitting map $\psi$, it follows the image of $\HY(\psi):\HY(\FT_2^k)\rightarrow \Q[\xx,\yy,\ttheta]$ equals the ideal generated by $(x_1-x_2)^{i}(y_1-y_2)^j$ ($i+j=k$) and $(x_1-x_2)^{i}(y_1-y_2)^j(\theta_1-\theta_2)$ ($i+j=k-1$).  
Define the ideal
\[
\AJ_2=(x_1-x_2,y_1-y_2,\theta_1-\theta_2).
\]
Since $(\theta_1-\theta_2)^2=0$, the image of $\HY(\psi)$ coincides with $\AJ_2^k$, and  we have just proven the following.

\begin{proposition}\label{prop:powersOfFT2}
For $k>0$, the image of $\HY(\FT_2^k)$ under the canonical embedding in $\Q[x_1,x_2,y_1,y_2,\theta_1,\theta_2]$ equals the ideal $\AJ_2^k$, as triply graded $\Q[\xx,\yy,\ttheta]$-modules.  Consequently $\HY(\FT_2^{k})$ is spanned by the tensor products of classes in $\HY(\FT_2)$.\qed
\end{proposition}

\begin{corollary}\label{cor:FT2ring}
The ring $\bigoplus_{k\geq 0}\HY(\FT_2^k)$ is isomorphic to $\Q[x_1,x_2,y_1,y_2,\a,\b]\otimes_\Q\Lambda[\theta_1,\theta_2,\gamma]$ modulo
\begin{equation}\label{eq:FT2relations}
(x_1-x_2)\b = (y_1-y_2)\a,\qquad (x_1-x_2)\gamma = (\theta_1-\theta_2)\a, \qquad (y_1-y_2)\gamma=(\theta_1-\theta_2)\b.
\end{equation}
This ring is quadruply graded, via $\deg(x_i)=(2,0,0,0)$, $\deg(y_i)=(-2,0,2,0)$, $\deg(\theta_i)=(-2,1,0,0)$, $\deg(\a)=(2,0,0,1)$, $\deg(\b)=(-2,0,2,1)$, $\deg(\gamma)=(-2,1,0,1)$.
\end{corollary}
\begin{proof}
The vector space $\bigoplus_{k\geq 0}\HY(\FT_2^k)$ has the structure of a quadruply graded algebra where three of the gradings are the usual ones in $\HY(\FT_n^k)$, and the fourth grading is by powers $k$.  The algebra structure is inherited from the composition
$$\HY\left(\FT_2^k\right)\otimes \HY\left(\FT_2^l\right)\to \HY\left(\FT_2^{k+l}\right).$$
According to Proposition \ref{prop:powersOfFT2}, this algebra is generated in degree 1, i.e.~is generated by $\HY(\FT_2)$ is generated as a module over $\Q[\xx,\yy,\ttheta]$ by $\a,\b$ and an odd variable $\gamma$, which map to $x_1-x_2$, $y_1-y_2$, respectively $(\theta_1-\theta_2)$ under the link splitting map $\HY(\FT_2)\rightarrow \Q[\xx,\yy,\ttheta]$.  The relations \eqref{eq:FT2relations} hold in $\Q[\xx,\yy,\ttheta]$, hence the also hold in $\HY(\FT_2)$ by injectivity of the link splitting map.  No further relations hold, again by examining the image of $\HY(\FT_2^k)$ under the link splitting map.
\end{proof}

We record the following important consequence of Proposition \ref{prop:powersOfFT2} below.

\begin{proposition}\label{prop:purebraidProps}
Let $\b\in \Br_n$ be expressed as a product of braids $A_{ij}$ (for instance $\b=\LC_1^{a_1}\cdots \LC_n^{a_n})$ with $a_1,\ldots,a_n\geq 0$), and let $\psi:\FY(\b)\rightarrow R[\yy]$ be the splitting map. If $\HY(\b)$ is parity, then
\begin{enumerate}
\item $\psi$ is injective in homology, and identifies $\HY(\b)$ with an ideal 
\[
\IY(\b)\subset \Q[\xx,\yy,\ttheta].
\]
\item $\HY(\b)\cong \IY(\b)$ is free as $\Q[y_1,\ldots,y_n]$-module.
\item $\IY(\b)$ is contained in the intersection of ideals $\bigcap_{i<j} \AJ_{ij}^{m_{ij}}$ where $\AJ_{ij}$ is the ideal generated by $x_i-x_j$, $y_i-y_j$ and $\theta_i-\theta_j$, and $m_{ij}$ is the multiplicity of $A_{ij}$ in the expression of $\b$ (that is, the linking number between the $i$th and $j$th strands).
\end{enumerate}
\end{proposition}
\begin{proof}
Statement (1) is an immediate consequence of Theorem \ref{thm:splittableInjectivity}, and (2) is an immediate consequence of Theorem \ref{th:flat vs collapsible}.  For (3), observe that the splitting map $\FY(\b)\rightarrow R[\yy]$ factors through the splitting map $\psi_{ij}^{\otimes m} :\FY(A_{ij}^{m})\rightarrow R[\yy]$, where $m=m_{ij}$ is the multiplicity.  Hence
\[
\IY(\b)\subset \IY(A_{ij}^m).
\]
On the other hand, $A_{ij}^m$ is conjugate to $\sigma_i^{2m}$, hence Proposition \ref{prop:powersOfFT2} implies that the map in homology induced by $\psi_{ij}^{\otimes m}$ has image
\[
\psi_{ij}^{\otimes m}(\HY(A_{ij}^m)) = \AJ_{ij}^m.
\]
Thus, $\IY(\b)$ is contained in the intersection of all such ideals, which is (3).
\end{proof}


\begin{remark}
The condition $a_i>0$ does not imply parity of $\LC_1^{a_1}\cdots \LC_n^{a_n}$.  For example, $\LC_3^2$ has homology in both even and odd homological degrees.  However, it is conjectured in \cite{GNR,OR3} that $\LC_1^{a_1}\cdots \LC_n^{a_n}$ is parity when $a_1\ge a_2\ge \cdots \ge a_n$.
\end{remark}




\subsection{Negative twists on two strands}
\label{subsec: negative}

In contrast with the positive powers of the  full twist, the homologies of the powers of negative full twist are rather subtle. In particular, they are not $y$-flat and do not inject into the homology of the unlink. To simplify the notations, we consider the Hochschild degree zero part of the homology of $\FT_2^{-2}$. The corresponding $y$-ified complex has the form:

\begin{tikzpicture}
\tikzstyle{every node}=[font=\small]
\node (a) at (0,0) {${R}[\yy]$};
\node (b) at (3,0) {$B[\yy]$};
\node (e) at (6,0) {$B[\yy]$};
\node (f) at (9,0) {$B[\yy]$};
\node (g) at (12,0) {$\underline{B}[\yy]$};
\path[->,>=stealth',shorten >=1pt,auto,node distance=1.8cm]
(a) edge node[above] {$b^*$}		(b)
([yshift=-2pt] e.west) edge node[below] {$y_1-y_2$}		([yshift=-2pt] b.east)
([yshift=3pt] b.east) edge node[above] {$x_1-x_1'$} ([yshift=3pt] e.west)
([yshift=-2pt] g.west) edge node[below] {$y_1-y_2$}		([yshift=-2pt] f.east)
(e) edge node[above] {$x_1-x_2'$} (f)
([yshift=3pt] f.east) edge node[above] {$x_1-x_1'$} ([yshift=3pt] g.west);
\end{tikzpicture}
 
After taking $\HH^0$ we get

\begin{tikzpicture}
\tikzstyle{every node}=[font=\small]
\node (a) at (0,0) {$\Q[\xx,\yy]$};
\node (b) at (3,0) {$\Q[\xx,\yy]$};
\node (e) at (6,0) {$\Q[\xx,\yy]$};
\node (f) at (9,0) {$\Q[\xx,\yy]$};
\node (g) at (12,0) {$\Q[\xx,\yy]$};
\path[->,>=stealth',shorten >=1pt,auto,node distance=1.8cm]
(a) edge node[above] {$1$}		(b)
(e) edge node[below] {$y_1-y_2$}		(b)
(g) edge node[below] {$y_1-y_2$}		(f)
(e) edge node[above] {$x_1-x_2$} (f);
\end{tikzpicture}
 
Therefore 
\[
\HY^0(\FT_2^{-2})\simeq \Q[\xx,\yy]\oplus \Q[\xx,\yy]/(x_1-x_2,y_2-y_2).
\]
With respect to the above decomposition, the class $(1,0)$ has tridegree $(-4,0,0)$ and the class $(0,1)$ has tridegree $(0,0,-1)$.  In particular, $\HY(\FT_2^{-2})$ is not free as a $\Q[y]$-module, and is not parity.

One can also consider the homology with  coefficients $\Q_{\nu}$. For $\nu_1=\nu_2$ (say, $\nu=0$)
one has
\[
\HY^0(\FT_2^{-2},\Q_{\nu_1=\nu_2})  \cong \HKR^0(\FT_2^{-2})\simeq \Q[\xx]\oplus \Q[\xx]/(x_1-x_2).
\] 
On the other hand,
\[
\HY^0(\FT_2^{-2},\Q_{\nu_1\neq \nu_2})\simeq \HKR^0(\one_2) \cong \Q[\xx].
\] 
For $\nu_1\neq \nu_2$, the spectral sequence from $\HKR(\FT_2^{-2})$ to $\HY^0(\FT_2^{-2},\Q_{\nu})$
does not collapse and has a nontrivial differential.	

\subsection{The Jucys-Murphy braids}
In this section we  identify $\HY(\LC_n)$ with an ideal inside $\Q[\xx,\yy]$.  Our goal is to prove the following.

\begin{proposition}
\label{prop:jucys murphy}
The homology $\HY(\LC_n)$ is isomorphic to the product of ideals $\AJ_{1,n}\AJ_{2,n}\cdots \AJ_{n-1,n}\subset \Q[\xx,\yy]$.  The graded dimension is
\[
\PC_{\LC_n}(q,t,a) = \frac{1+a}{(1-q)(1-t)}\left(\frac{q+t-qt+a}{(1-q)(1-t)}\right)^{n-1}.
\]
\end{proposition}
\begin{proof}
We first observe that $\IY(\LC_n)$ is contained in $\AJ_{1,n}\cap\AJ_{2,n}\cap \cdots \cap \AJ_{n-1,n}$, and contains the product $\AJ_{1,n}\cdots\AJ_{n-1,n}$ by Proposition \ref{prop:purebraidProps}, together with the fact that $\LC_n=A_{1,n}A_{2,n}\cdots A_{n-1,n}$.  The ideals $\AJ_{i,n}$ are independent from one another, hence the intersection of ideals here is equal to the product.

To compute the graded dimension, consider $\AJ_{1,n}\cdots\AJ_{n-1,n}$ as a triply graded $\Q[\xx,\yy]$-module.  It is generated by
\[
\{x_n,y_n,\theta_n\}\sqcup \{\a_1,\ldots,\a_{n-1},\b_1,\ldots,\b_{n-1},\gamma_1,\ldots,\gamma_{n-1}\}
\]
modulo the relations $(x_i-x_n)\b_i = (y_i-y_n)\a_i$.  Here, $\a_i=x_i-x_n$, $\b_i = y_i-y_n$, and $\gamma_i=\theta_i-\theta_n$.  There are no higher syzygies, hence the graded dimension is as claimed (the generators $\a_i,\b_i,\gamma_i$ have degree $q,t,a$ and the relations have degree $qt$).
\end{proof}

\section{The full twist ideals}
\label{sec:ftideals}

In this section, we compute the homology of the full twist and prove Theorem \ref{th: intro ft}.
We start with some general  results.

\subsection{Canonicalness of the splitting map $\FT_n\rightarrow \one$.}
The main goal of this section is to prove the following.

\begin{proposition}\label{prop:splittingmapuniqueness}
The Hochschild degree zero part of $\Hom_{\YC(\DC_n)}^{\Z\times \Z\times \Z}(\FT_n,\one)$ is isomorphic to $\Q[\xx,\yy]$, generated by the splitting map $\Psi:\FT_n\rightarrow \one$.
\end{proposition}
The splitting map $\Psi$ is clearly not null-homotopic because it induces an isomorphism $\HY(\FT_n,\Q_\nu)\rightarrow \HY(\one,\Q_\nu)$ for generic $\nu$.  Thus we need only prove the first statement. Given that $\FT_n$ is invertible, this is equivalent to showing that the Hochschild degree zero part of $\Hom_{\YC(\DC_n)}^{\Z\times \Z\times \Z}(\one,\FT_n\inv)$ is isomorphic to $\Q[\xx,\yy]$.

\begin{lemma}\label{lemma:neghom}
For any complexes $C,C'\in \KC^b(\SBim_{n-1})$ the complex of Homs satisfies
\[
\underline{\Hom}_{\KC(\SBim_n)}^{\Z\times \Z} (\one_n , (C\sqcup \one_1)\otimes F(\sigma_{n-1}) \otimes (C'\sqcup \one_1)) \simeq \underline{\Hom}_{\KC(\SBim_{n-1})}^{\Z\times \Z} (\one_{n-1} , C\otimes C'),
\]
\[
\underline{\Hom}_{\KC(\SBim_n)}^{\Z\times \Z} (\one_n , (C\sqcup \one_1)\otimes F(\sigma_{n-1}\inv) \otimes (C'\sqcup \one_1)) \simeq 0.
\]
\end{lemma}
\begin{proof}
These relations are obtained by taking Hochschild degree zero part of the Markov moves for Khovanov-Rozansky homology (see Lemma \ref{lemma:HHmarkov}).
\end{proof}

Throughout the remainder of this section, let $\LC_k= \sigma_{k-1}\cdots \sigma_2\sigma_1^2\sigma_2\cdots \sigma_{k-1}\in \Br_n$.  By abuse of notation, we also let $\LC_k$ denote the Rouquier complex $F(\LC_k)\in \KC^b(\SBim_n)$ or the $y$-ified Rouquier complex $\FY(\LC_k)$.  The meaning of $\LC_k$ will be clear from context.

\begin{lemma}\label{lemma:neghom2}
Let $\CB=(C,w,\Delta)\in \YC(\SBim_{n-1})$ be a $y$-ification.  Then $\underline{\Hom}_{\YC(\SBim)}(\one_n, (C\sqcup \one_1) \otimes \LC_n\inv) \simeq \underline{\Hom}_{\YC(\SBim)}(\one_n, C\sqcup \one_1)$.
\end{lemma}

\begin{proof}
We will prove by induction that
\begin{equation}\label{eq:inductivestep}
\underline{\Hom}_{\YC(\SBim)}(\one_n, (C\sqcup \one_1) \otimes (\one_{n-k}\sqcup \LC_{k}\inv)) \simeq \underline{\Hom}_{\YC(\SBim)}(\one_n, (C\sqcup \one_1)\otimes (\one_{n-k+1}\sqcup \LC_{k-1}\inv)).
\end{equation}
The base case $k=0$ is a tautology.  Now, fix $1\leq k\leq n-1$.  Recall diagram \eqref{eq:crossingchange}, which defines the splitting map $\psi:\FY(\sigma_{n-k})\rightarrow \FY(\sigma_{n-k}\inv)$.   The cone on this map is given in \eqref{eq:splitmapcone}.  Tensoring this map on the right with $\FT(\sigma_{n-k}^{-1})$ gives a map $\psi': \one_n\rightarrow \FY(\sigma_{n-k}^{-2})$ such that
\[
\Cone(\psi') \ \ \simeq  \ \ 
\left(\begin{tikzpicture}[baseline=-.5em]
\tikzstyle{every node}=[font=\small]
\node (a) at (0,0) {$\FY(\sigma_{n-k}\inv)[1](-1)$};
\node (b) at (4,0) {$\FY(\sigma_{n-k}\inv)(1)$};
\path[->,>=stealth',shorten >=1pt,auto,node distance=1.8cm,
  thick]
([yshift=3pt] a.east) edge node[above] {}		([yshift=3pt] b.west)
([yshift=-2pt] b.west) edge node[below] {}		([yshift=-2pt] a.east);
\end{tikzpicture}\right) =: Z'
\]
It follows that there is a distinguished triangle in $\YC(\SBim_n)$
\[
\one_n \rightarrow \FY(\sigma_{n-k}^{-2})\rightarrow Z' \rightarrow \one_n[1].
\]
Tensoring on the left with $\FY(\sigma_n\inv\cdots \sigma_{n-k+1}\inv)$ and on the right with $\FY(\sigma_{n-k+1}\inv\cdots \sigma_{n-1}\inv)$ gives a distinguished triangle in $\YC(\SBim_n)$:
\begin{equation}\label{eq:Ztriang}
\one_{n-k+1}\sqcup \LC_{k+1}\inv \rightarrow \one_{n-k}\sqcup \LC_k\inv \rightarrow Z\rightarrow (\one_{n-k+1}\sqcup \LC_{k+1}\inv)[1],
\end{equation}
where
\[
Z \ := \ 
\left(\begin{tikzpicture}[baseline=-.5em]
\tikzstyle{every node}=[font=\small]
\node (a) at (0,0) {$\FY(\b)[1](-1)$};
\node (b) at (4,0) {$\FY(\b)(1)$};
\path[->,>=stealth',shorten >=1pt,auto,node distance=1.8cm,
  thick]
([yshift=3pt] a.east) edge node[above] {}		([yshift=3pt] b.west)
([yshift=-2pt] b.west) edge node[below] {}		([yshift=-2pt] a.east);
\end{tikzpicture}\right),
\]
where we have abbreviated by writing $\b:=(\sigma_{n-1} \cdots \sigma_{n-k+1}  \sigma_{n-k}  \sigma_{n-k+1}  \cdots \sigma_{n-1})\inv$.  The precise form of the middle maps is not particularly important for now.

Note that $(\sigma_{n-1} \cdots \sigma_{n-k+1}  \sigma_{n-k}  \sigma_{n-k+1}  \cdots \sigma_{n-1})\inv \simeq (\sigma_{n-k}\cdots \sigma_{n-2}\sigma_{n-1}\sigma_{n-2}\cdots\sigma_{n-k})\inv$ as braids.  Graphically this is
\[
\ig{.8}{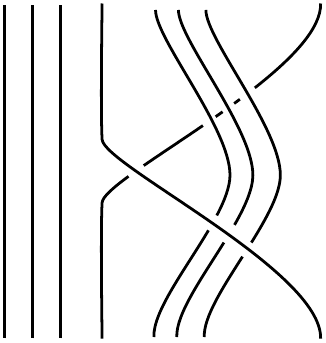} \ \ \simeq \ \ \ig{.8}{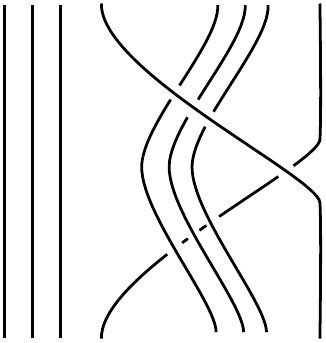}.
\]
Observe:
\begin{eqnarray*}
\b \simeq \sigma_{n-k}\inv\cdots \sigma_{n-2}\inv\sigma_{n-1}\inv\sigma_{n-2}\inv \cdots\sigma_{n-k}\inv
& \Rightarrow  & \underline{\Hom}_{\KC^b(\SBim_n)}^{\Z\times\Z}(\one_n,(C \sqcup \one_1)\otimes F(\b))\simeq 0\\
& \Rightarrow  & \underline{\Hom}_{\YC(\SBim_n)}^{\Z\times\Z}(\one_n,(\CB \sqcup \one_1)\otimes \FY(\b))\simeq 0\\
& \Rightarrow  & \underline{\Hom}_{\YC(\SBim_n)}^{\Z\times\Z}(\one_n,(\CB \sqcup \one_1)\otimes Z)\simeq 0\\
\end{eqnarray*}

The first implication holds by Lemma \ref{lemma:neghom}.  The second holds since the second line represents a $y$-ification of the contractible complex from the first line (compare with Remark \ref{rmk:HH(CB)}), together with Lemma \ref{lemma:contractible}.  The third implication holds since the third line can be expressed as a mapping cone constructed from two copies of the second line. 
 The equivalence \eqref{eq:inductivestep} follows from this together with the distinguished triangle \eqref{eq:Ztriang}.  Iterating this gives the Lemma.
\end{proof}

\begin{proof}[Proof of Proposition \ref{prop:splittingmapuniqueness}]
Since $\FT_n = \LC_2\LC_3\cdots \LC_n$, the fact that
\[
\underline{\Hom}_{\YC(\SBim_n)}^{\Z\times \Z}(\FT_n, \one)\simeq \underline{\Hom}_{\YC(\SBim_n)}^{\Z\times \Z}(\one, \FT_n\inv)\simeq \underline{\End}_{\YC(\SBim_n)}^{\Z\times \Z}(\one) = \Q[\xx,\yy]
\]
follows from repeated application of Lemma \ref{lemma:neghom2}.  Taking homology of these hom complexes gives the first statement. To see that the splitting map $\Psi$ is a generator of $\underline{\Hom}_{\YC(\SBim_n)}^{\Z\times \Z}(\one, \FT_n)$, it suffices to show that $\Psi$ is not null-homotopic.  But this follows from the fact that $\Psi$ becomes a homotopy equivalence after inverting each $y_i-y_j$.
\end{proof}

\subsection{The symmetric group action on $\HY(\FT_n^k)$}
\label{subsec: sn action}

Let $\b$ be a braid.  Since the full twist is central in the braid group, we have a homotopy equivalence $\tau_\b: \FY(\b)\otimes \FT_n \rightarrow \FT_n\otimes\FY(\b)$.   The group of automorphisms of $\FY(\b)\otimes \FT_n$ in $\YC(\DC_n)$ is isomorphic to $\Q^\times$, since $\FY(\b)\otimes \FT_n$ is invertible.  Thus the maps $\tau_\b$ are unique up to homotopy and unit scalar.  For the moment we will choose $\tau_\b$ arbitrarily, but we will see in a moment that canonicalness of the splitting morphism $\Psi:\FT_n\rightarrow \one$ will fix the scalar in a canonical way.

The morphisms which commute $\FY(\b)$ past $\FT_n$ can be composed, yielding homotopy equivalences $\tau_{\b,k}:\FY(\b)\otimes \FT_n^{\otimes k} \rightarrow \FT_n^{\otimes k} \otimes \FY(\b)$ for all $k$.  Then we have an automorphism $\phi_\b$ of $\Hom_{\YC(\DC_n)}^{\Z\times \Z\times \Z}(\one,\FT_n^{\otimes k})$ defined by
\begin{equation}\label{eq:braidaction}
\phi_\b(f) = \begin{diagram}[small] \one & \rTo^{\simeq} & \b \one \b\inv & \rTo^{\Id_\b f \Id_{\b\inv}} &  \b \FT_n^k \b\inv & \rTo^{\tau_{\b,k} \Id_{\b\inv}} &  \FT_n^k \b \b\inv & \rTo^{\simeq} &  \FT_n^k\end{diagram}
\end{equation}

\begin{remark}
There are several abuses of notation which we have made above, and will continue to make throughout this section.  First, we are denoting $\FY(\b)$ simply by $\b$.  Second, we are omitting the symbol $\otimes$ for brevity.  Finally, we are omitting all explicit occurrences of the grading shifts $(i,j)[k]$.   A more precise exposition would incorporate a coherent family of isomorphisms $C(i,j)[k]\cong \one(i,j)[k]\otimes C\cong C\otimes \one(i,j)[k]$.
\end{remark}

\begin{theorem}\label{thm:BraidGpAction}
The braiding morphisms $\tau_{\sigma_i^\pm}$ can be chosen so that $\b\mapsto \phi_\b$ defines an action of $\Br_n$ on $\Hom_{\YC(\DC_n)}^{\Z\times \Z\times \Z}(\one,\FT_n^{\otimes k})$ by degree preserving linear automorphisms.  This action factors through the natural quotient map $\pi:\Br_n\rightarrow S_n$.  In fact, if $\Psi:\FT\rightarrow \one$ denotes the splitting map, then post-composing with $\Psi^{\otimes k}$ defines an $S_n$-equivariant injective map
\[
\Hom_{\YC(\DC_n)}^{\Z\times \Z\times \Z}(\one,\FT_n^{\otimes k}) \hookrightarrow \End_{\YC(\DC_n)}^{\Z\times\Z\times \Z}(\one)=\Q[\xx,\yy,\ttheta].
\]
\end{theorem}

\begin{notation}
In what follows, $\b\in \Br_n$ will be a fixed braid and $w=\pi(\b)$ will be the corresponding permutation.  
\end{notation}

The proof of Theorem \ref{thm:BraidGpAction} will occupy the remainder of this section.  We will adopt the graphical notation for morphisms in a monoidal category.  The identity morphism of $\b$ (resp.~$\b\inv$) in $\YC(\DC_n)$ will be denoted by a dashed upward (resp.~downward) pointing line.  We will choose a homotopy equivalence $\one\rightarrow  \b\b\inv$ with inverse equivalence $\b\b\inv\rightarrow \one$.  These will be denoted diagrammatically by
\[
\ig{.8}{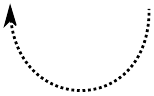} \qquad \text{ and } \qquad \ig{.8}{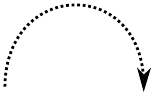},
\]
respectively.  These satisfy relations of the form
\[
\ig{.8}{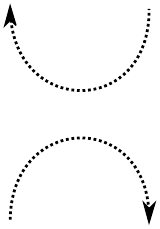} \ \ =\ \  \ig{.8}{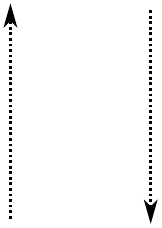},\qquad\qquad \ig{.8}{circle} \ \ = \ \ \Id_{\one}.
\]
We denote similarly a chosen pair of inverse equivalences $\b\inv \b \leftrightarrow \one$.  These satisfy the same relations as above, with the orientations reversed.  We may choose these morphisms so that the dashed lines the additional planar isotopy relations
\[
\ig{.8}{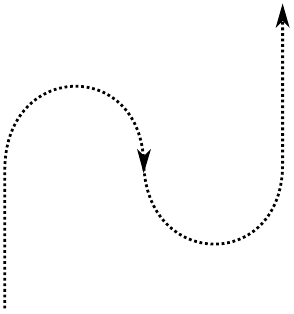} \ \  = \ \ \ig{.8}{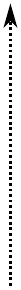} \ \ = \ \ \ig{.8}{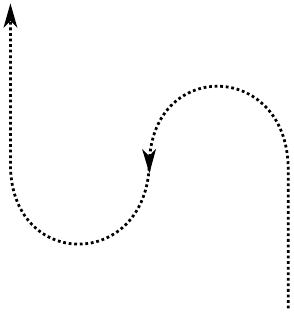},
\]
and similarly with the orientations reversed.

Polynomials $g\in \Q[\xx,\yy,\ttheta]$ are thought of as endomorphisms of $\one\in \YC(\DC_n)$, hence a string diagram with polynomial floating in the empty regions yields a well-defined morphism in $\YC(\DC_n)$.  We may define a braid group action on $\Q[\xx,\yy,\ttheta]$ by placing a polynomial $g$ inside a circle labeled by $\b$.  
\begin{lemma}\label{lemma:braidactiononpolys}
The braid group action on polynomials descends to the usual action of the symmetric group by permuting variables:
\[
\begin{minipage}{.5in}
\labellist
\small
\pinlabel $g$ at 25 25
\endlabellist
\includegraphics[scale=.8]{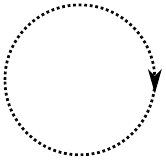}
\end{minipage}
 \ \ \simeq \ \ w(g)
\]
for all $g\in \Q[\xx,\yy,\ttheta]$.
\end{lemma}
\begin{proof}
Use the fact that $\Id_\b\otimes g \simeq w(g)\otimes \Id_{\b}$ and dashed circles evaluate to 1.
\end{proof}

The identity morphism of $\FT_n$ will be denoted by a solid, unoriented line.  The braiding morphisms $\tau_\b$, $\tau_{\b\inv}$, $\tau_\b\inv$, and $\tau_{\b\inv}\inv$ will be denoted by
\[
\ig{.8}{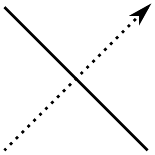},\qquad \qquad \ig{.8}{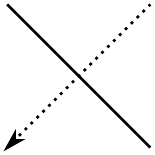},\qquad \qquad\ig{.8}{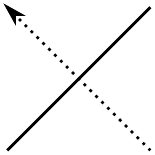},\qquad \qquad\ig{.8}{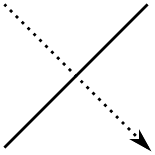}.
\]
Up to redefining $\tau_{\b\inv}$, we may assume that 
\begin{equation}\label{eq:twistedBraiding}
\ig{.8}{inversebraidinginverse} \ \ \simeq \ \ \ig{.8}{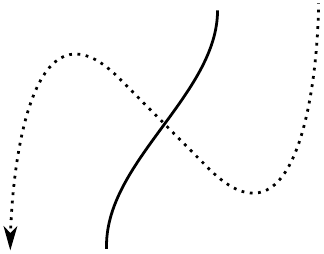}, \qquad\qquad \ig{.8}{inversebraiding} \ \ \simeq \ \ \ig{.8}{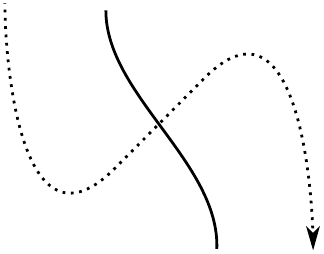}.
\end{equation}
These braiding morphisms and their inverses satisfy graphical relations:
\[
\ig{.8}{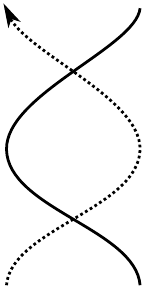} \ \ \simeq \ \ \ig{.8}{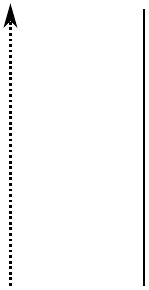} \qquad \qquad \ig{.8}{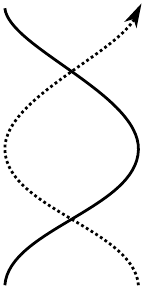} \ \ \simeq \ \ \ig{.8}{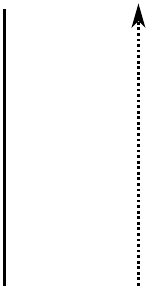},
\]
and similarly with the orientation on the dashed line reversed.

In this graphical notation, a morphism $f\in \Hom_{\YC(\DC_n)}^{\Z\times \Z\times \Z}(\one,\FT_n)$ will be denoted by a labelled univalent vertex:
\[
\begin{minipage}{.3in}
\labellist
\small
\pinlabel $f$ at 12 12
\endlabellist
\includegraphics[scale=.8]{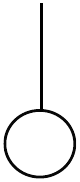}
\end{minipage} \ \ \in \ \  \Hom_{\YC(\DC_n)}^{\Z\times \Z\times \Z}(\one,\FT_n) \ \ 
\]

Recall the automorphism $\phi_\b$ of $\Hom_{\YC(\DC_n)}^{\Z\times \Z\times \Z}(\one,\FT_n)$ defined in \eqref{eq:braidaction}, for $k=1$.  Graphically, we have
\[
\phi_\b(f)  \ \  =  \ \ \begin{minipage}{.6in}
\labellist
\small
\pinlabel $f$ at 24 23
\endlabellist
\includegraphics[scale=.8]{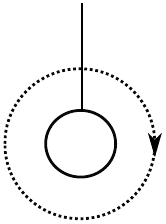}
\end{minipage}.
\]
\begin{remark}
Strictly speaking, to make sense of the horizontal dashed line meeting the vertical solid line, one must perturb so that the dashed line is not horizontal at the crossing.  There are two ways of doing this, but they are homotopic by \eqref{eq:twistedBraiding}:
\[
\begin{minipage}{.6in}
\labellist
\small
\pinlabel $f$ at 24 23
\endlabellist
\includegraphics[scale=.8]{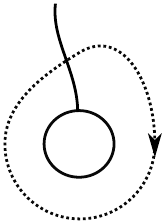}
\end{minipage}
\ \ \simeq \ \ 
\begin{minipage}{.6in}
\labellist
\small
\pinlabel $f$ at 24 23
\endlabellist
\includegraphics[scale=.8]{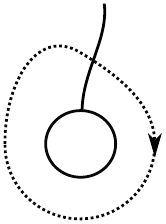}
\end{minipage}.
\]
\end{remark}

We state the following for $k=1$ for simplicity, though its statement and proof extend immediately to a more general situation (for instance $\FT^{\otimes k}_n$). 
\begin{lemma}
The map $\phi_\b$ is an automorphism of $\Hom_{\YC(\DC_n)}^{\Z\times \Z\times \Z}(\one,\FT_n)$, with  inverse $\phi_{\b\inv}$.  Furthermore, for each $f\in \Hom_{\YC(\DC_n)}^{\Z\times \Z\times \Z}(\one,\FT_n)$, the following square commutes up to homotopy:
\[
\begin{diagram}
\b \FT_n & \rTo^{\tau_\b} & \FT_n \b \\
\uTo^{\Id_\b f} && \uTo^{\phi_\b(f)\Id_\b}\\
\b & \rTo^{=} &  \b 
\end{diagram},\qquad\qquad\qquad
\begin{minipage}{1in}
\labellist
\small
\pinlabel $f$ at 49 42
\endlabellist
\includegraphics[scale=.8]{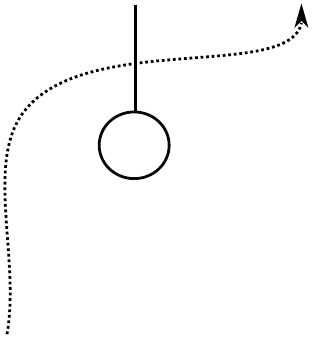}
\end{minipage}
\ \ \simeq \ \ 
\begin{minipage}{1in}
\labellist
\small
\pinlabel $f$ at 65 41
\endlabellist
\includegraphics[scale=.8]{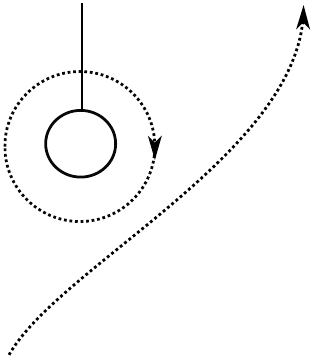}
\end{minipage}.
\]
\end{lemma}
\begin{proof}
The following shows that $\phi_{\b\inv}(\phi_\b(f))\simeq f$:
\[
\begin{minipage}{.8in}
\labellist
\small
\pinlabel $f$ at 34 33
\endlabellist
\includegraphics[scale=.8]{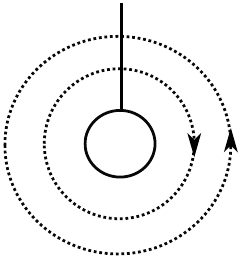}
\end{minipage}
 \ \ \simeq \ \ 
 \begin{minipage}{.8in}
\labellist
\small
\pinlabel $f$ at 32 33
\endlabellist
\includegraphics[scale=.8]{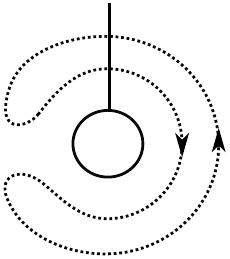}
\end{minipage}
 \ \ \simeq \ \ 
 \begin{minipage}{.8in}
\labellist
\small
\pinlabel $f$ at 32 34
\endlabellist
\includegraphics[scale=.8]{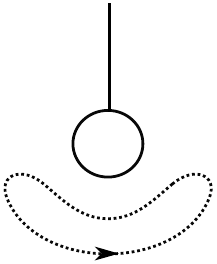}
\end{minipage}
\ \ \simeq \ \ 
\begin{minipage}{.3in}
\labellist
\small
\pinlabel $f$ at 12 12
\endlabellist
\includegraphics[scale=.8]{fig/ftclass}
\end{minipage} .
\]
A similar argument shows that $\phi_\b(\phi_{\b\inv}(f))\simeq f$.  The second statement follows from
\begin{equation}\label{eq:mapslideproof}
\begin{minipage}{1in}
\labellist
\small
\pinlabel $f$ at 23 61
\endlabellist
\includegraphics[scale=.8]{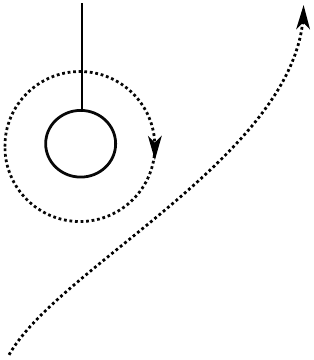}
\end{minipage}
\ \ \simeq \ \ 
\begin{minipage}{1in}
\labellist
\small
\pinlabel $f$ at 23 61
\endlabellist
\includegraphics[scale=.8]{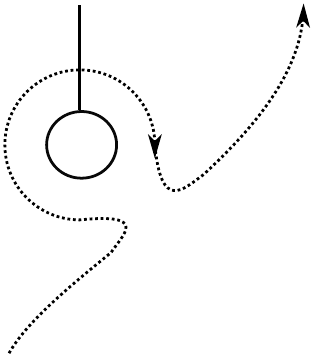}
\end{minipage}
\end{equation}
followed by an isotopy.
\end{proof}

\begin{proof}[Proof of Theorem \ref{thm:BraidGpAction}]
Let $\Psi:\FT_n\rightarrow \one$ denote the splitting morphism.   The space of maps $\FT_n\rightarrow \one$ modulo homotopy is one dimensional generated by $\Psi$ by Proposition \ref{prop:splittingmapuniqueness},  hence any automorphism of $\Hom_{\YC(\DC_n)}^{\Z\times \Z\times \Z}(\FT_n,\one)$ sends $\Psi$ to a unit multiple of itself.  It follows that, for each braid $\b\in \Br_n$ there is a scalar $0\neq c_\b\in \Q$ such that
\[
\begin{minipage}{.6in}
\labellist
\small
\pinlabel $\Psi$ at 24 41
\endlabellist
\includegraphics[scale=.8]{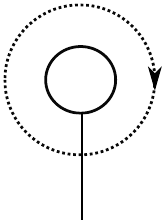}
\end{minipage}
\ \ \simeq \ \ c_\b \begin{minipage}{.6in}
\labellist
\small
\pinlabel $\Psi$ at 24 41
\endlabellist
\includegraphics[scale=.8]{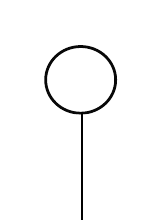}
\end{minipage}
\]
By replacing $\tau_\b: \b\FT_n\rightarrow \FT_n\b$ by $c_\b \tau_\b$ if necessary, we may assume that $c_\b=1$ for all $\b$. An argument similar to \eqref{eq:mapslideproof} shows that
\begin{equation}\label{eq:psislide}
\begin{minipage}{1.1in}
\labellist
\small
\pinlabel $\Psi$ at 49 42
\endlabellist
\includegraphics[scale=.8]{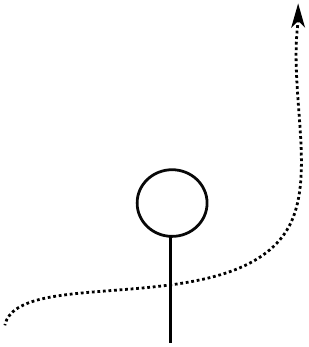}
\end{minipage}
\ \ \simeq \ \ 
\begin{minipage}{1.1in}
\labellist
\small
\pinlabel $\Psi$ at 66 42
\endlabellist
\includegraphics[scale=.8]{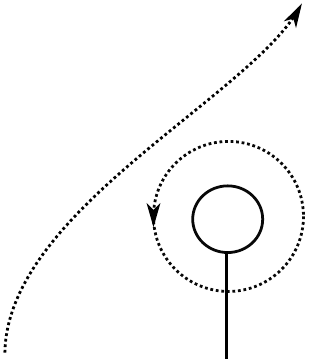}
\end{minipage}
\ \ \simeq \ \ 
\begin{minipage}{1.1in}
\labellist
\small
\pinlabel $\Psi$ at 66 42
\endlabellist
\includegraphics[scale=.8]{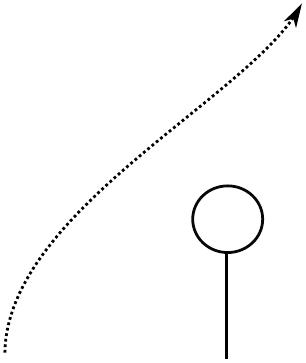}
\end{minipage}.
\end{equation}

The braid group action on $\Hom_{\YC(\DC_n)}^{\Z\times \Z\times\Z}(\one,\FT_n^{\otimes k})$ is defined by
\[
\begin{minipage}{1.6in}
\labellist
\small
\pinlabel $f$ at 72 24
\pinlabel $\cdots$ at 85 62
\endlabellist
\includegraphics[scale=.8]{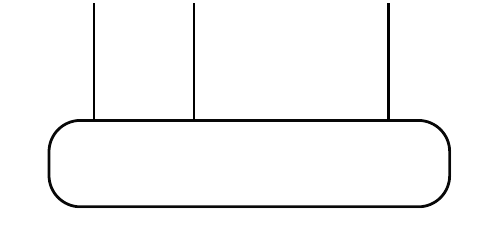}
\end{minipage}
\ \ \ \mapsto  \ \ \ \begin{minipage}{1.7in}
\labellist
\small
\pinlabel $\b$ at 150 40
\pinlabel $f$ at 72 24
\pinlabel $\cdots$ at 85 62
\endlabellist
\includegraphics[scale=.8]{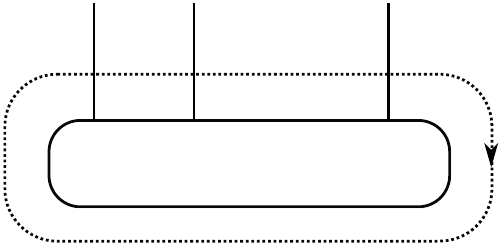}
\end{minipage} \ \  =: \ \phi_\b(f)
\]
for all $\b\in \Br_n$ and all $f\in \Hom_{\YC(\DC_n)}^{\Z\times \Z\times\Z}(\one,\FT_n^{\otimes k})$.  Now we can compute
\[
\Psi^{\otimes k}\circ \phi_\b(f) \ \ = \ \ 
\begin{minipage}{1.7in}
\labellist
\small
\pinlabel $f$ at 72 24
\pinlabel $\Psi$ at 27 81
\pinlabel $\Psi$ at 56 81
\pinlabel $\Psi$ at 113 81
\pinlabel $\cdots$ at 85 81
\endlabellist
\includegraphics[scale=.8]{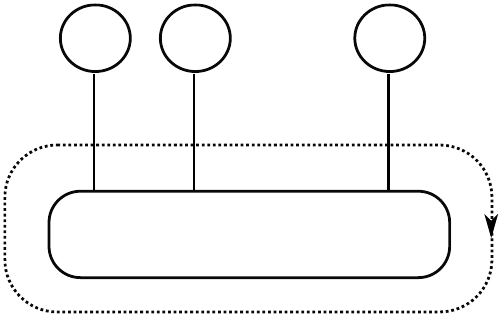}
\end{minipage}
\ \ \simeq \ \ 
\begin{minipage}{1.8in}
\labellist
\small
\pinlabel $f$ at 72 24
\pinlabel $\Psi$ at 27 62
\pinlabel $\Psi$ at 56 62
\pinlabel $\Psi$ at 113 62
\pinlabel $\cdots$ at 85 62
\endlabellist
\includegraphics[scale=.8]{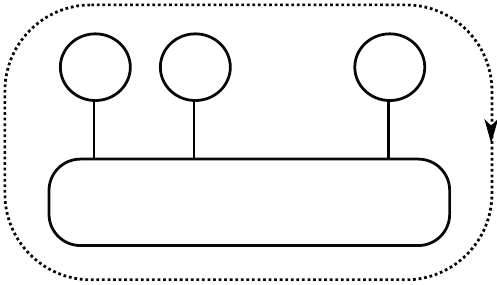}
\end{minipage}
\ \ \simeq \ \ w(\Psi^{\otimes k}\circ f).
\]
In the second equivalence we used \eqref{eq:psislide}, and in the last equivalence we used Lemma \ref{lemma:braidactiononpolys}.

This proves that $\Psi\circ \phi_\b(f)\simeq w(\Psi\circ f)$.  Since post-composing with $\Psi$ defines an injective map $\HY(\FT_n)\rightarrow \Q[\xx,\yy,\ttheta]$, we have $\phi_\b(f)\simeq \phi_{\b'}(f)$ whenever $\b$ and $b'$ are associated to the same permutation.  This completes the proof.
\end{proof}

\subsection{The full twist ideal}
In \cite{EH} Elias and Hogancamp computed $\HKR(\FT_n)$ (as triply graded vector spaces), later in \cite{Hog17a} Hogancamp generalized this computation to $\HKR(FT_n^k)$ for all $n,k>0$. In particular, one gets the following.

\begin{theorem}[\cite{Hog17a}]
\label{th: ft parity}
The Khovanov-Rozansky homology $\HKR(FT_n^k)$ is concentrated in even homological degrees for all $n,k>0$.
\end{theorem}

Still, the structure of $\HKR(\FT^k_n)$ as a module over $\Q[\xx,\ttheta]$ as well as the description of multiplication maps $\HKR(\FT_n^{k})\otimes_\Q \HKR(\FT_n^{l})\to \HKR(\FT_n^{k+l})$ remained unknown. The results of Section \ref{sec:ideals} immediately show the relevance of $y$-ified homology in understanding these problems.  We summarize this below.

\begin{proposition}\label{prop:FTbasicProps}
The following are true:
\begin{enumerate}
\item The link splitting map is injective in homology, hence identifies $\HY(\FT_n^k)$ with $\IY(\FT_n^k)\subset \Q[\xx,\yy,\thth]$.
\item The $y$-ified homology $\HY(\FT_n^{k})$ is free as $\Q[\yy]$-module, and
\[
\HKR(\FT_n^{k})\simeq \HY(\FT_n^{k})/(y_1,\ldots,y_n)\HY(\FT_n^{k}).
\]
\item The multiplication maps fit into a commutative diagram of the form
\[
\begin{diagram}
\HY(\FT_n^{\otimes k})\otimes_\Q \HY(\FT_n^{\otimes l})  & \rTo & \IY(\FT_n^k) \otimes_\Q\IY(\FT_n^l) \\
\dTo && \dTo \\
\HY(\FT_n^{\otimes k+l}) & \rTo & \IY(\FT_n^{k+l}) 
\end{diagram}
\]
\end{enumerate}

\end{proposition}

\begin{proof}
The claims (1) and (2) follow from Theorem \ref{thm:splittableInjectivity}.  For (3) we note that the link splitting map $\FT_n^{k+l} \rightarrow \one$ is the tensor product, over $\one\in \YC(\DC_n)$, of the link splitting maps $\FT_n^k\rightarrow \one$ and $\FT_n^l\rightarrow \one$.  Taking the tensor product of a class $f\in\Hom_{\YC(\DC_n)}^{\Z\times \Z\times \Z}(\one,\FT_n^k)$ and a class $g\in \Hom_{\YC(\DC_n)}^{\Z\times \Z\times \Z}(\one,\FT_n^l)$ and post-composing with the splitting map $\FT_n^{k+l}\rightarrow \one$ naturally yields the product of $\Psi^{\otimes k}\circ f$ and $\Psi^{\otimes l}\circ g$ in $\Q[\xx,\yy,\ttheta]$.  This proves (3). 
\end{proof}

Motivated by the work of Haiman \cite{haiman2001hilbert}, we introduce another family of ideals.  
\begin{definition}
Let $J_n\subset \Q[\xx,\yy]$ and $\AJ_n\subset  \Q[\xx,\yy,\ttheta]$ denote the ideals generated by the alternating (i.e.~anti-symmetric) polynomials with respect to the diagonal action of $S_n$.
\end{definition}

Note that $J_n$ is the $a$-degree zero component of $\AJ_n$.  

\begin{theorem}[\cite{haiman2001hilbert}]
\label{th:haiman}
\begin{enumerate}
\item[(a)]For all $k>0$ the ideal $J_n^{k}$ is free as a $\Q[y_1,\ldots,y_n]$-module.
\item[(b)]For all $k>0$ one has 
\[
J_n^k=\bigcap_{i\neq j}(x_i-x_j,y_i-y_j)^{k}.
\]
\end{enumerate}
\end{theorem}
In fact, Haiman showed that (a) implies (b), and the really hard part is the proof of (a). 
The proof of the implication (a)$\Rightarrow$ (b) (see \cite[Proposition 6.1]{haiman1999macdonald}) uses the following result from commutative algebra, which we will need.  

\begin{lemma}
Let $M$ be a free module over $\C[\yy]$, and let $\MC$ denote the corresponding coherent sheaf on $\C^n$.
Let $Z=\{y_1=\ldots=y_n\}$ and $U=\C^n\setminus Z$. Then for $n\ge 3$ one has
\[
M=H^0(\C^n,\MC)=H^0(U,\MC).
\]
\end{lemma}

\begin{proof}
Since $M$ is free, the sequence $y_1-y_2,\ldots,y_1-y_n$  is regular for $M$. Therefore (see e.~g.~ \cite[Theorem 5.8]{huneke}) the local cohomology $H^i_Z(\MC)$ vanish for $i<n-1$. Since $n\ge 3$, we have $H^0_Z(\MC)=H^1_Z(\MC)=0$. Now the long exact sequence
\[
0=H^0_Z(\MC)\to H^0(\C^n,\MC)\to H^0(U,\MC)\to H^1_Z(\MC)=0
\]
implies the desired isomorphism.
\end{proof}

\begin{corollary}
\label{cor: extension to diagonal}
In the assumptions of the lemma, suppose that $M\subset N$ for some module $N$ and the restrictions of the corresponding sheaves to $U$ coincide: $\MC|_U=\NC|_U$. Then $M=N$ and, in particular, $N$ is free. 
\end{corollary}

\begin{proof}
We have $N\subset H^0(U,\NC)=H^0(U,\MC)=M$. Since $M\subset N$, we have $M=N$.
\end{proof}

We are ready to prove the following result, motivated by some of the conjectures in \cite{GNR}.

\begin{theorem}
\label{th: homology full twist}
For all $k,n>0$ one has 
\[
\HY(\FT_n^{k}) \cong \IY(\FT_n^k) \cong \bigcap_{i\neq j}(x_i-x_j,y_i-y_j,\theta_i-\theta_j)^k\subset \Q[\xx,\yy,\ttheta].
\]
\end{theorem}

\begin{proof}
Without loss of generality we can change the field of definition from $\Q$ to $\C$. 
The isomorphism $\HY(\FT_n^{k})\cong \IY(\FT_n^k)$ holds by Proposition \ref{prop:FTbasicProps}.    We can write $\FT_n=\prod_{i<j}A_{ij}$, where the ordering of the factors is
 \[
 \FT_n = (A_{12}) (A_{13}A_{23}) \cdots (A_{1,n}\cdots A_{2,n}).
 \]
By Proposition \ref{prop:purebraidProps} we have an inclusion:
\begin{equation}
\label{eq:full twist inclusion}
\IY(\FT_n^k) \subset \bigcap_{i\neq j}(x_i-x_j,y_i-y_j,\theta_i-\theta_j)^k.
\end{equation}
It is sufficient to prove that this inclusion is an isomorphism after localization at every point $\yy=(\nu_1,\ldots,\nu_n)\in \C^n$.
We prove this by induction in $n$. For $n=1$ it is obvious, for $n=2$ it follows from Proposition \ref{prop:powersOfFT2}.

Suppose that $n\ge 3$. For all $\nu_i$ distinct the inclusion \eqref{eq:full twist inclusion} is an isomorphism by Corollary \ref{cor:linksplittingproperty}.  For arbitrary, $\nu$, let $\{1,\ldots,n\}=\Pi_1\sqcup \Pi_2\sqcup \ldots \sqcup \Pi_r$ be the set partition such that $\nu_i=\nu_j$ if and only if $i,j$ are in the same block $\Pi_a$.  

Let $\Psi:\FT_n\rightarrow \one$ be the splitting map.  In the course of proving Theorem \ref{thm:BraidGpAction}, we showed that conjugation by any braid $\FY(\b)$ fixes $\Psi$ up to homotopy.  Since conjugation also permutes the labels on the components, we may as well assume that
\[
\nu_1=\cdots=\nu_{n_1},\qquad \nu_{n_1+1}=\cdots=\nu_{n_1+n_2},\qquad\cdots\qquad \nu_{n_1+\cdots + n_{r-1}+1}=\cdots = \nu_{n_1+\cdots + n_{r}}
\]
for some integers $n_1,\ldots,n_r$ with $\sum_a n_a=n$.  Then by changing crossings between the components in different blocks, we can factor the map $\Psi^{\otimes k}$ as
\[
\CY(\FT_n^k)\to \CY(\FT_{n_1}^k)\otimes_\Q \cdots \otimes_\Q \CY(\FT_{n_r}^k)  \to \Q[\xx,\yy,\ttheta]
\]
The first arrow is an equivalence after localization at all $y_i-y_j$ for $i$ and $j$ in different blocks. Therefore at such $\nu$ one has
\[
\HY(\FT_n^k,\C_{\nu})\cong \bigotimes_a \HY(\FT_{n_a}^k,\C_{\nu}).
\]
A similar factorization holds for the right hand side of \eqref{eq:full twist inclusion}, so we conclude that \eqref{eq:full twist inclusion}
is an isomorphism for all $\nu$ outside of the ``diagonal'' $Z=\{\nu_1=\ldots=\nu_n\}$. 

Finally, by Theorem \ref{th: ft parity} and Theorem \ref{th:flat vs collapsible} $\HY(\FT_n^k,\C)$ is a free $\C[\yy]$-module.
By Corollary \ref{cor: extension to diagonal} the inclusion \eqref{eq:full twist inclusion} is an isomorphism.
\end{proof}

\begin{corollary}
\label{cor: full twist a=0}
One has $\HY^0(\FT_n^{k})\cong J_n^k$ and $\HKR^0(\FT_n^k)\cong J_n^k/\yy J_n^{k}.$
\end{corollary}

\begin{proof}
The odd variables $\theta_i$ have $a$-degree 1, so in $a$-degree 0 we just have $x_i,y_i$ and
\[
\HY^0(\FT_n^{k})\cong \bigcap_{i\neq j}(x_i-x_j,y_i-y_j).
\]
By Theorem \ref{th:haiman} this ideal is isomorphic to $J_n^k$ and free as $\Q[\yy]$-module.
Therefore 
\[
\HKR^0(\FT_n^k)=\HY(\FT_n^k)/\yy\HY(\FT_n^k)\cong J_n^k/\yy J_n^{k}.
\]
\end{proof}

\begin{remark}
In \S \ref{subsec:polygraphs} below we show that $\HY(\FT_n^{k})\cong \AJ_n^{ k}$, hence $\HKR(\FT_n^{k})\cong \AJ_n^{k}/(\yy)\AJ_n^k$.
\end{remark}


\section{Hilbert scheme of points}
\label{sec:Hilbert}

In this section we collect several facts on the geometry of the Hilbert scheme of points, mostly due to Haiman \cite{haiman2001hilbert,haiman2002vanishing}.  We extend some of Haiman's work to give a description of $\JC_n^k\subset \Q[\xx,\yy,\ttheta]$ and relate to the homology of full twists.

\subsection{Geometry}

The Hilbert scheme of $n$ points on $\C^2$ is defined as the moduli space of ideals $I\subset \C[x,y]$ such that $\dim \C[x,y]/I=n$.
It is a smooth quasiprojective complex algebraic variety of dimension $2n$. There is a natural Hilbert-Chow map $\Hilb^n(\C^2)\to \Sym^n(\C^2)$ 
which sends an ideal to its support. 

The {\em isospectral Hilbert scheme} $X_n$ is defined as the reduced fibered product:
\begin{center}
\begin{tikzcd}
X_n \arrow[r]\arrow[d] & (\C^2)^n\arrow[d]\\
\Hilb^n(\C^2) \arrow[r] & \Sym^n(\C^2).\\
\end{tikzcd}
\end{center}
Both $\Hilb^n(\C^2)$ and $X_n$ can be explicitly described as blowups of certain ideals.
As above, let $A^{\C}_n\subset \C[\xx,\yy]$ denote the space of antisymmetric polynomials with respect to the diagonal action of $S_n$,
and let $J^{\C}_n$ be the ideal generated by $A^{\C}_n$. Then the following holds:
\begin{equation}
\Hilb^n(\C^2)=\Proj\left(\bigoplus_{k=0}^{\infty} (A^{\C}_n)^k\right),\ X_n=\Proj\left(\bigoplus_{k=0}^{\infty}(J_n^{\C})^k\right).
\end{equation}
As a consequence, both $X_n$ and $\Hilb^n(\C^2)$ carry the natural line bundle $\OC(1)$ and
\[
H^0(\Hilb^n(\C^2),\OC(k))=(A_n^{\C})^k,\ H^0(X_n,\OC(k))=(J_n^{\C})^k.
\]
More generally, both $X_n$ and $\Hilb^n(\C^2)$ carry the tautological vector bundle $T$ with fiber $\C[x,y]/I$.
One can prove that $\OC(1)=\det(T)$. The following result is a key theorem of \cite{haiman2002vanishing}:

\begin{theorem}
For all $l>0$ one has
\[
H^i(X_n,T^l)=0\ \text{for}\ i>0,\ H^0(X_n,T^l)=R(n,l),
\]
where $R(n,l)$ is the {\em polygraph ring} defined in \S \ref{subsec:polygraphs}.
\end{theorem}

\begin{corollary}
\label{cor:wedges}
For all $l,m>0$ one has
\[
H^i(X_n,\wedge^l(T)\otimes \OC(m))=0\ \text{for}\ i>0,\ H^0(X_n,\wedge^l(T)\otimes \OC(m))=R(n,l+mn)^{\rho},
\]
where $\rho$ is the sign representation of $S_l\times (S_n)^m$ on $R(n,l+mn)$.
\end{corollary}

\subsection{Polygraphs}
\label{subsec:polygraphs}

The polygraph ring $R(n,l)$ is defined as follows. Consider an $(n+l)$-tuple of points on $\C^2$:
\[
P_1=(x_1,y_1),\ldots,P_n=(x_n,y_n),Q_1=(a_1,b_1),\ldots,Q_l=(a_l,b_l).
\]
We define the {\em polygraph} $Z(n,l)\subset \C^{2n+2l}$ by the equation
\[
Z(n,l)=\{(P_1,\ldots,P_n,Q_1,\ldots Q_l)\ :\ Q_i\in \{P_1,\ldots,P_n\}\ \text{for all}\ i\}.
\]
Observe that $Z(n,l)$ is a union of $n^l$ subspaces parametrized by the functions $f:\{1,\ldots,l\}\to \{1,\ldots,n\}$:
\[
Z(n,l)=\bigcup_{f}W_f,\qquad W_f=\{Q_i=P_{f(i)}\ \text{for all}\ i\}.
\]
Each $W_f$ has dimension $2n$ and coordinates $x_1,\ldots,x_n,y_1,\ldots,y_n$.
Let $I(n,l)$ be the defining ideal of $Z(n,l)$.  The polygraph ring is defined as a ring of functions on $Z(n,l)$:
\[
R(n,l)=\C[x_1,\ldots,x_n,y_1,\ldots,y_n,a_1,\ldots,a_l,b_1,\ldots,b_l]/I(n,l).
\]
\begin{example}
For $l=1$ we get 
\[
I(n,1)=\bigcap_{i=1}^{n}(x_i-a_1,y_i-b_1)\subset \C[x_1,\ldots,x_n,y_1,\ldots,y_n,a_1,b_1].
\]
By Proposition \ref{prop:jucys murphy} one has 
\[
I(n,1)\cong \HY(\LC_{n+1},\C),
\]
where $\LC_{n+1}$ is the Jucys-Murphy braid and we identify $x_{n+1}=a_1,y_{n+1}=b_1$.
\end{example}

\begin{theorem}(\cite{haiman2001hilbert})
\label{th: polygraph free}
The ring $R(n,l)$ is a free $\C[\yy]$-module for all $n$ and $l$.
\end{theorem}

Given a function $g\in R(n,l)$, one can restrict it to each subspace $W_f\cong (\C^2)^n$ and obtain 
a collection of functions in $\C[\xx,\yy]$:
\begin{equation}
R(n,l)\rightarrow \bigoplus_{f}\C[\xx,\yy],\qquad g\mapsto (g|_{W_f}).
\end{equation}
By construction, this map is injective. 

As above, let $\AJ_n^{\C}=\AJ_n\otimes_{\Q}\C$ denote the ideal in $\C[\xx,\yy,\thth]$ generated by the antisymmetric polynomials.

\begin{lemma}
\label{lem:generating set for AJ}
The ideal $\AJ_n^{\C}$ is generated as a $\C[\xx,\yy]$-module by polynomials of the form
\begin{equation}
\label{generating set for AJ}
\Alt\left(g(x_{k+1},\ldots,x_{n},y_{k+1},\ldots,y_{n})\theta_1\ldots\theta_k\right)
\end{equation}
for all $k$ and all antisymmetric polynomials $g(x_1,\ldots,x_{n-k},y_1,\ldots,y_{n-k})\subset A_{n-k}$.
\end{lemma}

\begin{proof}\MH{come back}
For each subset $S\subset \{1,\ldots,n\}$, let $\theta_S = \theta_{i_1}\cdots \theta_{i_k}$ where $S=\{i_1<\cdots <i_k\}$.  The antisymmetric component of $\C[\xx,\yy,\thth]$ is spanned by the antisymmetrizations of monomials.  In fact we restrict to monomials of the form $m(\xx,\yy)\theta_1\cdots \theta_k$.  We can first antisymmetrize separately in variables from $S=\{1,\ldots,k\}$ and from its complement $S^c$. The result is a polynomial of the form $g   h \theta_S $, where $g$ is an antisymmetric element of $\C[x_i,y_i]_{i\not\in S}$ and $h$ is a symmetric element of $\C[x_i,y_i]_{i\in S}$.   One can check that the space of such polynomials is also spanned by the polynomials of the form
$g h' \theta_S$, where now $h'\in \C[\xx,\yy]^{S_n}$ is symmetric in all variables.  Now we antisymmetrize with respect to the whole group $S_n$ and obtain
\[
\Alt(m(\xx,\yy)\theta_S)=\Alt(g h'\theta_S )=h'\Alt(g \theta_S)=\]
\[
\pm h'\Alt(g(x_{k+1},\ldots,x_n,y_{k+1},\ldots,y_n)\theta_1\cdots\theta_k).
\]
This shows that the space of antisymmetric functions in $\C[\xx,\yy,\thth]$ is generated over $\C[\xx,\yy]$ by \eqref{generating set for AJ}. 

Finally, one can check that $\C[\xx,\yy]$--submodule generated by \eqref{generating set for AJ} is stable with respect to the multiplication by $\theta_i$. Indeed, let $G$ be an antisymmetric polynomial of the form \eqref{generating set for AJ}. We can write
\[
G=\sum_{a,b} x_i^{a}y_i^bg_{a,b}+\theta_ig',
\]
where $g_{a,b}$ and $g'$ do not depend on $x_i,y_i$ or $\theta_i$. 
Then 
\[
\theta_iG=\theta_i\sum_{a,b} x_i^{a}y_i^{b}g_{a,b}=\sum_{a,b}x_i^{a}y_i^{b}\Alt(\theta_i g_{a,b}).
\]
Indeed, for $j\neq i$ the right hand side has a summand
\[
\theta_j\sum_{a,b}x_i^{a}y_i^{b}\cdot s_{ij}(g_{a,b})=\theta_{j}s_{ij}\left[\sum_{a,b}x_j^{a}y_j^{b}g_{a,b}\right]=
\]
\[
\theta_js_{ij}\left[G(x_i=x_j,y_i=y_j,\theta_i=\theta_j)\right]=0,
\]
where $s_{ij}=(i\ j)$.  
\end{proof}

\begin{lemma}
\label{lem:antizymmetrization polygraph}
Let $\rho=\sgn\boxtimes\cdots \boxtimes \sgn$ denote the sign representation of $S_l\times (S_n)^k$. Then
\[
\bigoplus_{l=0}^{n} R(n,l+kn)^{\rho}\cong \AJ^{\C}_{n}\cdot (J_n^{\C})^k.
\]
\end{lemma}

\begin{proof}
Let us fix some $l$.  Let $\mathbf{a}=(a_1,\ldots,a_{l+kn})$, $\mathbf{b}=(b_1,\ldots,b_{l+kn})$ are sets of formal variables.  The left hand side is the $\C[\xx,\yy]$-submodule  of $R(n,l+kn)\cong \C[\xx,\yy,\mathbf{a},\mathbf{b}]/I(n,l+kn)$ generated by the products of antisymmetric functions 
\[
G=g_1(a_1,\ldots,a_l,b_1,\ldots,b_l)\cdot \prod_{j=0}^{k-1} g_j(a_{l+jn+1},\ldots,a_{l+(j+1)n},b_{l+jn+1},\ldots,b_{l+(j+1)n}).
\]
Let us describe the restriction of $G$ to the subset $W_f$. If $f(i)=f(j)$ for some $1\le i<j\le l$ then the restriction vanishes. 
Similarly, the values $f(l+jn+1),\ldots,f(l+j(n+1))$ should be  all distinct for each $0\le j\le k-1$. If $f([1,l])=S$, we get
\[
G|_{W_f}=\pm g_1(S)\cdot \prod_{j=0}^{k-1} g_j(x_1,\ldots,x_n,y_1,\ldots,y_n).
\]
By Lemma \ref{lem:generating set for AJ} this collection of restrictions can be compared to the generator of $\AJ\cdot (J_n)^k$ of the form 
\[
\Alt(g_1(S)\theta_{\overline{S}})\cdot \prod_{j=0}^{k-1} g_j(x_1,\ldots,x_n,y_1,\ldots,y_n).
\]
\end{proof}


\begin{proposition}
\label{prop: AJ vs intersection}
For all $k>0$ one has
\[
\AJ^{\C}_{n}\cdot (J_n^{\C})^{k-1}=\bigcap_{i\neq j}(x_i-x_j,y_i-y_j,\theta_i-\theta_j)^{k}.
\]
\end{proposition}

\begin{proof}
Let $g\in \C[\xx,\yy,\thth]$ be an antisymmetric polynomial. Let us prove that 
$g\in (x_i-x_j,y_i-y_j,\theta_i-\theta_j)$ for all $i$ and $j$. 
The terms with $\theta_i\theta_j$ are divisible by $\theta_i-\theta_j$. The remaining terms have the form
\[
h\theta_i-(i\ j)h\theta_j=h(\theta_i-\theta_j)+(h-(i\ j)h)\theta_j.
\]
Therefore
\[
\AJ_n\subset (x_i-x_j,y_i-y_j,\theta_i-\theta_j) 
\]
for all $i$ and $j$, and
\[
\AJ^{\C}_{n}\cdot (J_n^{\C})^{k-1}\subset (x_i-x_j,y_i-y_j,\theta_i-\theta_j)^{k}.
\]
By Lemma \ref{lem:antizymmetrization polygraph} and  Theorem \ref{th: polygraph free} the left hand side is a free $\C[y]$-module. Similarly to the proof of Theorem \ref{th: homology full twist} (see also \cite[Proposition 6.1]{haiman1999macdonald}) we check by induction that the ideals coincide outside $\{y_1=\ldots=y_n\}$ and by Corollary \ref{cor: extension to diagonal}  they coincide on the diagonal as well.
\end{proof}

Note that setting $k=1$ in the above and changing $\C$ to $\Q$ yields  $\AJ_n = \bigcap_{i\neq j}(x_i-x_j,y_i-y_j,\theta_i-\theta_j)$.  Thus $\AJ_n^k \subset \bigcap_{i\neq j}(x_i-x_j,y_i-y_j,\theta_i-\theta_j)^k$.  We also have the opposite containment by the above lemma: $\bigcap_{i\neq j}(x_i-x_j,y_i-y_j,\theta_i-\theta_j)^k=\AJ_{n}\cdot (J_n)^{k-1}\subset  \AJ_n^k$.  Thus we have the following.

\begin{corollary}\label{cor:AJk is intersection}
We have $\AJ_n^k =\AJ_nJ_n^{k-1} = \bigcap_{i\neq j}(x_i-x_j,y_i-y_j,\theta_i-\theta_j)^k  \subset \Q[\xx,\yy,\ttheta]$ for all $n,k\geq 0$.\qed
\end{corollary}
Combining this with the computation from Theorem \ref{th: homology full twist} now gives us one of our main results.
\begin{corollary}
\label{cor: full twist from AJ}
For all $k\ge 0$ one has $\HY(\FT_n^k)\cong\AJ_n^{k}$.\qed
\end{corollary}

\begin{corollary}
\label{cor: full twist from isospectral}
The $y$-ified homology of $\FT_n^k$ (with complex coefficients) is isomorphic to the space of global sections of the vector bundle $\Lambda(T^*)\otimes \OC(k)$
on the isospectral Hilbert scheme $X_n$, where $\Lambda=\bigoplus_i \Lambda^i$ denotes the exterior algebra.
\end{corollary}

 \begin{proof}
First observe that $\Lambda^i(T^*)\otimes \OC(k)=\wedge^{n-i}(T)\otimes \OC(k-1)$, so
\[
H^0(X_n,\Lambda(T^*)\otimes \OC(k))\simeq  H^0(X_n,\Lambda(T)\otimes \OC(k-1)).
\]
By Corollary \ref{cor:wedges} and Lemma \ref{lem:antizymmetrization polygraph} this is isomorphic to 
\[
R(n,l+(k-1)n)^{\rho}\cong \AJ^{\C}_{n}\cdot (J_n^{\C})^{k-1}.
\]
By Corollary \ref{cor:AJk is intersection} this is isomorphic to $\HY(\FT_n^k,\C)$. 
\end{proof}

\subsection{Combinatorics}
\label{sec:combinatorics}

In this appendix we use the combinatorics of Macdonald polynomials to compute the triply graded character of the powers of the full twist.
Since we do not expect a reader to be familiar with this subject and feel that it is quite tangential to the subject of this paper, we keep the exposition brief and refer the reader to \cite{haiman2002vanishing,GN,GNR, Wilson} for definitions and more details. 

Let $\Lambda_{q,t}$ denote the ring of symmetric functions with coefficients in $\Q(q,t)$.  Let $\Ht_{\lambda}\in \Lambda_{q,t}$ denote the modified Macdonald polynomial corresponding to the Young diagram $\lambda$ \cite{haiman2002vanishing}, and let $\nabla$ denote the Bergeron-Garsia operator such that
\[
\nabla \Ht_{\lambda}=\prod_{\square \in \lambda} \Ht_{\lambda}.
\]

Let $\ev:\Lambda_{q,t}\rightarrow \Q(q,t)[a]$ denote the algebra homomorphism sending the elementary symmetric functions $e_n$ to $1+a$ for all $a$.  Alternatively, 
\[
\ev(f) = \sum_i (f, h_ke_{n-k})a^k,
\]
where $h_k$ is the complete homogeneous symmetric function and $(\  , \ )$ is the Hall inner product.  In terms of plethystic notation, $\ev(f) = \omega(f)[1-a]|_{a\to -a}$ where $\omega$ is the classical involution sending $e_n\leftrightarrow h_n$ and $p_n\mapsto (-1)^{n-1}p_n$.
\begin{theorem}
\label{th:poincare nabla}
The Poincar\'e series of $\HY(\FT_n^k)$ equals $\ev(\nabla ^k p_1^n)$.
\end{theorem}

\begin{proof}
By Corollary \ref{cor: full twist from isospectral} we get
\[
\HY(\FT_n^k,\C)\cong H^0(X_n,\Lambda(T^*)\otimes \OC(k)),
\]
where $X_n$ is the isospectral Hilbert scheme and $T$ is the tautological bundle. 
By \cite{haiman2002vanishing}, the higher homology $H^i(X_n,\Lambda(T^*)\otimes \OC(k))$ ($i>0$) vanishes, so the Poincar\'e series of $\HY(\FT_n^k)$
equals the Euler characteristic
\[
\PC(\HY(\FT_n^k))=\chi(X_n,\Lambda(T^*)\otimes \OC(k))=\chi(\Hilb^n(\C^2),\pi_*(\Lambda(T^*)\otimes \OC(k))).
\]
The pushforward of the structure sheaf of $X_n$ to the Hilbert scheme is the celebrated {\em Procesi bundle} $P$. 
By  projection formula, we get
\[
\PC(\HY(\FT_n^k))=\chi(\Hilb^n(\C^2),P\otimes\Lambda(T^*)\otimes \OC(k)).
\]
One can identify the equivariant $K$-theory of $\Hilb^n(\C^2)$ with the space of symmetric functions of degree $n$.
The Procesi bundle $P$ corresponds to $p_1^n$, and the operator of multiplication by $\OC(1)$ corresponds to $\nabla$,
so $P\otimes \OC(k)$ corresponds to $\nabla^k(p_1^n)$. 

Finally, the functional $\FC\mapsto \chi(\Hilb^n(\C^2),\FC\otimes \Lambda(T^*))$ corresponds to $f\mapsto \ev(f)$.  See e.~g.~ \cite{GN}.
\end{proof}

Following \cite{Hog17a}, we introduce a family of rational functions $f_{v}(q,t,a)$ labeled by sequences $v\in \{0,\ldots,k\}^n$ and satisfying the following recursion relations:
\begin{equation}
f_{\emptyset}=1,\ f_{0,v} = (t^{\sharp\{ i : v_i<k\} }+ a)f_v,
\end{equation}
\begin{equation}
f_{j,v}=t^{\sharp\{ i : v_i<j\}}f_{v,j-1}\ \text{for}\ 1\le j\le k-1,
\end{equation}
\begin{equation}
f_{k,v}=f_{v,k-1}+qf_{v,k}.
\end{equation}

We refer to \cite{Hog17a} for more details on functions $f_{v}$ and their connections to link invariants.

\begin{theorem}
\label{th: full magic formula}
The following identity holds for all $n,k\ge 1$:
\[
f_{k,\ldots,k}=(1-t)^n\ev(\nabla^k p_1^n).
\]
\end{theorem}

For $k=1$ this identity was conjectured in \cite[Conjecture 1.15]{EH} and \cite[eq. (48)]{Wilson}.

\begin{proof}
By \cite[Theorem 1.3]{Hog17a} the trigraded Poincar\'e series of $\HKR(\FT_n^k)$ equals $f_{k,\ldots,k}$.
Since $\FT_n^k$ is $y$-flat, the Poincar\'e polynomial of $\HY(\FT_n^k)$ equals $f_{k,\ldots,k}/(1-t)^n.$
On the other hand, by Theorem \ref{th:poincare nabla} it equals $\ev(\nabla ^k p_1^n).$
\end{proof}

Wilson in \cite[Conjecture 4.1]{Wilson} also conjectured an explicit combinatorial formula for $\nabla p_1^n$.  The above description of $\ev(\nabla p_1^n)$ recovers only the coefficients of hook Schur functions in $\nabla p_1^n$.  To our knowledge, Wilson's conjecture remains open.  

\section{Symmetry properties}
\label{sec:symmetry}


\begin{definition}\label{def:rev}
Let $\rev$ denote the automorphism of $\Z\times \Z\times \Z$ sending
\[
i(1,0,0)+j(-2,1,0)+k(-1,0,1) \leftrightarrow k(1,0,0)+j(-2,1,0)+i(-1,0,1).
\]
If $V$ is a triply graded vector space, we let $V_{\rev}$ be the triply graded vector space with $V_{\rev}^{(i,j,k)}=V^{\rev(i,j,k)}$.  Let $\tau:\Q[\xx,\yy,\ttheta]\rightarrow \Q[\xx,\yy,\ttheta]_{\rev}$ be the algebra automorphism sending $x_i\leftrightarrow y_i$ and $\theta_i\mapsto \theta_i$.
\end{definition}

\begin{conjecture}\label{conj:Hsymmetry}
For each oriented link $L\in \R^3$ we have $\HY(L)\rightarrow \HY(L)_{\rev}$ which exchanges the actions of $x_c$ and $y_c$ and is equivariant with respect to the action of $\theta_c$.
\end{conjecture}

\begin{example}
By Proposition \ref{prop:jucys murphy} and Theorem \ref{th: homology full twist} the conjecture is true for Jucys-Murphy braids and for the powers of the full twist.
\end{example}

This symmetry would be a categorical analogue of the $Q\leftrightarrow -Q\inv$ symmetry in the HOMFLY-PT polynomial.
For reduced Khovanov-Rozansky homologies of knots this symmetry was conjectured by Gukov, Dunfield and Rasmussen in \cite{DGR}.
For $y$-flat links with several components we get the following result.

\begin{proposition}
Suppose that Conjecture \ref{conj:categoricalSymmetry} holds for a $y$-flat link $L$ with $r$ components. Then 
\[
(1-q)^r\PC_{L}(q,t,a) = (1-t)^r\PC_{L}(t,q,a)
\]
where $\PC_{L}$ is the Poincar\'e polynomial of $\HKR(L)$.
\end{proposition}

\begin{proof}
Since the link is $y$-flat, the Poincar\'e polynomial of $\HY(L)$ equals $\PC_{L}(q,t,a)/(1-t)^r$, and by Conjecture \ref{conj:Hsymmetry} we get
\[
\PC_{L}(q,t,a)/(1-t)^r=\PC_{L}(t,q,a)/(1-q)^r.
\]
\end{proof}

Moreover, we expect that the symmetry can be seen on the level of complexes. 

\begin{conjecture}\label{conj:categoricalSymmetry}
There exists a triangulated auto-equivalence of $\YC(\SBim_n)$ with the following properties:
\begin{enumerate}\setlength{\itemsep}{3pt}
\item $\tau(\FY(\b))\simeq \FY(\b)$.  That is $\tau$ fixes Rouquier complexes.
\item $\tau(\CB(1))\cong \tau(\CB)(-1)[1]$.
\item $\tau(\CB[1])\cong \tau(\CB)[1]$. 
\item $\tau(\CB\otimes \DB) \cong \tau(\CB)\otimes\tau(\DB)$.
\item the action of $\tau$ on $\End^{\Z\times \Z\times \Z}(\one)$ agrees with the definition of $\tau$ in Definition \ref{def:rev}.
\end{enumerate}
The statements (1), (2), (3) are supposed hold for all $\CB,\DB\in \YC(\SBim_n)$.  We also conjecture that this autoequivalence extends to the derived version of $\SBim_n$, and behaves as follows with respect to the Hochschild degree shift:
\[
\tau(\CB(-2,1)) \cong \tau(\CB)(-2,1).
\]
\end{conjecture}
  
We plan to investigate this symmetry and its connection to various Koszul dualities in \cite{AMRW} in a future work.

\appendix
\section{Curved complexes and homotopy lemmas}
\label{sec:appendix}

Let $\CC$ be an additive category, and let $Z\in \End(\Id_{\CC})$ be an element in the center of the category (natural endomorphism of the identity functor).  That is, for each $A\in \CC$ we have an endomorphism $Z_A:A\rightarrow A$, and $Z_C$ commute with all morphisms in $\CC$.  A \emph{curved complex with curvature $Z$} is a pair $(C,\Delta)$ with $C\in \CC$ and $\Delta^2=Z_C$.  A morphism of curved complexes is a morphism in $\CC$ which commutes with the differentials $\Delta$.  Two morphisms, $f,g:(C,\Delta_C)\rightarrow (D,\Delta_D)$ are \emph{homotopic} if there exists $h\in \Hom_{\CC}(C,D)$ such that $f-g=\Delta h+ h\Delta$.   A curved complex is \emph{contractible} if its identity endomorphism is homotopic to zero.

The homotopy category of curved complexes with curvature $Z\in \End(\Id_{\CC})$ will be denoted by $\KC(\Fac(Z))$.  We let $\FC:\KC(\Fac(Z))\rightarrow \CC$ denote the functor which forgets the differential $\Delta$.

\begin{remark}
One may refine the above definitions to include extra data, such as gradings or linearity with respect to the action of some ground ring.  In particular, chain complexes and matrix factorizations arise as special cases of curved complexes.
\end{remark}

\subsection{Mapping cones}
\label{subsec:cones}
If $(A,\Delta_A)$ is a curved complex, then we let $A[1]$ denote the curved complex which equals $A$ with $\Delta_{A[1]}=-\Delta_A$.  Let $(A,\Delta_A)$ and $(B,\Delta_B)$ be  curved complexes in $\CC$ with curvature $Z$.  If $f:A\rightarrow B$ is a morphism of curved complexes, then the \emph{mapping cone} of $f$ is the curved complex
\[
\Cone(f) = \left(A[1]\oplus B , w , \smMatrix{-\Delta_{A} & 0 \\ f & \Delta_{B}}\right).
\]

\begin{lemma}\label{lemma:coneDetectsHomotopyEquiv_appendix}
A morphism $f:(A,d_A)\rightarrow (B,d_B)$ of curved complexes is a homotopy equivalence if and only if $\Cone(f)$ is contractible.
\end{lemma}

\begin{proof}
Throughout we suppress any and all gradings if there are any.  Assume that $\FC(\Cone(f))=A\oplus B$ has an endomorphism of the form
\[
H = \begin{bmatrix}-h & g \\ m & k\end{bmatrix}
\]
Any easy computation shows that
\[
\sqmatrix{-d_A & 0 \\ f & d_B}\sqmatrix{-h & g \\ m & k} + \sqmatrix{-h & g \\ m & k}\sqmatrix{-d_A & 0 \\ f & d_B} = \sqmatrix{dh+hd + gf & -dg + gd \\ dm-md-(fh-kf) & dk+kd+fg}
\]
This equals the identity matrix if and only if
\begin{eqnarray}
 [d,g] &=& 0 \label{eq-dg0}\\
 {}[d,h] &=& \Id_A - gf  \label{eq-dh} \\
 {}[d,k]  &=& \Id_B - fg \label{eq-dk}\\
  {}[d,m] &=& kf - fh \label{eq-dm}
\end{eqnarray}

The first three of these relations show that if $\Cone(f)$ is contractible, then $f$ is a homotopy equivalence.  Suppose conversely that $f$ is a homotopy equivalence.  Let $g\in \Hom^0(B,A)$, $h\in \End^{-1}(A)$, and $k\in \End^{-1}(B)$ be such that
\[
[d,g]=0 \hskip.5in \Id_B-fg = [d,k] \hskip.5in Id_A - gf = [d,h]
\]
so that $g$ and $f$ are inverse equivalences.  Set $z:=kf - fh\in \Hom^{-1}(A,B)$, and consider the matrix
\[
H' = \sqmatrix{-(h + gz) & g \\ kz & k} =: \sqmatrix{-h' & g \\ m & k} 
\]
An easy computation shows that $H'$ is a null-homotopy for $\Cone(f)$.  Indeed, note that $z$ is a cycle, since
\[
[d,z] = [d,k]f - f[d,h] = (\Id - fg)f - f(\Id-gf) = 0
\]
Thus the top left corner satisfies the relation $[d,h+gz] = [d,h] = \Id - gf$.  The only remaining relation to check is (\ref{eq-dm}).   The bottom left corner satisfies
\[
[d,kz] = [d,k]z = z - fgz 
\]
On the other hand,
\[
kf - fh' = kf - fh -fgz = z-fgz
\]
This shows that $[d,m]=kf - fh'$, and completes the proof.
\end{proof}

\subsection{Gaussian elimination}
\label{subsec:gauss}
Gaussian elimination is traditionally stated as follows.  Let $\AC$ be an additive category, and suppose we have a complex of objects in $\AC$ of the form
\[
\begin{diagram}
\cdots &\rTo & A^{k-1} &\rTo^{\sqmatrix{\a \\ \b}} & A^{k}\oplus B &\rTo^{\sqmatrix{\gamma & \d \\ \e& \zeta}} & A^{k+1}\oplus B'&\rTo^{\sqmatrix{\eta & \theta}} & A^{k+2} &\rTo &\cdots 
\end{diagram}
\]
where $\zeta:B\rightarrow B'$ is an isomorphism.  Then this complex splits as a direct sum of
\[
\begin{diagram}
\cdots &\rTo & A^{k-1} &\rTo^{\a} & A^{k} &\rTo^{\gamma - \e\zeta\inv \d} & A^{k+1} &\rTo^{\eta} & A^{k+2} &\rTo &\cdots 
\end{diagram}
\]
plus the contractible complex
\[
\begin{diagram}[small]
0 &\rTo & B &\rTo^{\zeta} & B' &\rTo & 0.
\end{diagram}
\]

Now we consider a generalization of this idea, suppressing all gradings for the time being.  Let $\CC$ be an additive category and $Z\in \End(\Id_{\CC})$ an element in the center.   Let $A,B\in \CC$, and suppose we have a curved complex $(A\oplus B, \Delta)$.   We will write $\Delta$ in terms of its components:
\[
\Delta = \sqmatrix{\Delta_{AA} & \Delta_{AB}\\ \Delta_{BA} & \Delta_{BB}},\qquad \qquad \Delta^2 = \sqmatrix{Z_A&0\\0&Z_B}.
\]

Suppose there exists $h\in \End(B)$ such that
\begin{subequations}
\begin{equation}\label{eq:Bhomotopy}
\Delta_{BB}\circ h + h\circ \Delta_{BB}=\Id_B
\end{equation}
\begin{equation}\label{eq:homotopysquared}
h^2=0.
\end{equation}
\end{subequations}
Assume also that $\Delta_{BB}^2=0$.  In other words, $(B,\Delta_{BB})$ is a contractible complex in the usual sense.

\begin{remark}
Equation \eqref{eq:homotopysquared} can always be arranged: if $h'$ satisfies \eqref{eq:Bhomotopy}, and $\Delta_{BB}^2=0$, then $h:=h'\circ  \Delta_{BB}\circ  h'$ satisfies \eqref{eq:Bhomotopy} and \eqref{eq:homotopysquared}.
\end{remark}

\begin{lemma}
Under the above assumptions, $(A\oplus B, \Delta)$ is isomorphic to the direct sum of $(A,\Delta_{AA}-\Delta_{AB}\circ h\circ \Delta_{BA})$ and the contractible curved complex $(B,\Delta_{BB}+Z\circ h)$.
\end{lemma}
\begin{proof}
Let $\Phi\in \End(A\oplus B)$ be defined by
\[
\Phi=\sqmatrix{\Id_A & -\Delta_{AB}\circ h \\ h\circ \Delta_{BA} & \Id_B}
\]
Then $\Phi$ is an isomorphism, with inverse
\[
\Phi\inv=\sqmatrix{\Id_A & \Delta_{AB}\circ h \\ -h\circ \Delta_{BA} & \Id_B}
\]
We leave it as an exercise to show that
\[
\Phi\circ \Delta \circ \Phi\inv = \sqmatrix{\Delta_{AA}-\Delta_{AB}\circ h\circ \Delta_{BA} & 0\\ 0& \Delta_{BB}+\Delta_{BA}\circ \Delta_{AB}\circ h}
\]
Now, since we assume that $\Delta_{BB}^2=0$ and $\Delta^2=Z_{A\oplus B}$ it follows that $\Delta_{BA}\Delta_{AB}=Z$, as claimed.  Finally, we remark that $(B,\Delta_{BB}+Z\circ h)$ is contractible since
\[
\Delta_{BB}+Z\circ h)\circ h+ h\circ (\Delta_{BB}+Z\circ h) = \Delta_{BB}\circ h+h\circ \Delta_{BB}\circ h +Z\circ h^2+h\circ Z\circ h = \Id_B
\]
by \eqref{eq:Bhomotopy} and \eqref{eq:homotopysquared}, using the fact that $Z$ is central.
\end{proof}

\begin{example}\label{ex:gaussExample}
Consider the cone of the morphism \eqref{eq:crossingchange}, which is illustrated below:
\[
C \  = \ \left(\begin{tikzpicture}[baseline=-4.5em]
\tikzstyle{every node}=[font=\small]
\node (a) at (0,0) {${B_i}[\yy]$};
\node (b) at (4.5,0) {$\underline{R}[\yy](1)$};
\node (c) at (0,-3) {$R[\yy](-1)$};
\node (d) at (4.5,-3) {$\underline{B_i}[\yy]$};
\path[->,>=stealth',shorten >=1pt,auto,node distance=1.8cm,
  thick]
([yshift=3pt] a.east) edge node[above] {$-b\otimes 1$}		([yshift=3pt] b.west)
([yshift=-2pt] b.west) edge node[below] {$b^{\ast}\otimes (y_i-y_{i+1})$}		([yshift=-2pt] a.east)
([yshift=3pt] c.east) edge node[above] {$b^\ast \otimes 1$}		([yshift=3pt] d.west)
([yshift=-2pt] d.west) edge node[below] {$-b \otimes (y_i-y_{i+1})$}		([yshift=-2pt] c.east)
(a) edge node[xshift=1cm,yshift=-.6cm] {$\Id\otimes 1$} (d);
\draw[frontline,->,>=stealth',shorten >=1pt,auto,node distance=1.8cm,thick]
(b) to node[xshift=-3.7cm,yshift=0cm] {$-\Id\otimes (y_i-y_{i+1})$} (c);
\end{tikzpicture}\right)
\]
We let $B$ denote the direct sum of the terms in the Northwest and Southeast corners, and we let $A$ denote the direct sum of the terms in the Northeast and Southwest corners.  Then $(B,\Delta_{BB})$ is the mapping cone on the identity morphism $B_i[\yy]\rightarrow B_i[\yy]$.  This complex is contractible, with null-homotopy given by the identity morphisms in the opposite direction, pictured below:
\[
\left(\begin{tikzpicture}[baseline=-4.5em]
\tikzstyle{every node}=[font=\small]
\node (a) at (0,0) {${B_i}[\yy]$};
\node (b) at (4.5,0) {$\underline{R}[\yy](1)$};
\node (c) at (0,-3) {$R[\yy](-1)$};
\node (d) at (4.5,-3) {$\underline{B_i}[\yy]$};
\path[->,>=stealth',shorten >=1pt,auto,node distance=1.8cm,
  thick]
(d) edge node[xshift=1cm,yshift=-.6cm] {$\Id\otimes 1$} (a);
\end{tikzpicture}\right)
\]
Then $C$ is isomorphic to a direct sum
\[
C \cong \left( \begin{tikzpicture}[baseline=-.5em]
\tikzstyle{every node}=[font=\small]
\node (a) at (0,0) {$R[\yy](-1)$};
\node (b) at (4,0) {$\underline{R}[\yy](1)$};
\path[->,>=stealth',shorten >=1pt,auto,node distance=1.8cm,
  thick]
([yshift=3pt] a.east) edge node[above] {$bb^\ast\otimes 1$}		([yshift=3pt] b.west)
([yshift=-2pt] b.west) edge node[below] {$-\Id\otimes (y_i-y_{i+1})$}		([yshift=-2pt] a.east);
\end{tikzpicture}\right)
 \ \oplus \ 
 \left(
  \begin{tikzpicture}[baseline=-.5em]
\tikzstyle{every node}=[font=\small]
\node (a) at (0,0) {$B_i[\yy]$};
\node (b) at (4,0) {$\underline{B_i}[\yy]$};
\path[->,>=stealth',shorten >=1pt,auto,node distance=1.8cm,
  thick]
([yshift=3pt] a.east) edge node[above] {$\Id\otimes 1$}		([yshift=3pt] b.west)
([yshift=-2pt] b.west) edge node[below] {$-b^\ast b\otimes (y_i-y_{i+1})$}		([yshift=-2pt] a.east);
\end{tikzpicture}
\right).
\]
Observe that the second summand is contractible.
\end{example}

\printbibliography

\end{document}